\documentclass[11pt]{amsart}

\usepackage{graphicx,hyperref}
\usepackage[T1]{fontenc}
\usepackage[cp1250]{inputenc}

\usepackage{amsaddr}

\title[Analytic sphere eversion]{Analytic sphere eversion using ruled surfaces}
\author{Adam Bednorz, Witold Bednorz}
\address{A. Bednorz: Faculty of Physics, University of Warsaw, ul. Pasteura 5, 02-093 Warsaw, Poland}
\address{W. Bednorz: Faculty of Mathematics, Informatics and Mechanics, University of Warsaw, ul. Banacha 2, 02-093 Warsaw, Poland}
\email{abednorz@fuw.edu.pl, wbednorz@mimuw.edu.pl}

\keywords{sphere eversion, Boy surface, ruled surfaces}

\begin{document}
\begin{abstract}
Sphere eversions have been described so far by either  pictures with minimal topological complexity, numerical evolution or complex equations.
We write down relatively simple explicit formulas for the whole eversion, both analytic and topologically simpler, including also Boy surface (real projective plane), using a family of ruled surfaces. We show their usefulness in visualizing the process using commonly available modeling software.
\end{abstract}
\maketitle

\section{Introduction}

Over 50 years ago Smale proved that a sphere ($S^2$ in space $\mathbb R^3$) can be everted in continuous way \cite{smale}. More precisely the set of all sphere immersions,
smooth functions $\vec{R}\to \vec{r}(\vec{R})$  for $\vec{R}\in S^2$ is connected. Since $\vec{r}(\vec{R})=-\vec{R}$ is an immersion, it means that the sphere can be 
continuously turned inside out, without crease, although allowing for self-intersections, with some continuous $\vec{r}(\vec{R},t)$ such that $\vec{r}(\vec{R},t_-)=\vec{R}$ and $\vec{r}(\vec{R},t_+)=-\vec{R}$ ($t$ can be thought of as time). Unfortunately, the Smale's proof gives little 
hint how to visualize the process. Only later detailed models of eversion \cite{phil,morin1,morin2,morin3,morin4,morin5,mfilm, topol}, including discussion of critical points and halfway models \cite{morin2,morin3,apery}, appeared. Computer era offered new tools of presentation \cite{levy,sull}, like numerical evolving from a halfway model \cite{fsk}. The list of very clever approaches to sphere eversion is left open \cite{cherit,atch,neve}. 

Both numerical and pictorial approaches lack full analyticity, achieved only in the original proof \cite{smale} and Morin model \cite{morin1} (and halfways \cite{apery,sull,fsk}). Morin model fails to keep the minimum number of topological events, achieved in contrast in numerical and pictorial models \cite{morin2,morin3,apery,fsk} but without analytic proof.
Here we close this gap, by presenting a set of formulas to describe the complete eversion process in the topologically simplest way, based on a family of ruled surfaces.
The formulas are maximally simplified and shown to preserve expectations from eversion: smoothness and the simplest set of critical topological events. They are also useful in a modeling of the eversion by commonly available computer software (here we used Mathematica).

Our work is organized as follows. We show eversion of a twisted cylinder -- topological annulus (not yet a sphere), discussing all relevant features, with the critical halfway surface, and generalization to nonorientable Boy surface, using a family of ruled surfaces. Next, we transform the cylinder into a sphere by a kind of inverse stereographic projection. The process is heavily illustrated with Mathematica pictures and the rigorous proofs of smoothness and other properties are left in Appendices.

\section{Cylinder eversion}

\begin{figure}
\includegraphics[scale=0.3]{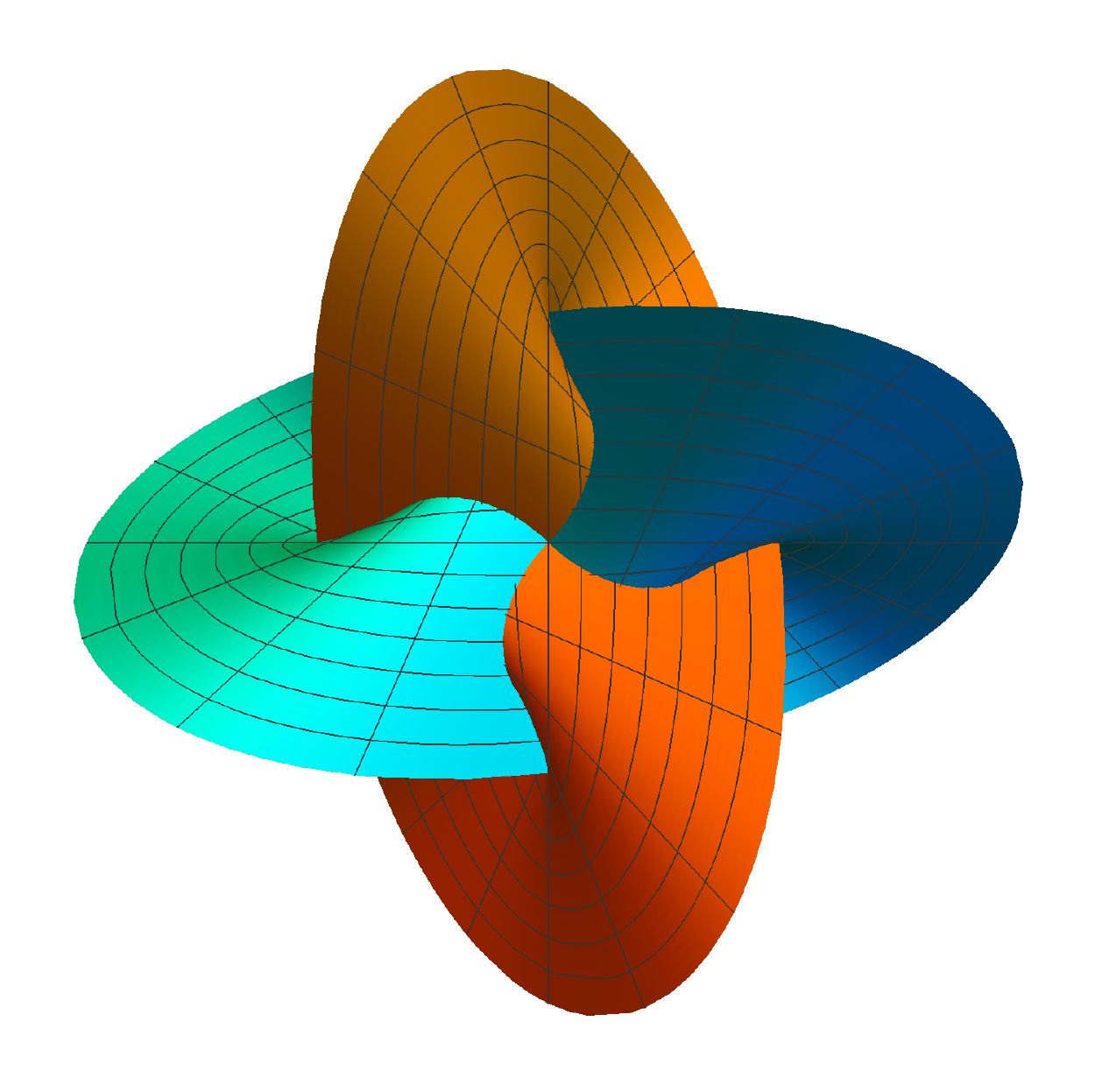}
\includegraphics[scale=0.3]{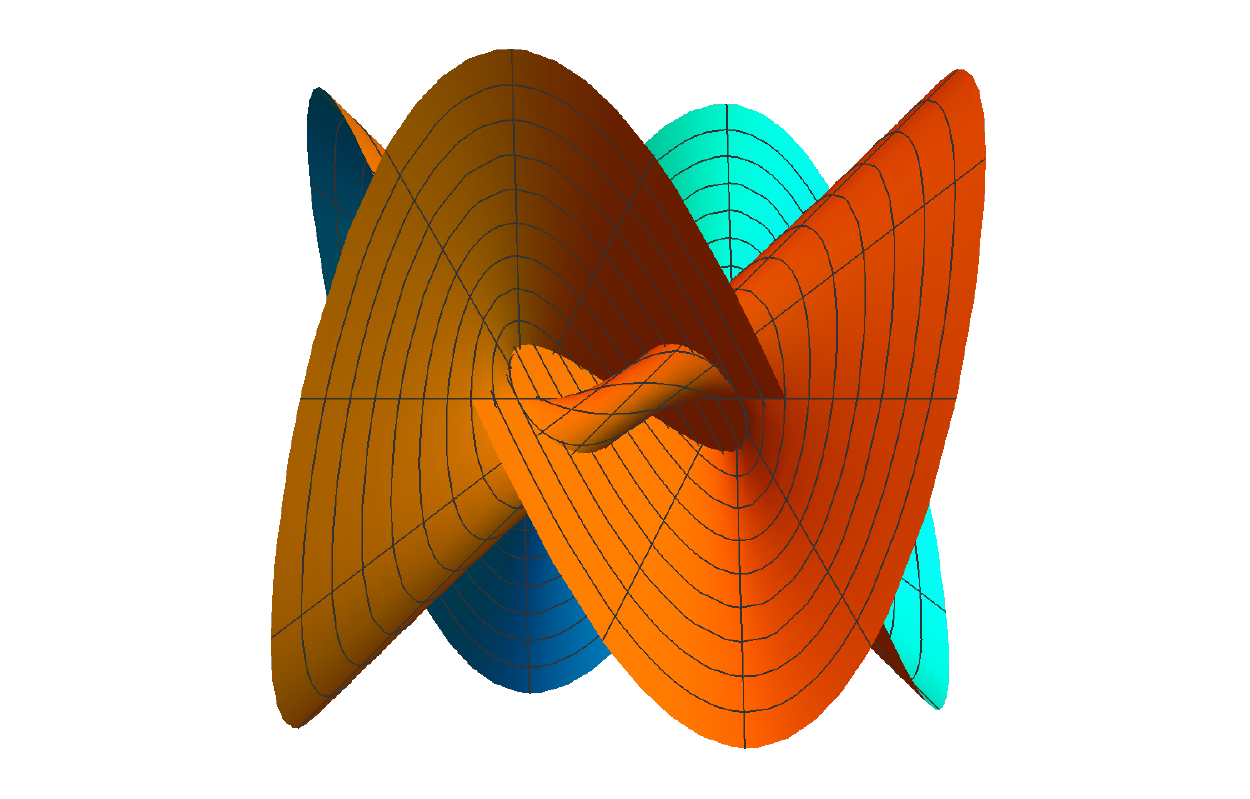}
\includegraphics[scale=0.3]{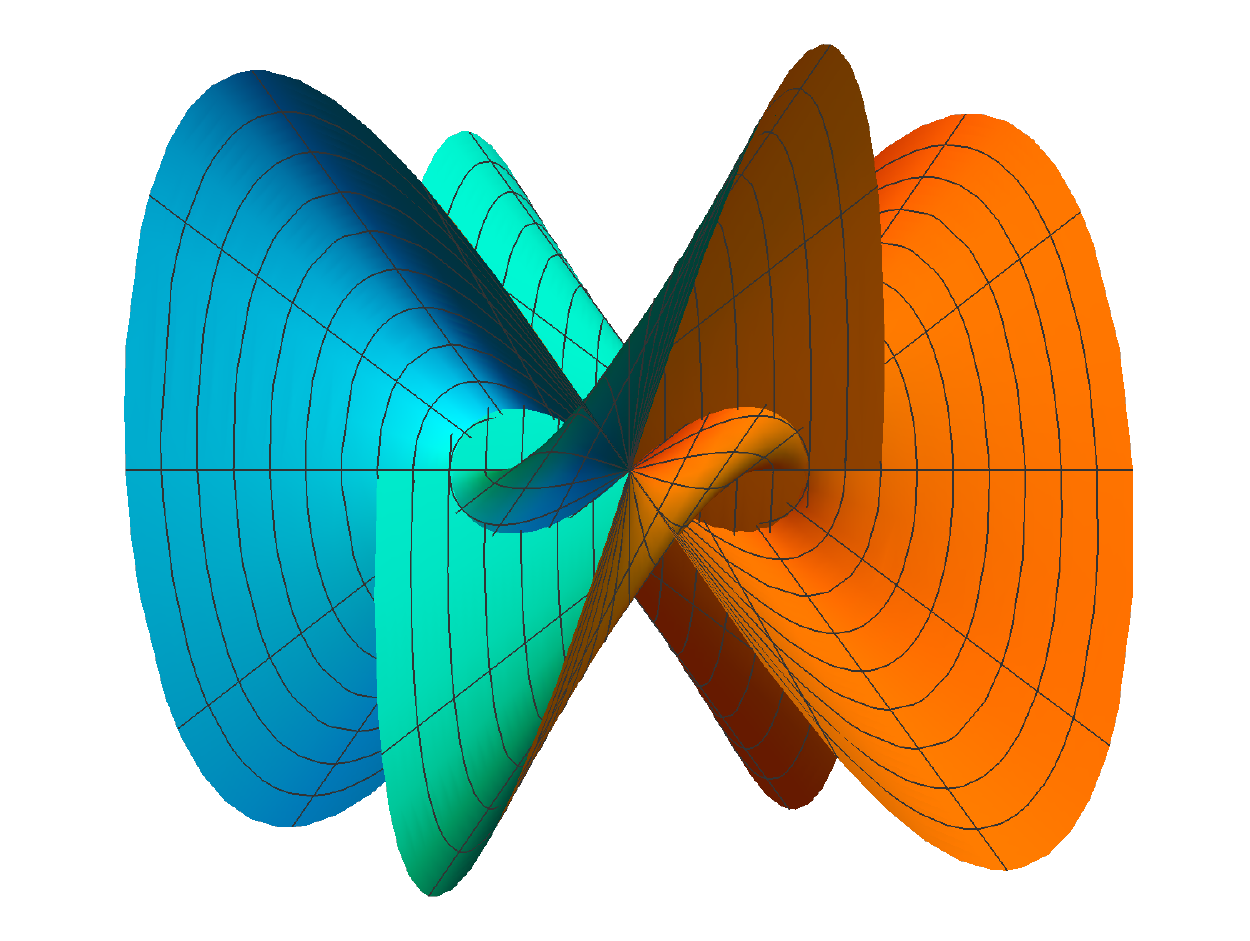}
\caption{The halfway $n=2$ surface (\ref{halfn}) viewed from different angles with different colors at opposite sides. Left shows the central $Q$ point, middle shows pairs of $D_1$ points, right shows a pair of $D_1$ point while an additional $D_1$, $Q$ and the last $D_1$ lie on the central line of perpendicular view}
\label{hh}
\end{figure}

Before eversion of  the complete sphere, we cut out the poles of the sphere and consider a cylinder, which we extend to infinity, similarly as done by Morin \cite{morin1}.  We set the workspace as $(\vec{r},t)\in \mathbb R^3\times \mathbb R$ with $\vec{r}=(x,y,z)$ and time $t$.
Let us take a very special immersion of such a cylinder, $(h,\phi)\in \mathbb R\times S^1$ (here $S^1$ is parametrized by a real variable of the period $2\pi$). The immersion is given by
\begin{equation}
x=\sin(n-1)\phi - h\sin\phi,\:y=\cos(n-1)\phi+h\cos\phi,\:z=h\sin n\phi\label{halfn}
\end{equation}
for a natural $n\geq 2$. The surface is always smooth, see Appendix A. For odd $n$ in (\ref{halfn}), $h\to -h$, $\phi\to \phi+\pi$ gives the same point but opposite oriented.

 The case $n=2$ is depicted in Fig. \ref{hh}. It is quite clear that it there is no privileged side of the surface, making it the best candidate for a halfway model, in the middle of the eversion. In Appendix B we show that it is sextic (degree 6).
The surface contains also critical topological events: point $Q=(0,0,0)$ for $(h,\phi)=(1,\pi/2),(1,-\pi/2),(-1,0),(-1,\pi)$ (four-fold intersection), inevitable in sphere eversion \cite{bmax,hugh}, and four $D_1$ points (crossing of saddles), see details in Appendix C. 

The $D_1$ points are located at: $(\pm \sqrt{2},0,0)$ and $(0,\pm\sqrt{2},0)$ while
the lines of self-intersections are located at two straight lines $x=z=0$ and $y=z=0$ and the curve
\begin{equation}
x=\sqrt{2}\cos 2\varphi\cos\varphi,\:y=\sqrt{2}\cos 2\varphi\sin\varphi,\:z=-(1/2)\sin 4\varphi
\end{equation}
(with $\varphi$ covering interval of length $4\pi$), depicted in Fig. \ref{qf}, see Appendix D.
The projection of this curve onto $xy$ plane is known as quadrifolium or four-leaved clover (rose curve of order 2).

\begin{figure}
\includegraphics[scale=0.3]{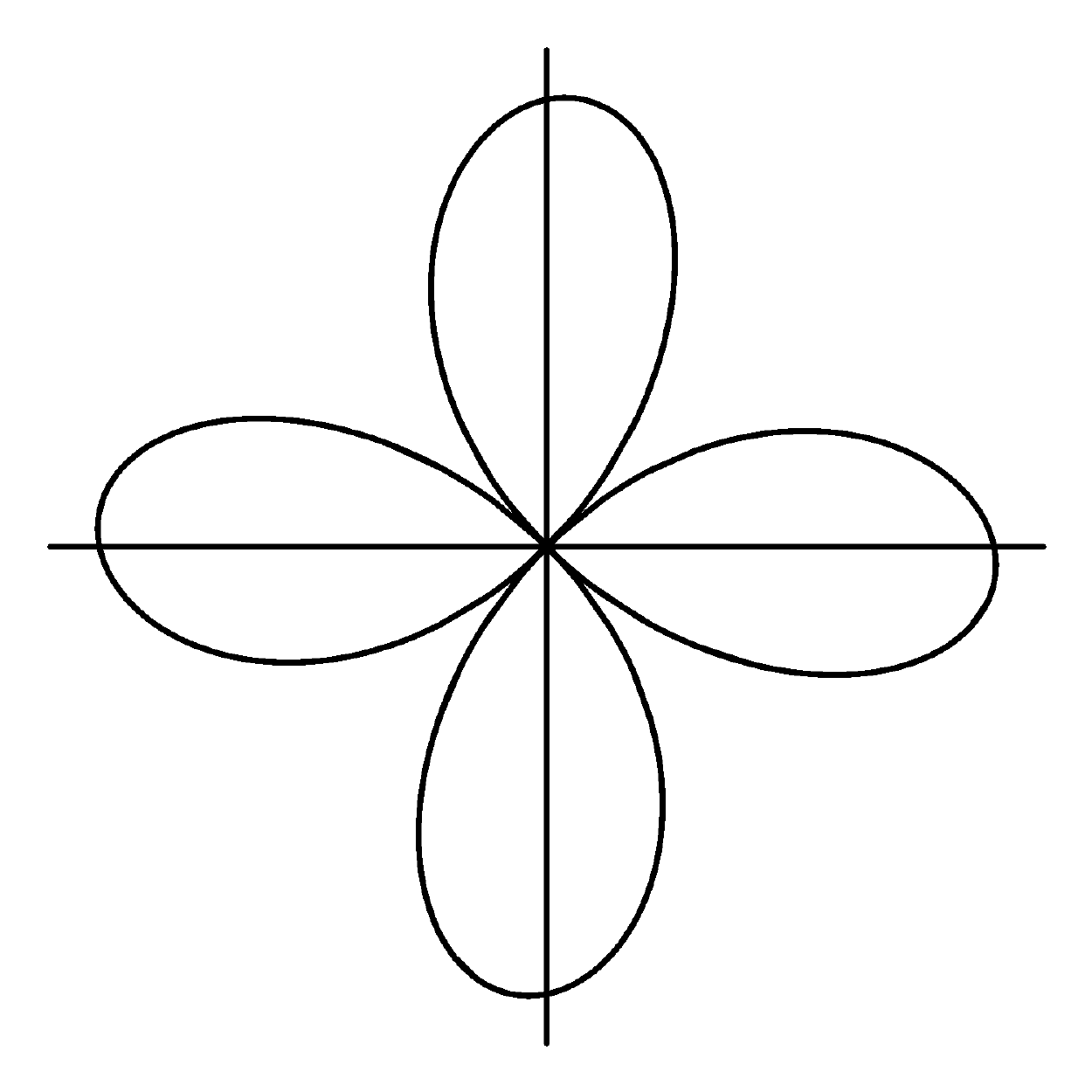}
\includegraphics[scale=0.3]{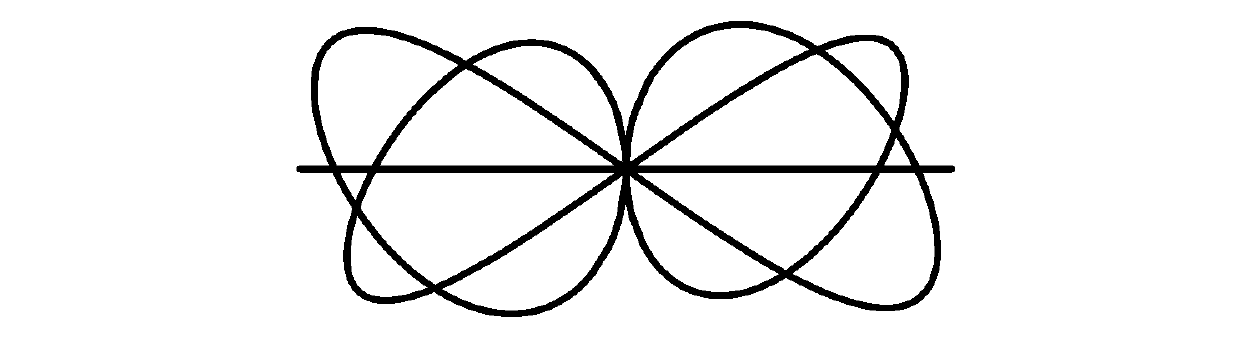}
\includegraphics[scale=0.3]{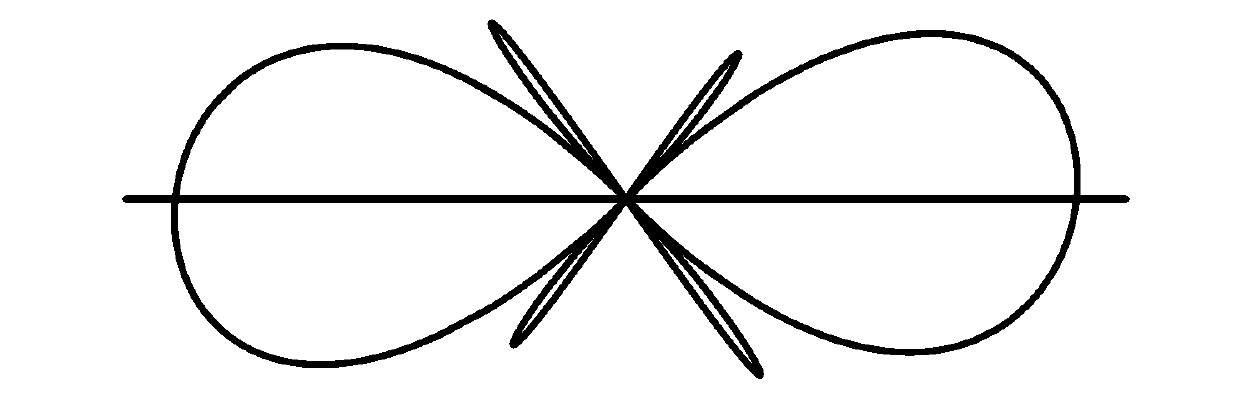}
\caption{The lines of self-intersections of the halfway $n=2$ surface (\ref{halfn}) viewed from the same angles as in Fig. \ref{hh}.}
\label{qf}
\end{figure}

\begin{figure}
\includegraphics[scale=0.3]{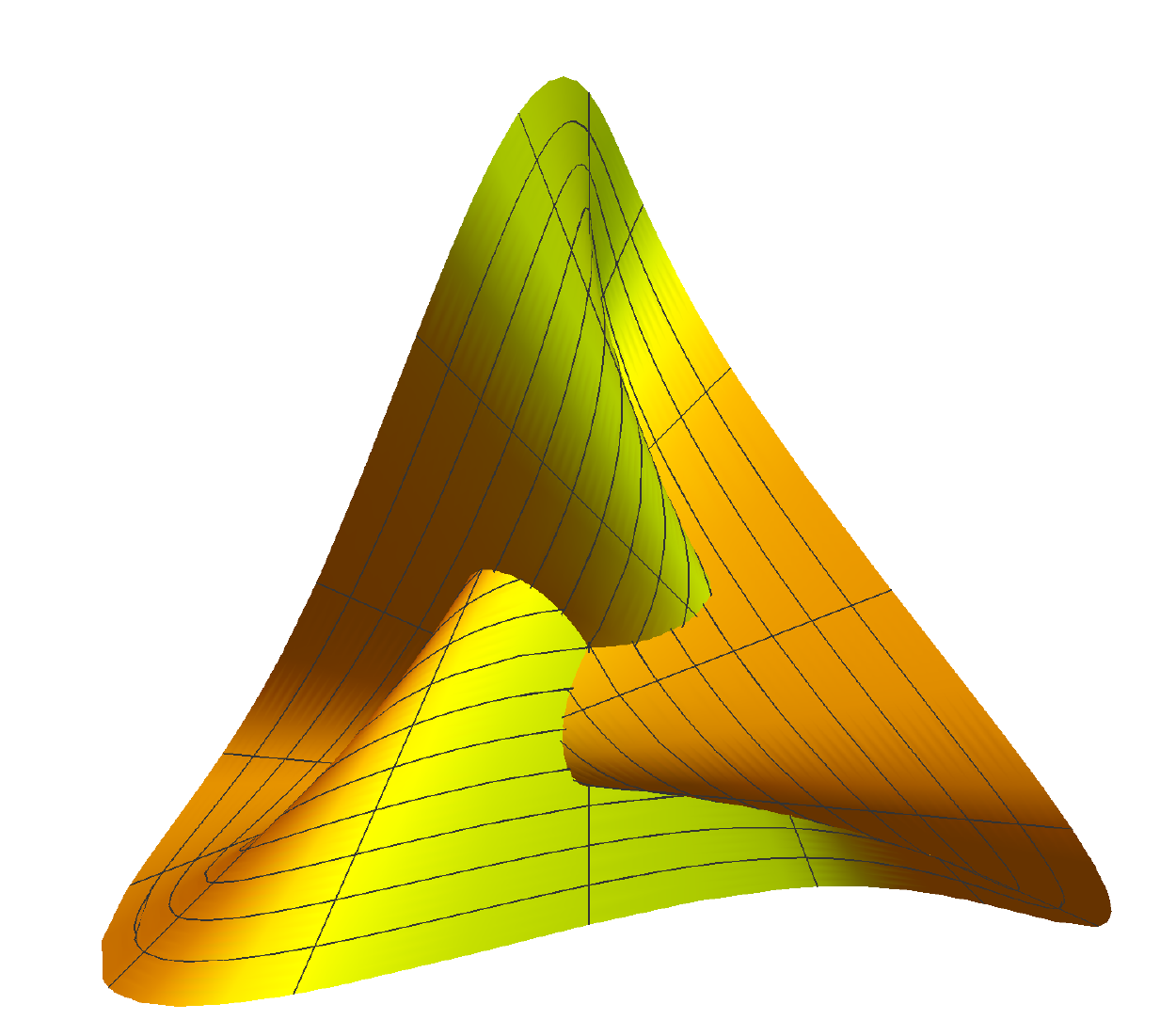}
\includegraphics[scale=0.3]{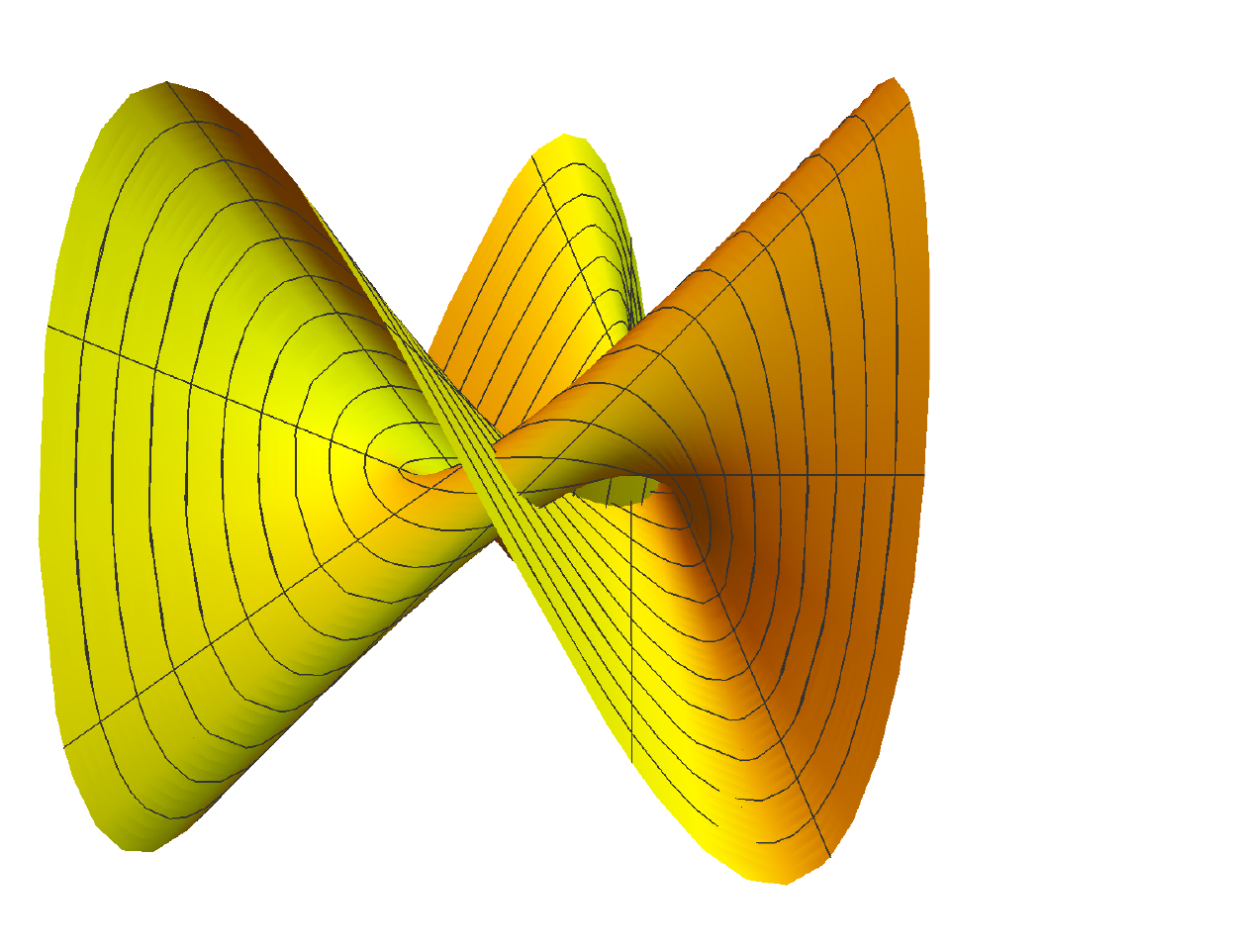}
\includegraphics[scale=0.3]{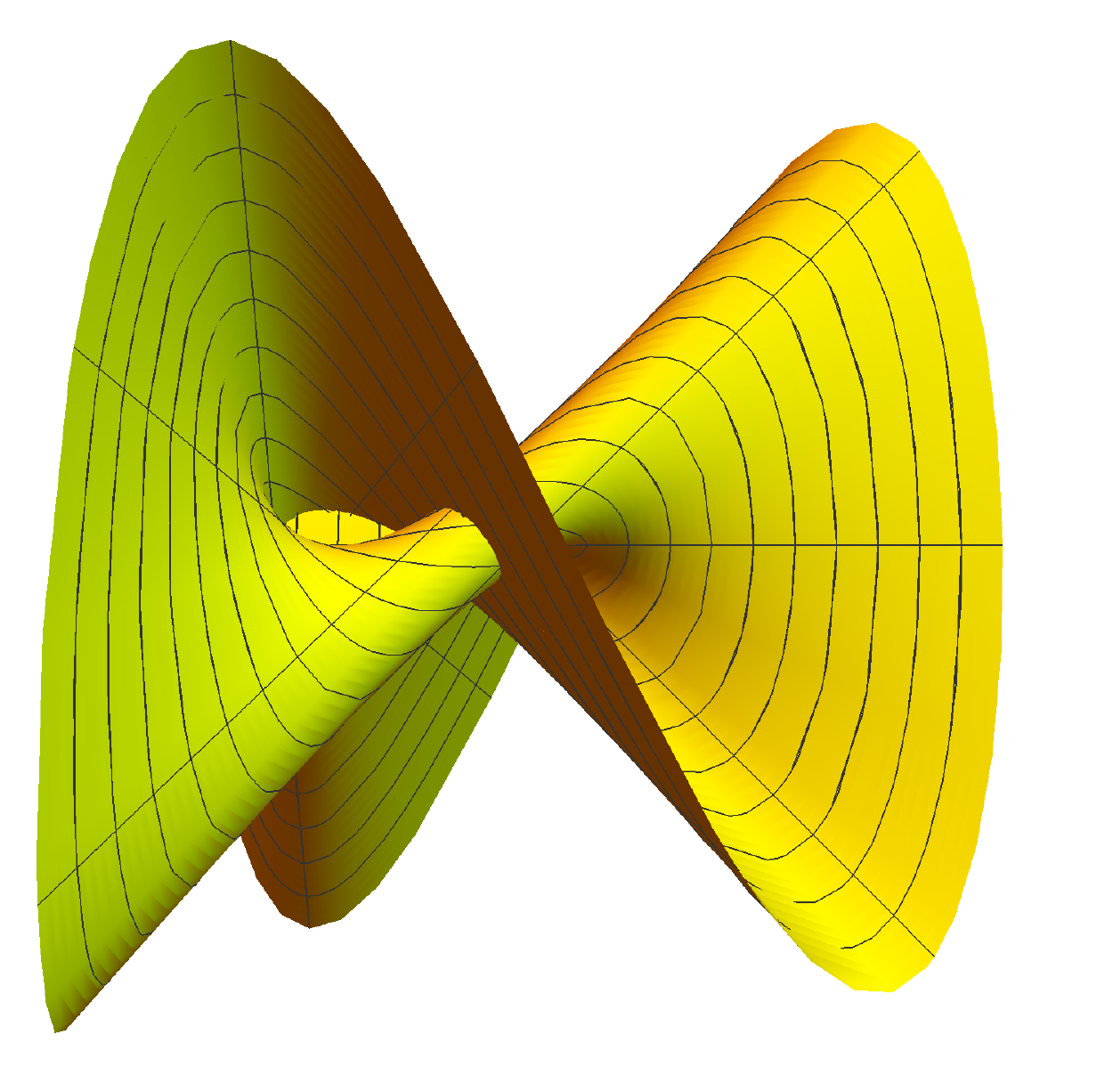}
\caption{The Boy surface (\ref{halfn}) for $n=3$  viewed as in Fig. \ref{hh}.}
\label{boy}
\end{figure}

The case $n=3$ of (\ref{halfn}) gives Boy surface (can be easily closed at infinity),
see Fig. \ref{boy}. The Boy surface is an smooth immersion of real projection plane ($S^2$ with $\vec{R}\equiv -\vec{R}$) in $\mathbb R^3$, also parameterized by Morin \cite{morin1} with $n=3$. It is not orientable, has a single three-fold intersection at $(0,0,0)$ and trifolium self-intersection,
see Fig. \ref{boyi} and Appendix D.
The (Boy) surface at $n=3$ is quintic (degree 5), see Appendix E.

\begin{figure}
\includegraphics[scale=0.3]{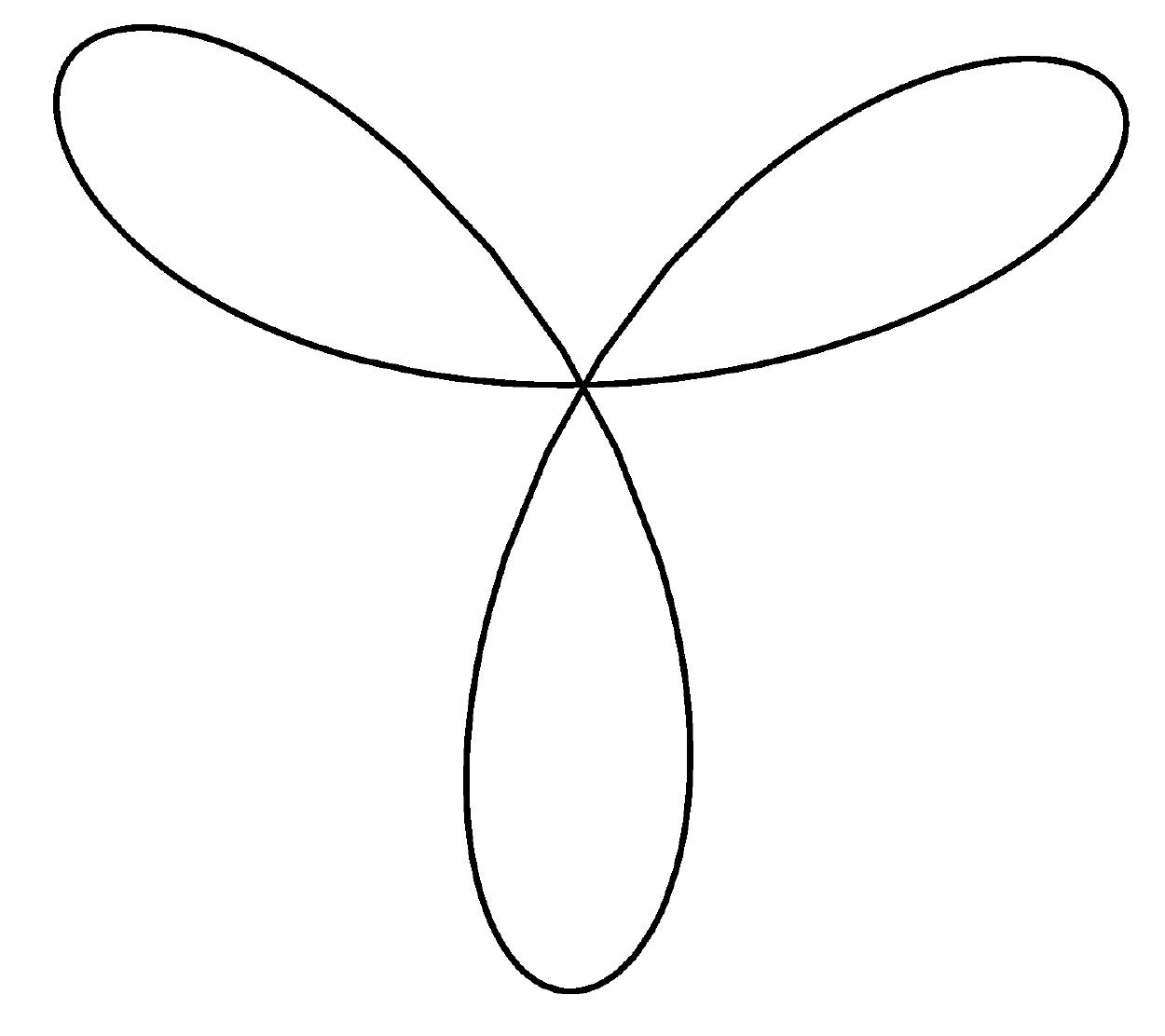}
\includegraphics[scale=0.3]{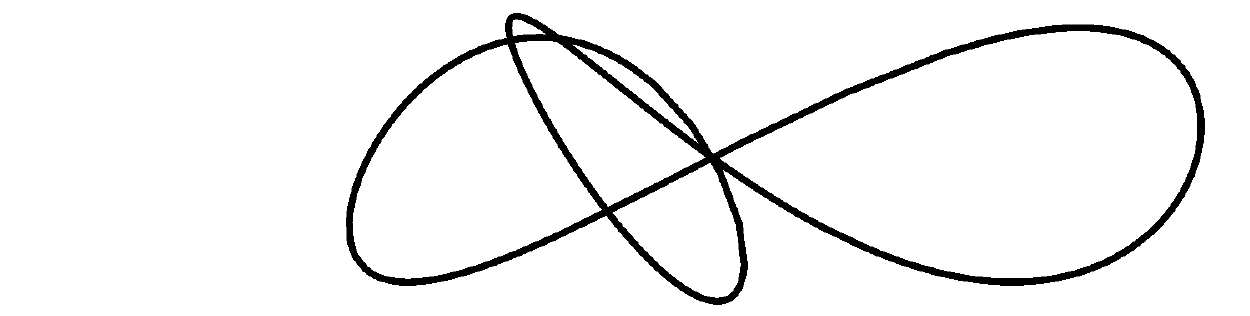}
\includegraphics[scale=0.3]{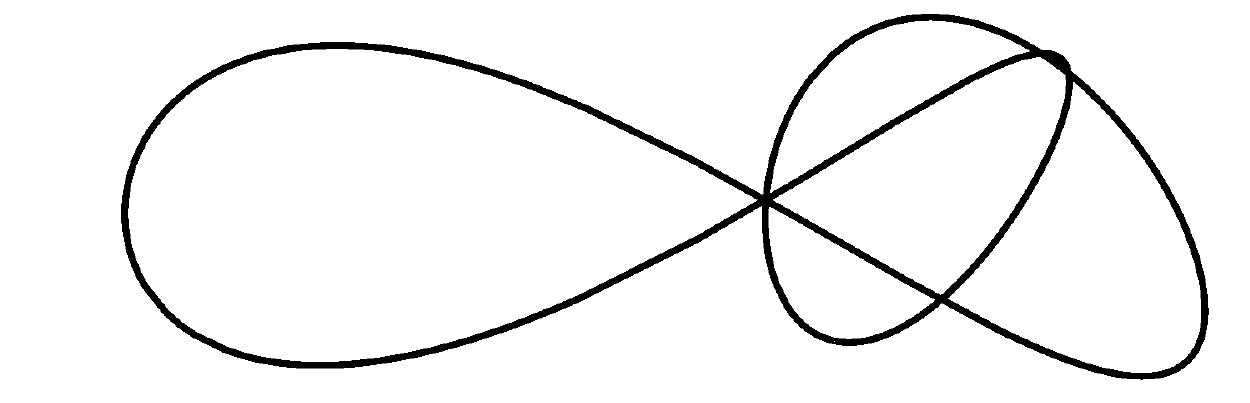}
\caption{The lines of self-intersections of the Boy surface (\ref{halfn}) for $n=3$ viewed from the same angles as in Fig. \ref{hh}.}
\label{boyi}
\end{figure}

Now we generalize the surfaces into  time-dependent ones
\begin{equation}
\begin{matrix}
x=t\cos\phi+\sin (n-1)\phi-h\sin\phi\\
y=t\sin \phi+\cos(n-1)\phi+h\cos\phi\\
z=h\sin n\phi-(t/n)\cos n\phi\end{matrix}\label{twormn}
\end{equation} 
For each $n\geq 2$ one can use (\ref{twormn})
to perform the central step of sphere eversion. This is true even in the case $n=3$ where $t=0$ means two overlapping Boy surfaces.
The case $n=2$ is still sextic, see Appendix B.

\begin{figure}
\includegraphics[scale=0.3]{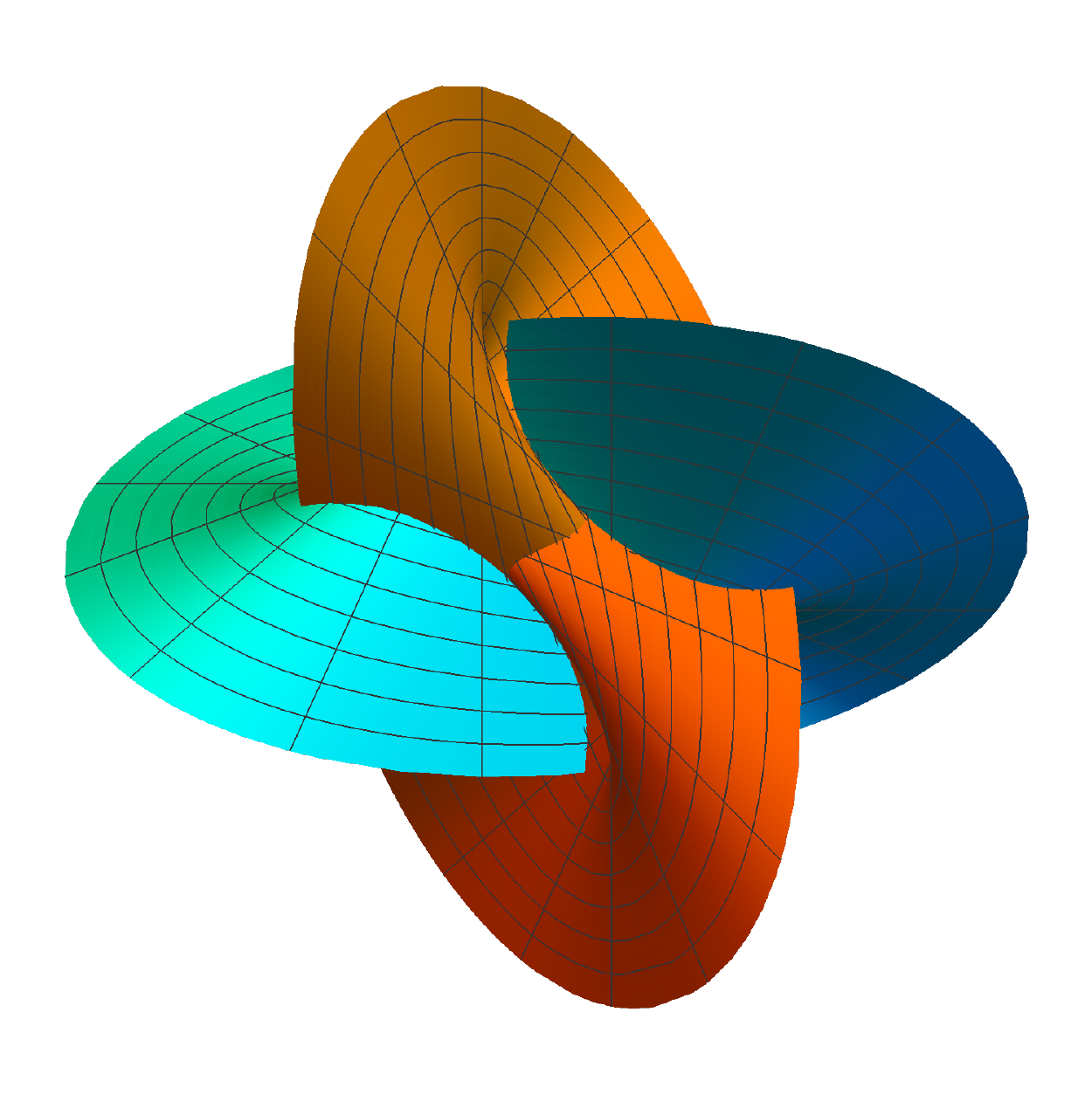}
\includegraphics[scale=0.3]{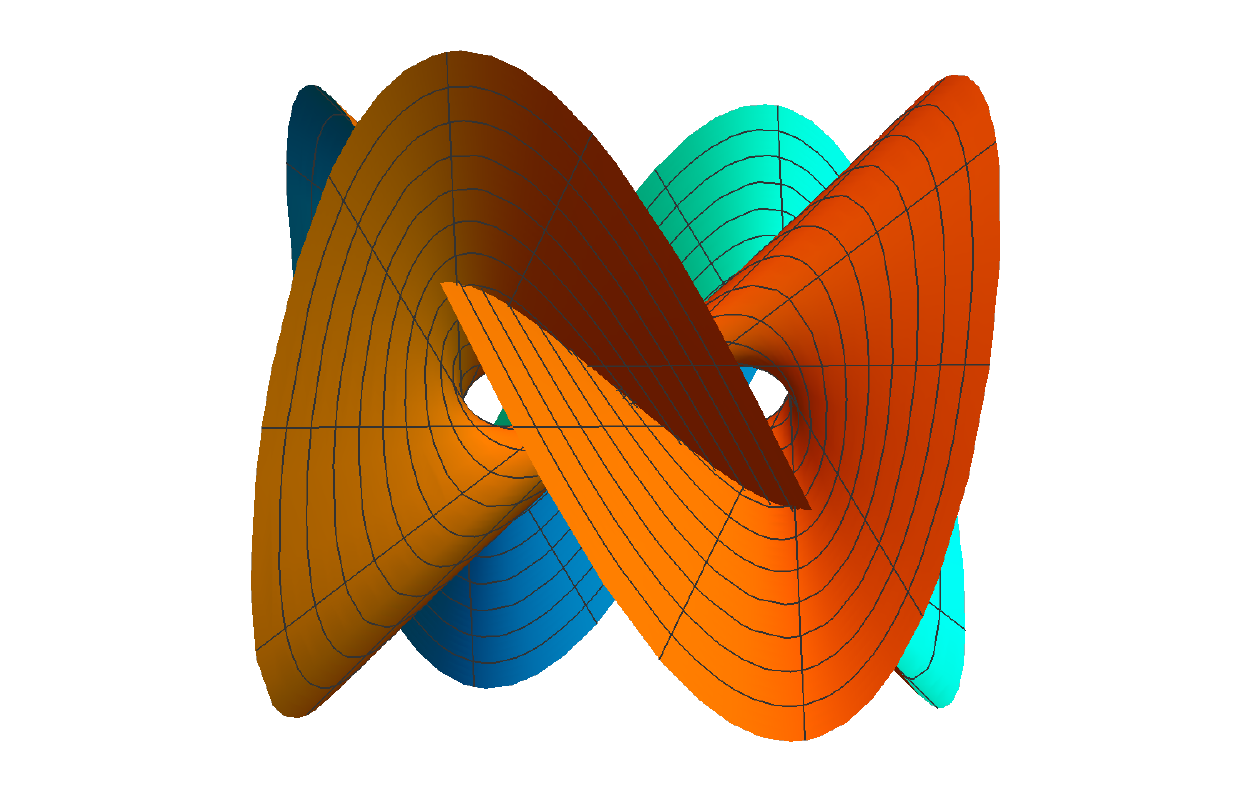}
\includegraphics[scale=0.3]{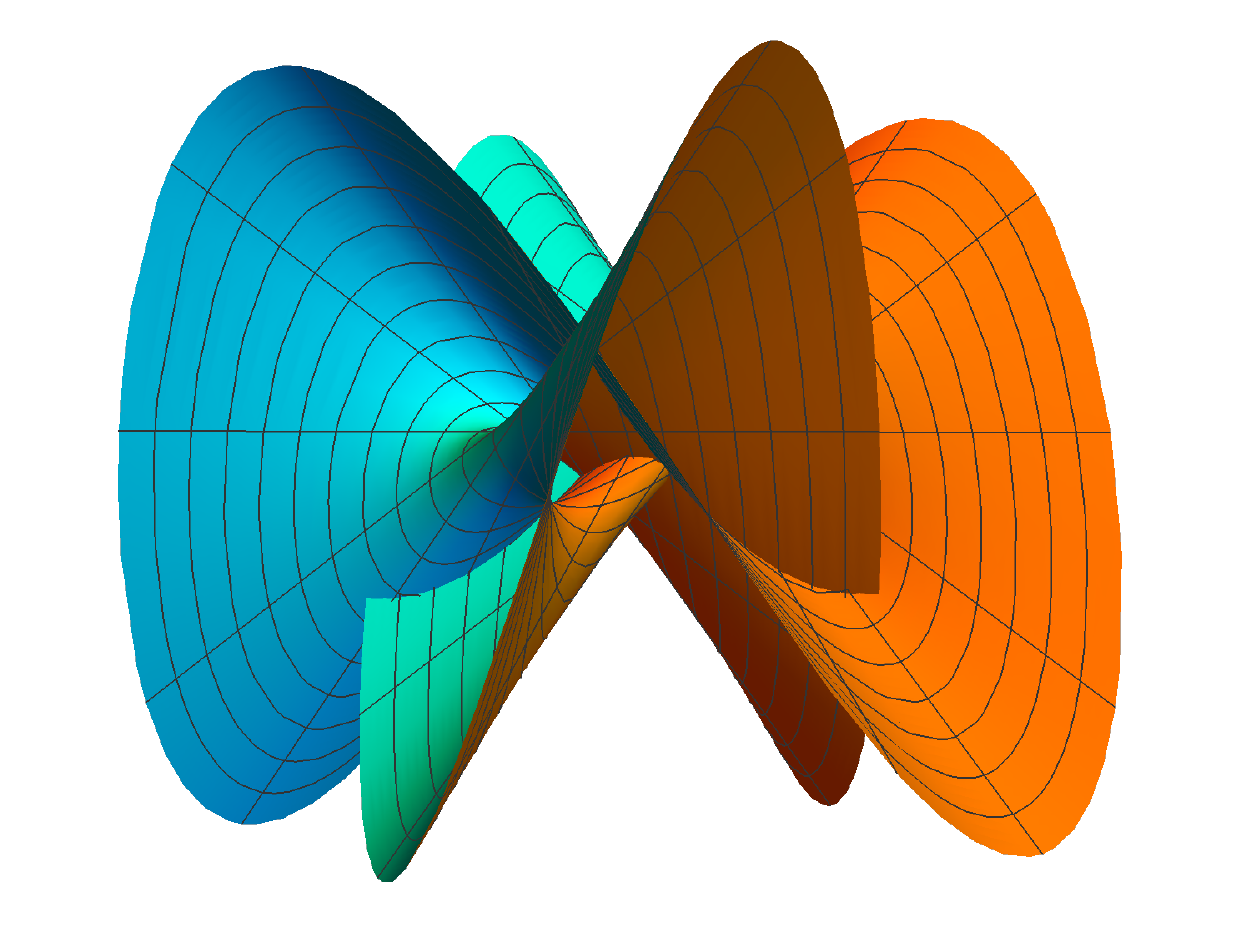}
\includegraphics[scale=0.3]{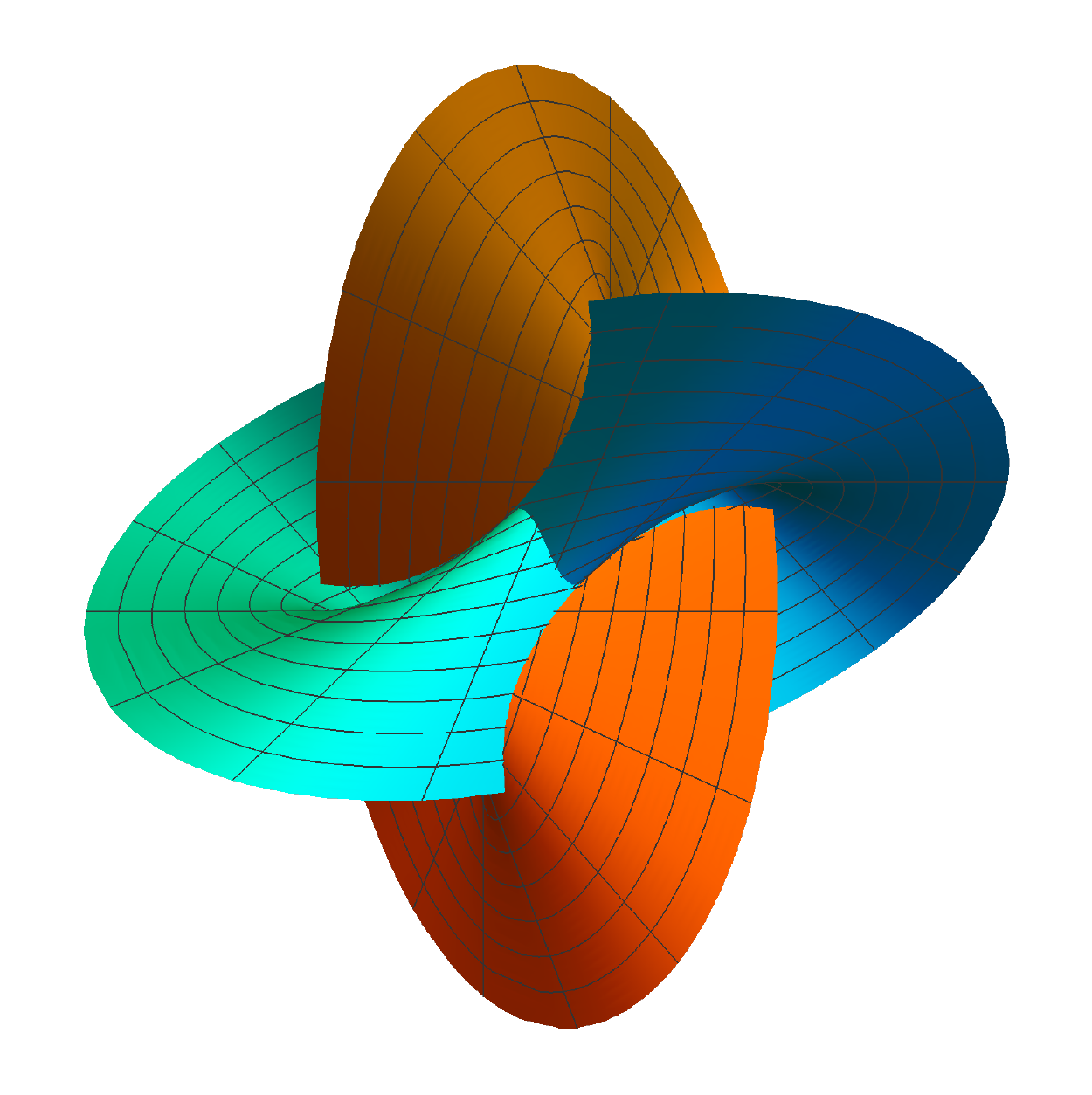}
\includegraphics[scale=0.3]{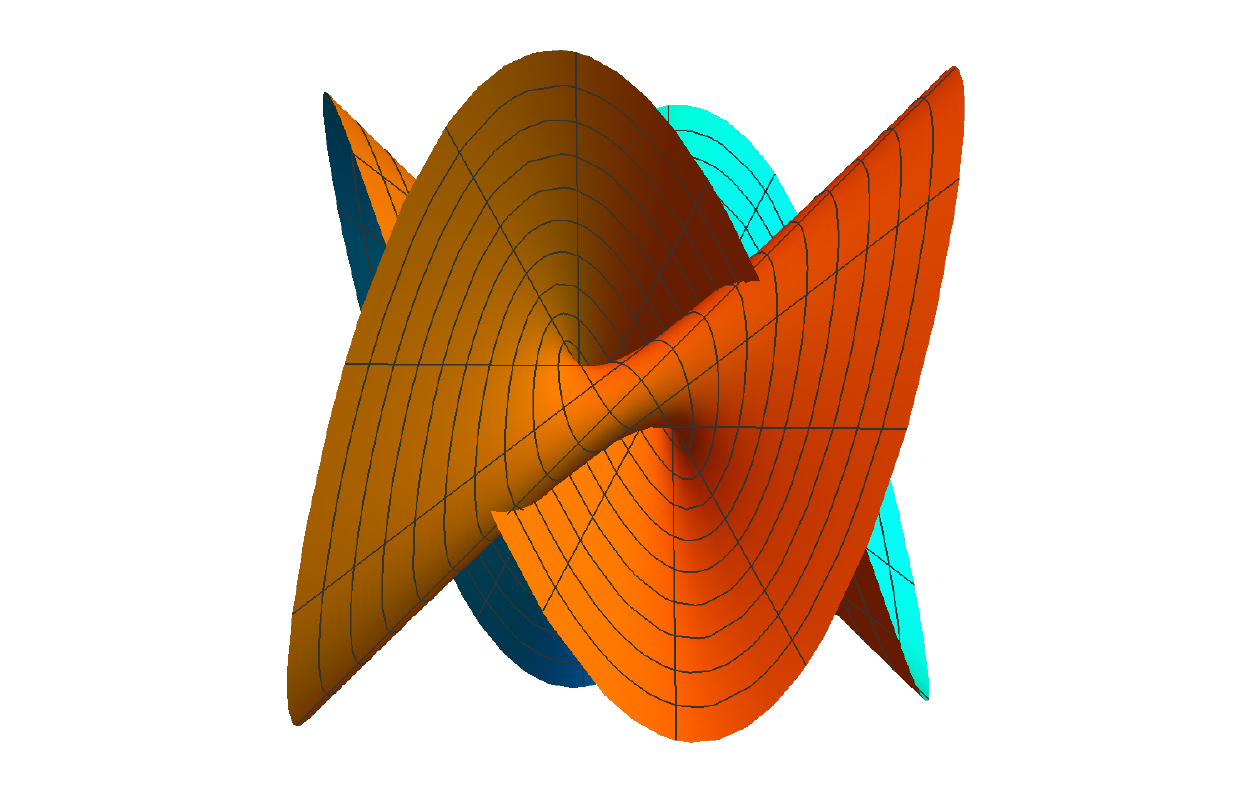}
\includegraphics[scale=0.3]{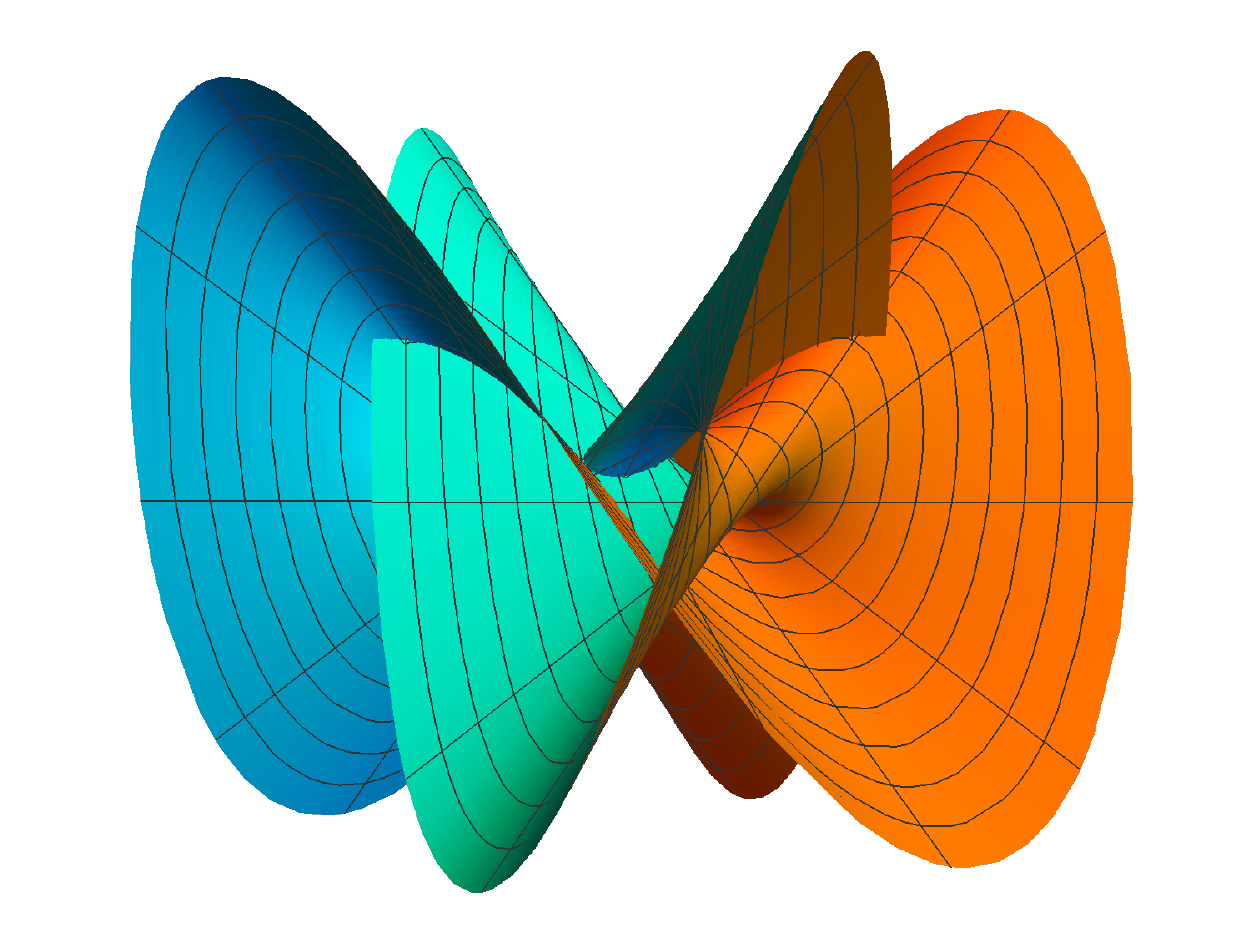}
\caption{The $n=2$ surface (\ref{twormn}) at pair of $T_+$ (top) and $T_-$ (bottom), $t=\mp(\sqrt{17}-3)/2$ viewed from the same angles as in Fig. \ref{hh}. The $T_\pm$ points are best visible on the left.}
\label{tt}
\end{figure}

\begin{figure}
\includegraphics[scale=0.3]{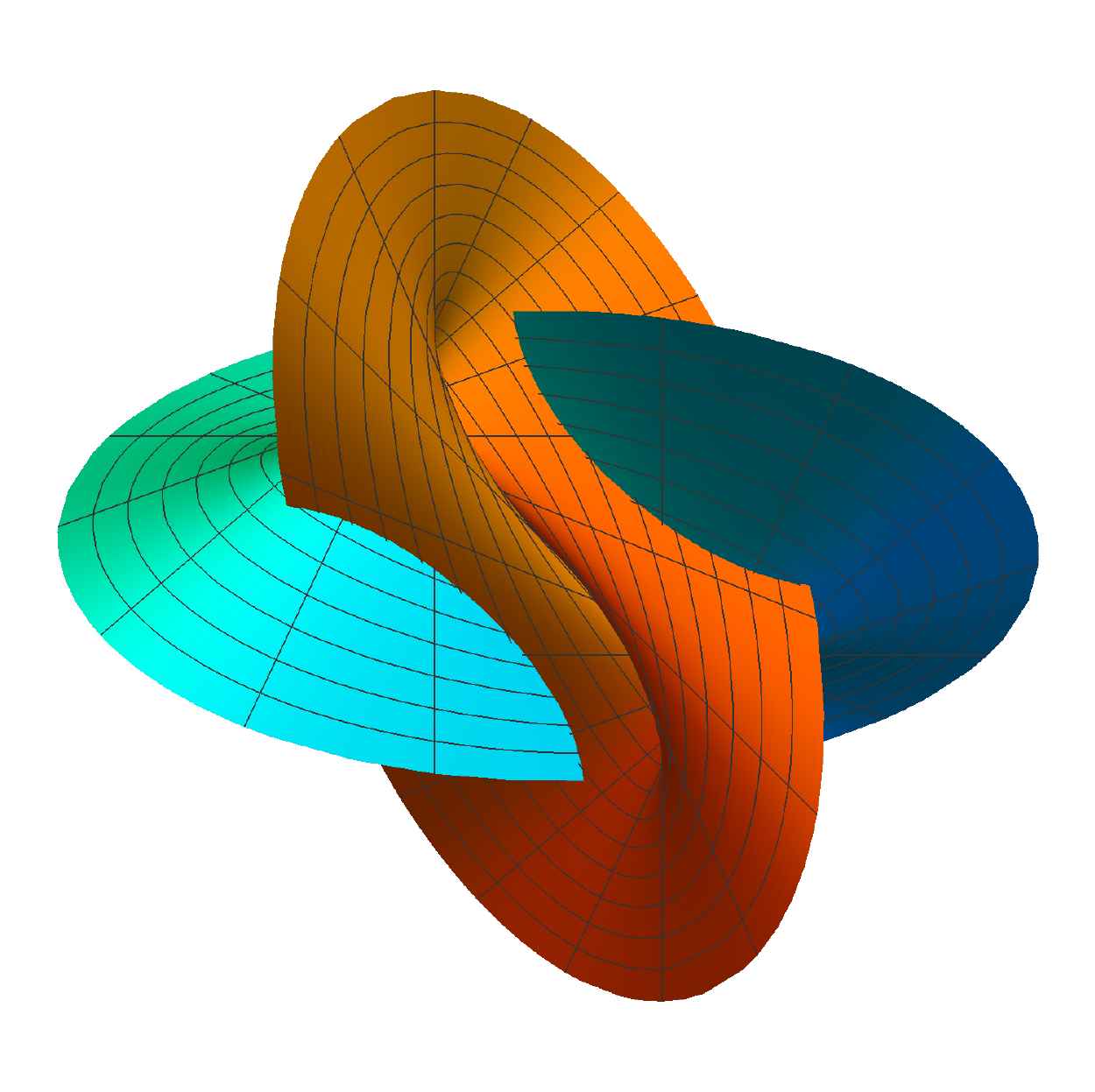}
\includegraphics[scale=0.3]{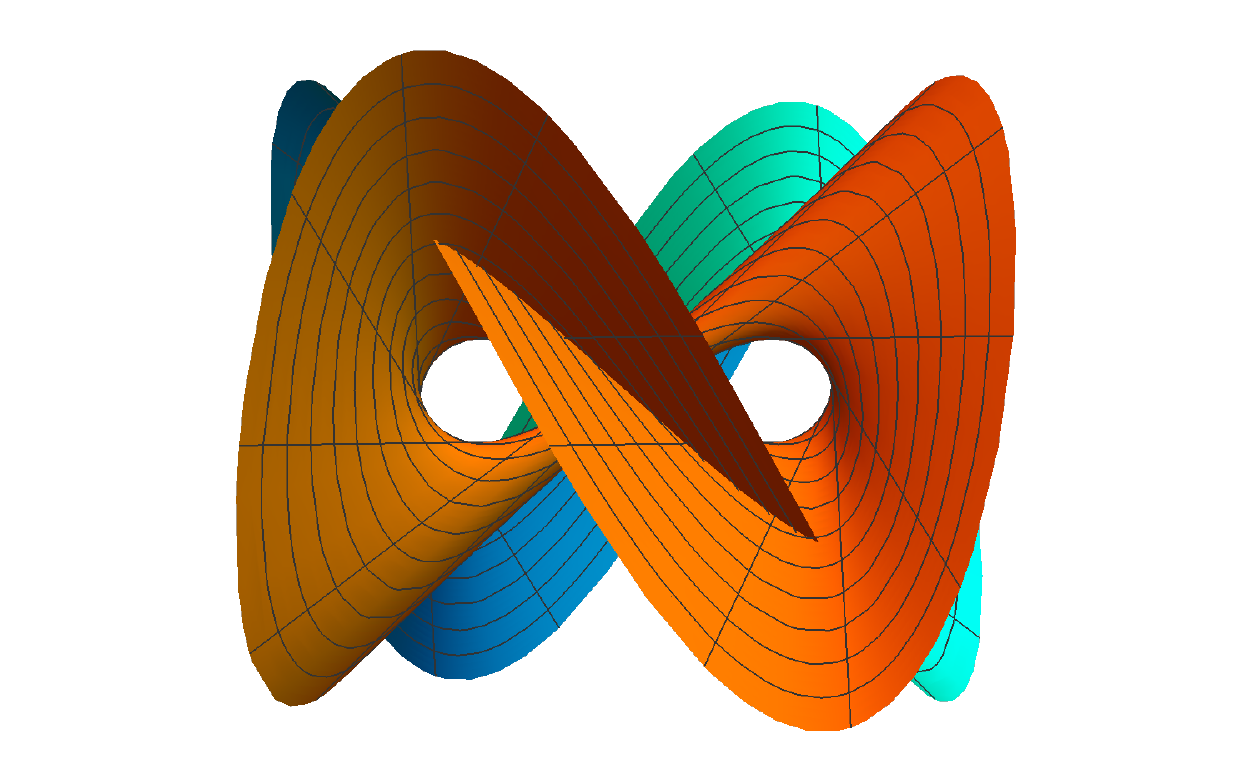}
\includegraphics[scale=0.3]{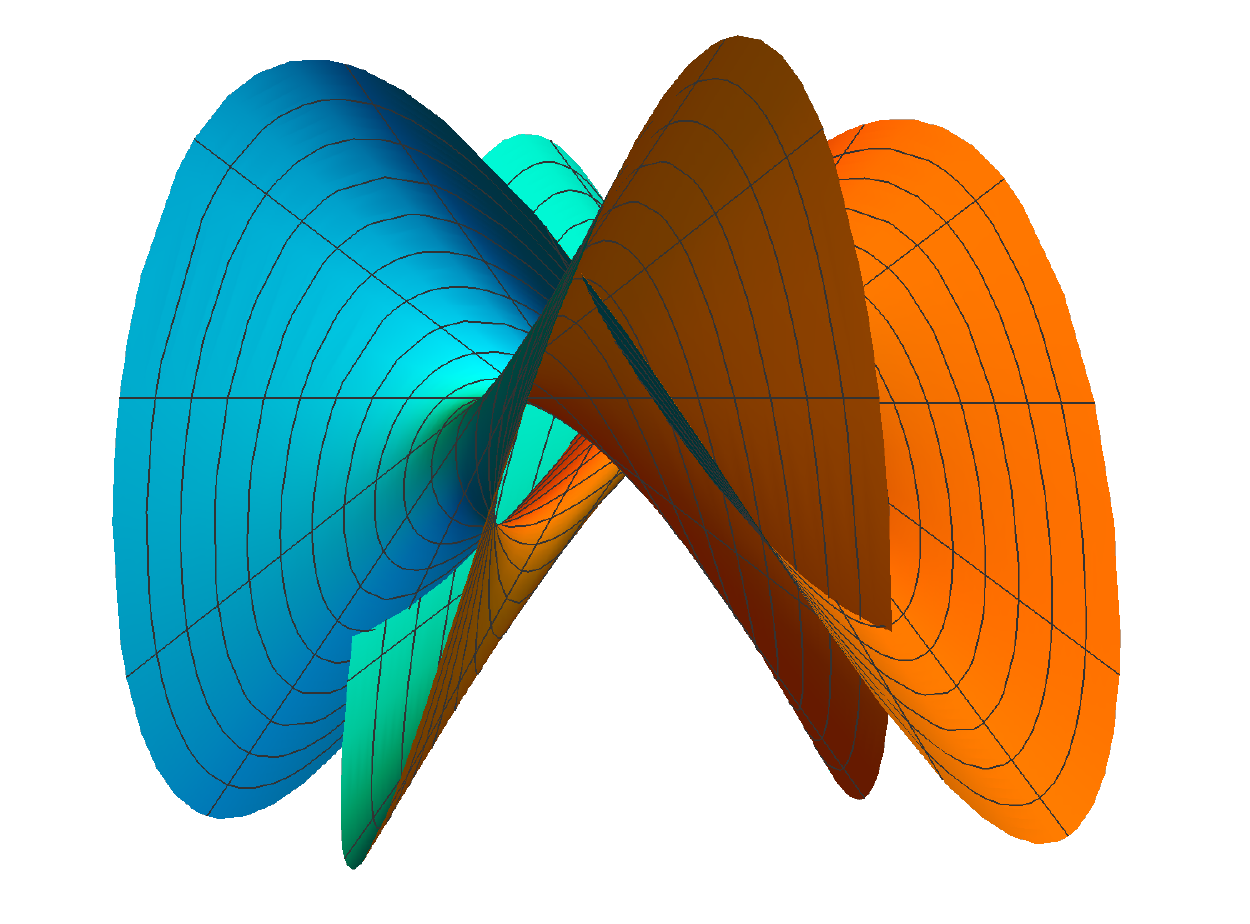}
\includegraphics[scale=0.3]{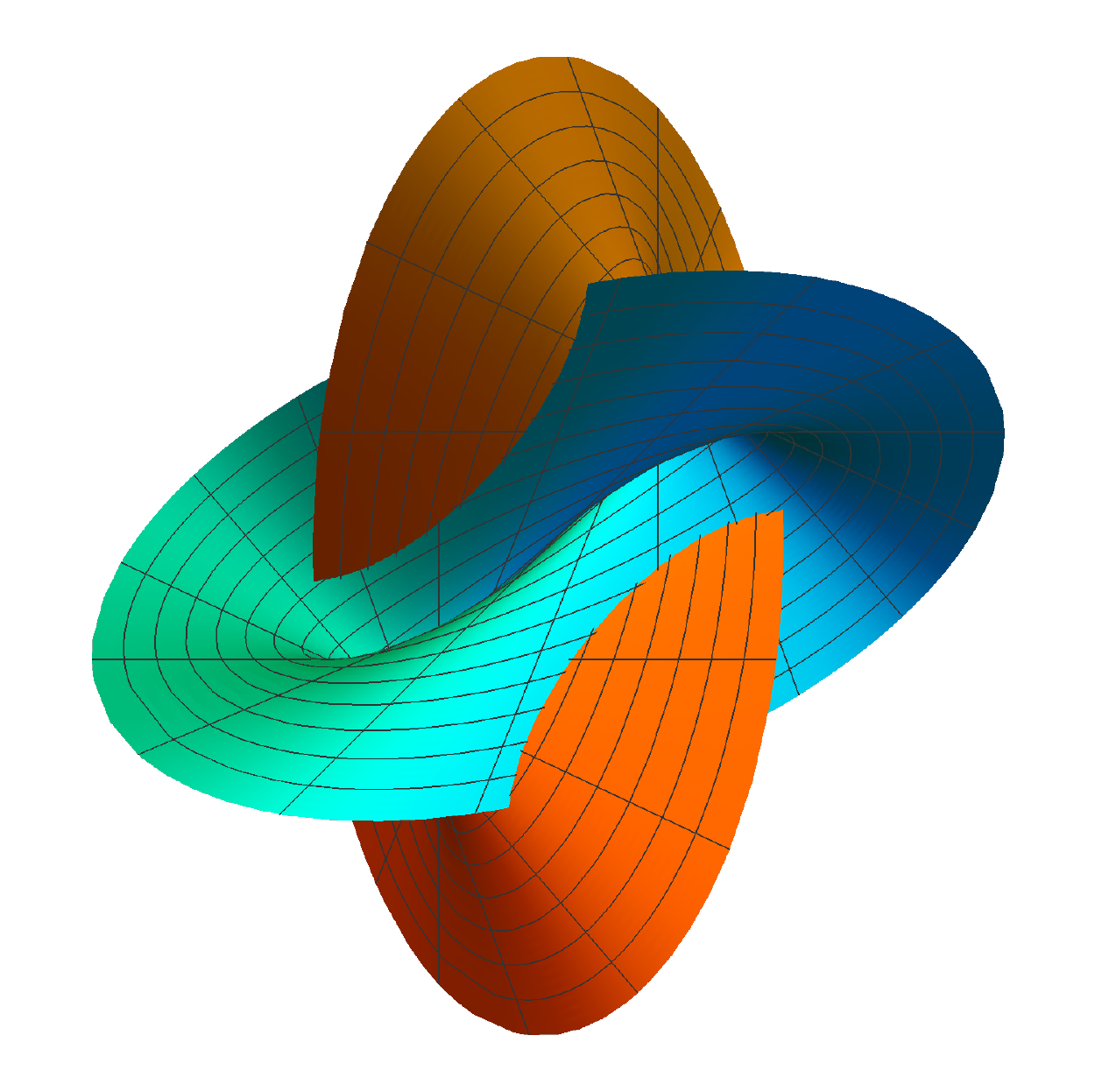}
\includegraphics[scale=0.3]{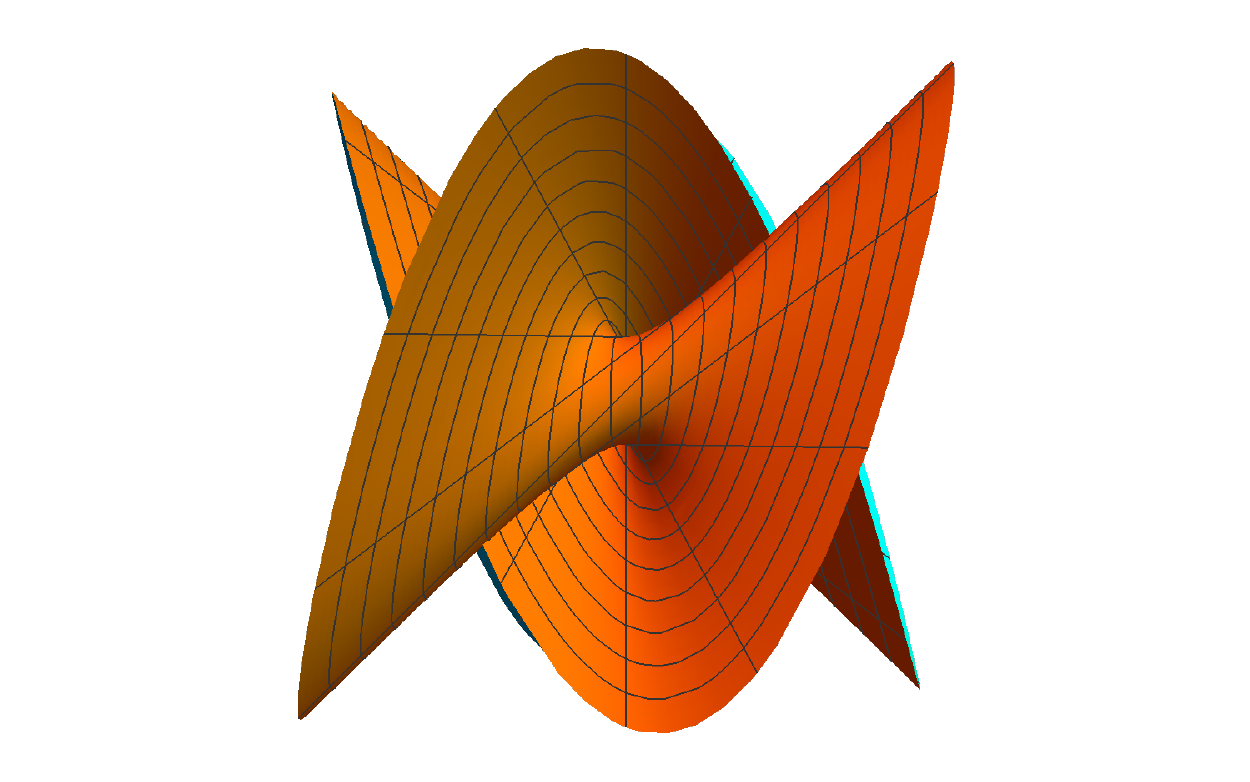}
\includegraphics[scale=0.3]{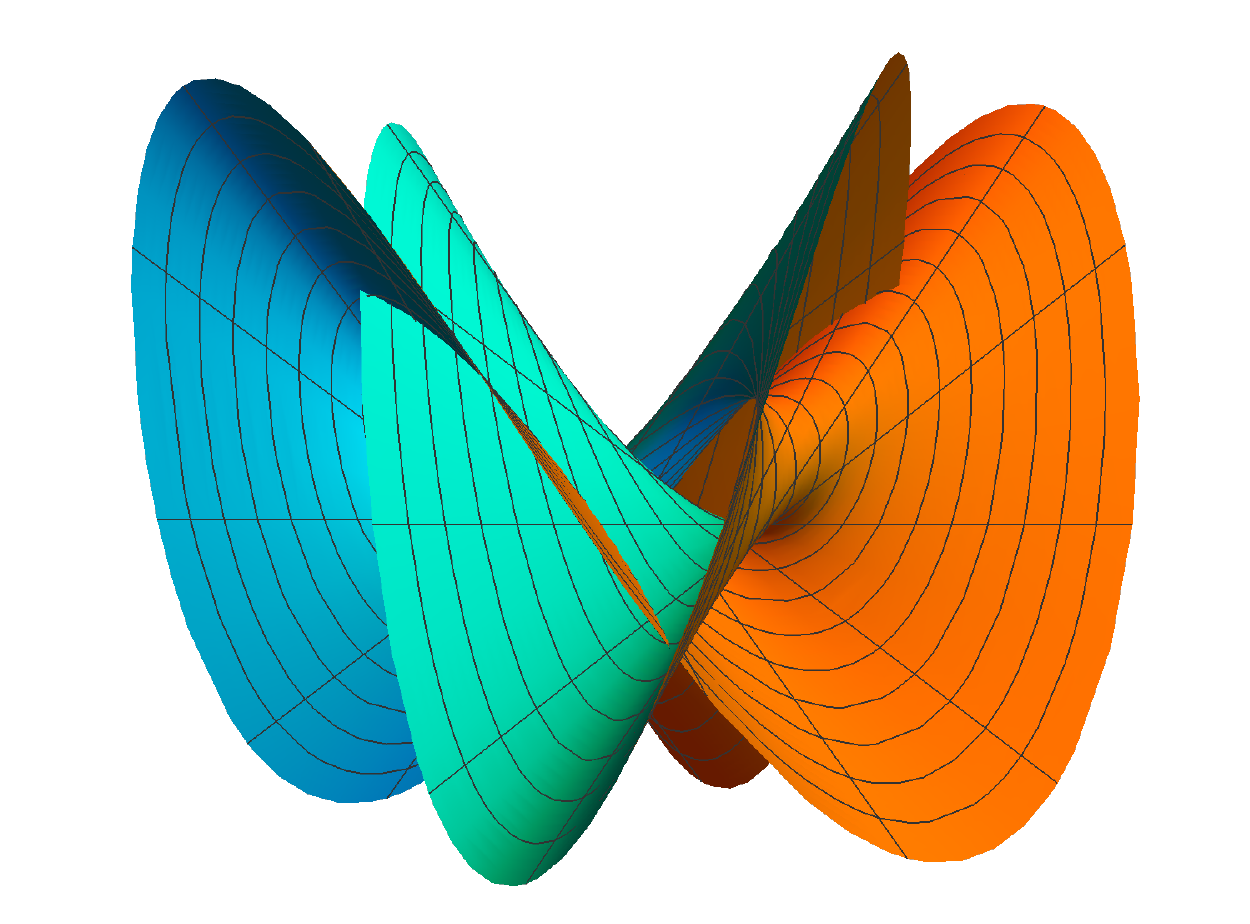}
\caption{The $n=2$ surface (\ref{twormn}) at $D_0$ (top) and $D_2$ (bottom), $t=\mp 1$ viewed from the same angles as in Fig. \ref{hh}. Again they are best visible on the left.}
\label{dd}
\end{figure}

The $n=2$ surface captures several other topological events, depicted in Fig. \ref{tt} and \ref{dd}. We have $D_0=D_2=(0,0,0)$ at $t=\pm 1$. The points 
$T_\pm$ (birth/death of three-fold intersections, Appendix C) occur at $t=\pm(\sqrt{17}-3)/2\simeq \pm 0.56$. They are located at $z=0$ and  $x=-y=\pm(\sqrt{17}-5)/2\sqrt{2}$ at $T_+$ while
$x=y=\pm(\sqrt{17}-5)/2\sqrt{2}$ at $T_-$ (see details in Appendix F).

\begin{figure}
\includegraphics[scale=0.3]{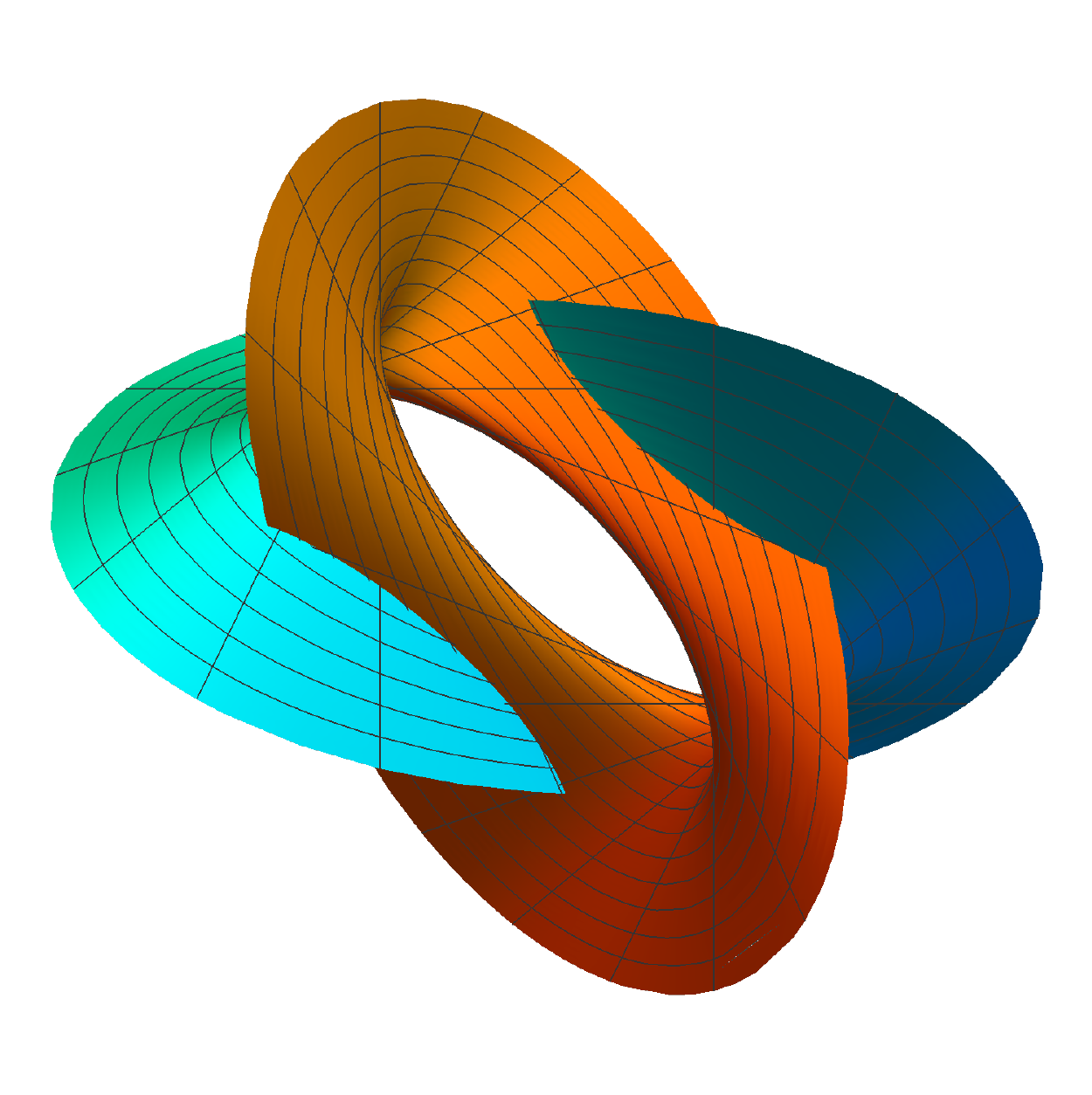}
\includegraphics[scale=0.3]{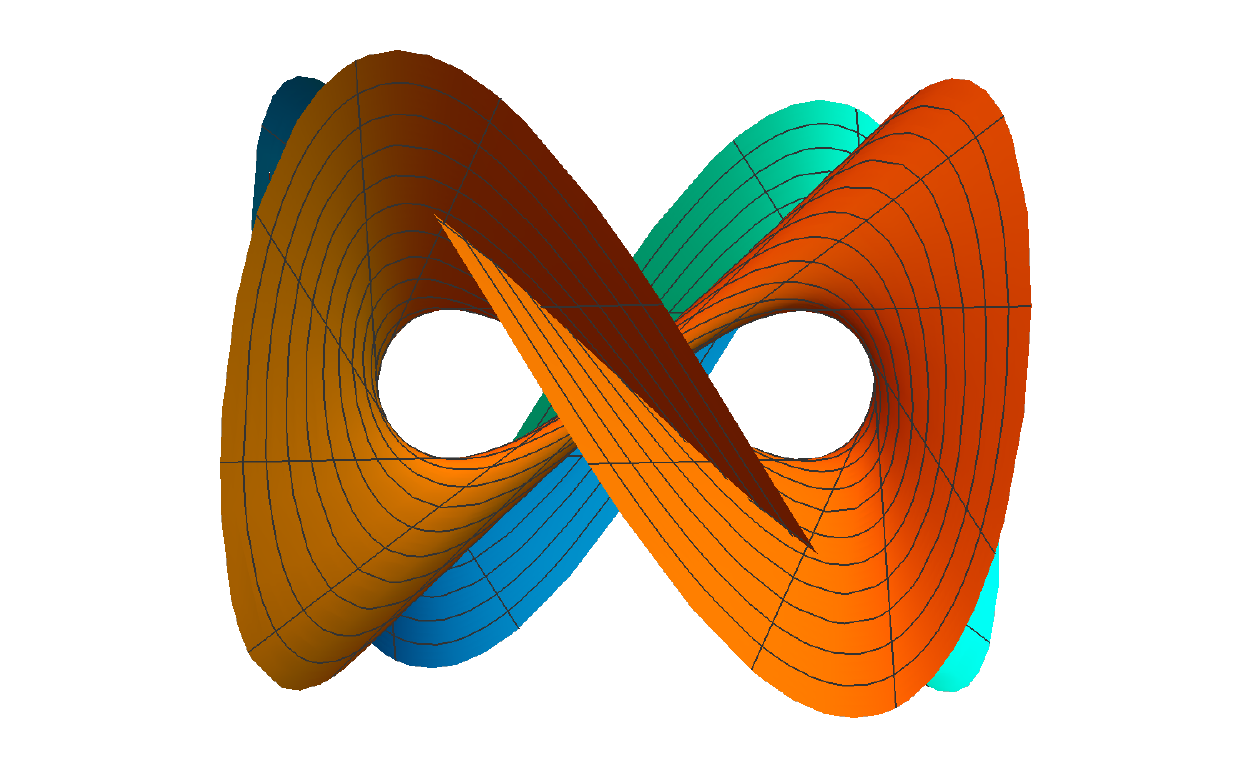}
\includegraphics[scale=0.3]{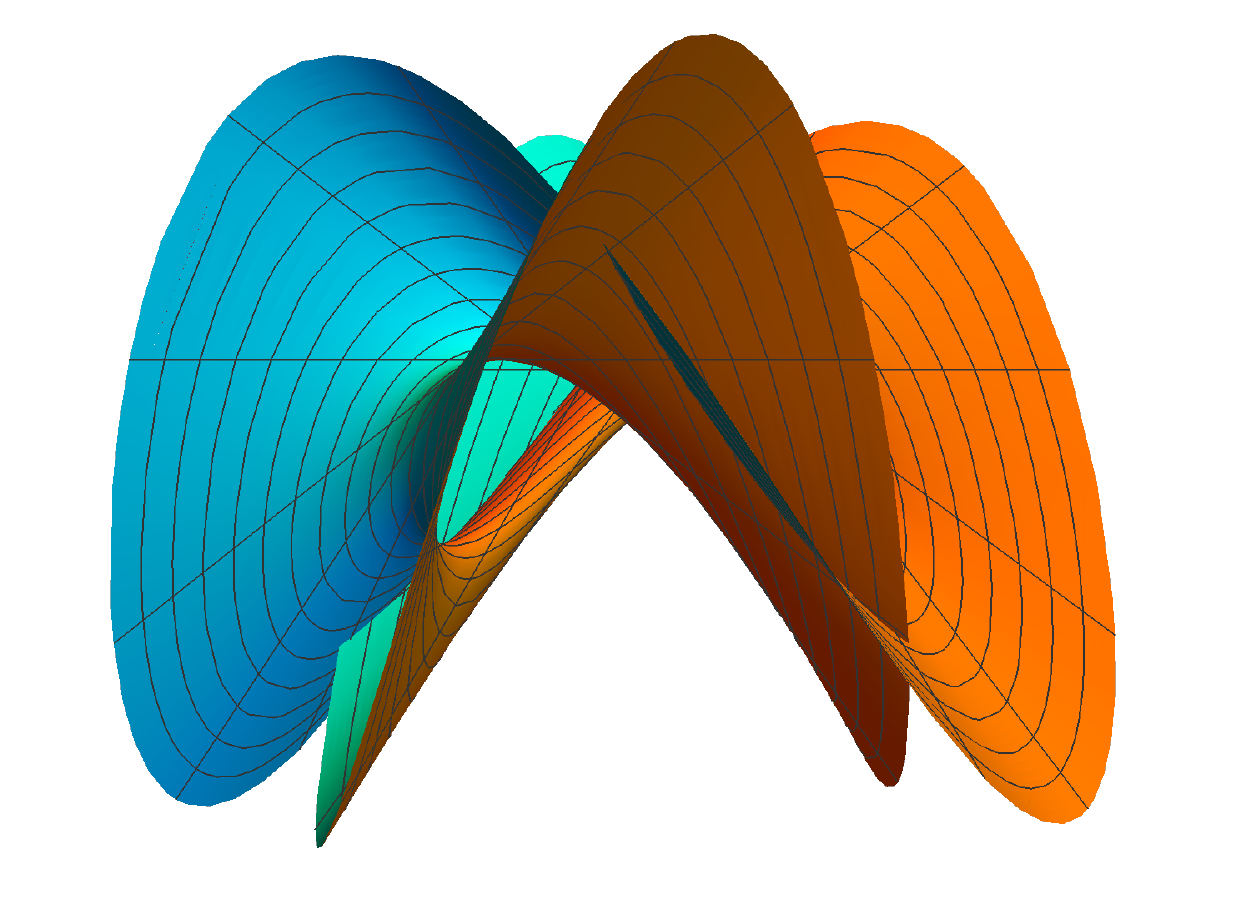}
\includegraphics[scale=0.3]{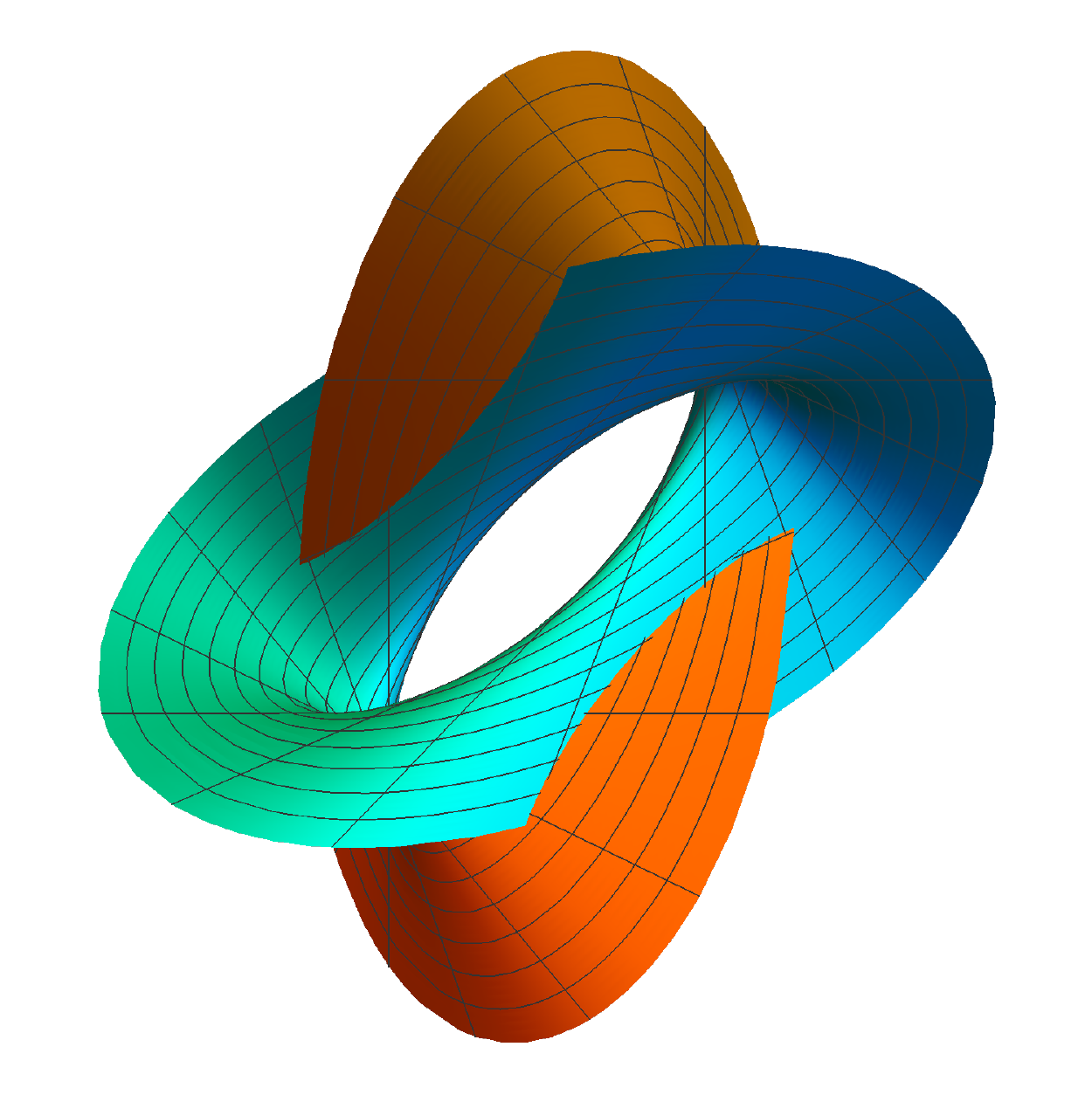}
\includegraphics[scale=0.3]{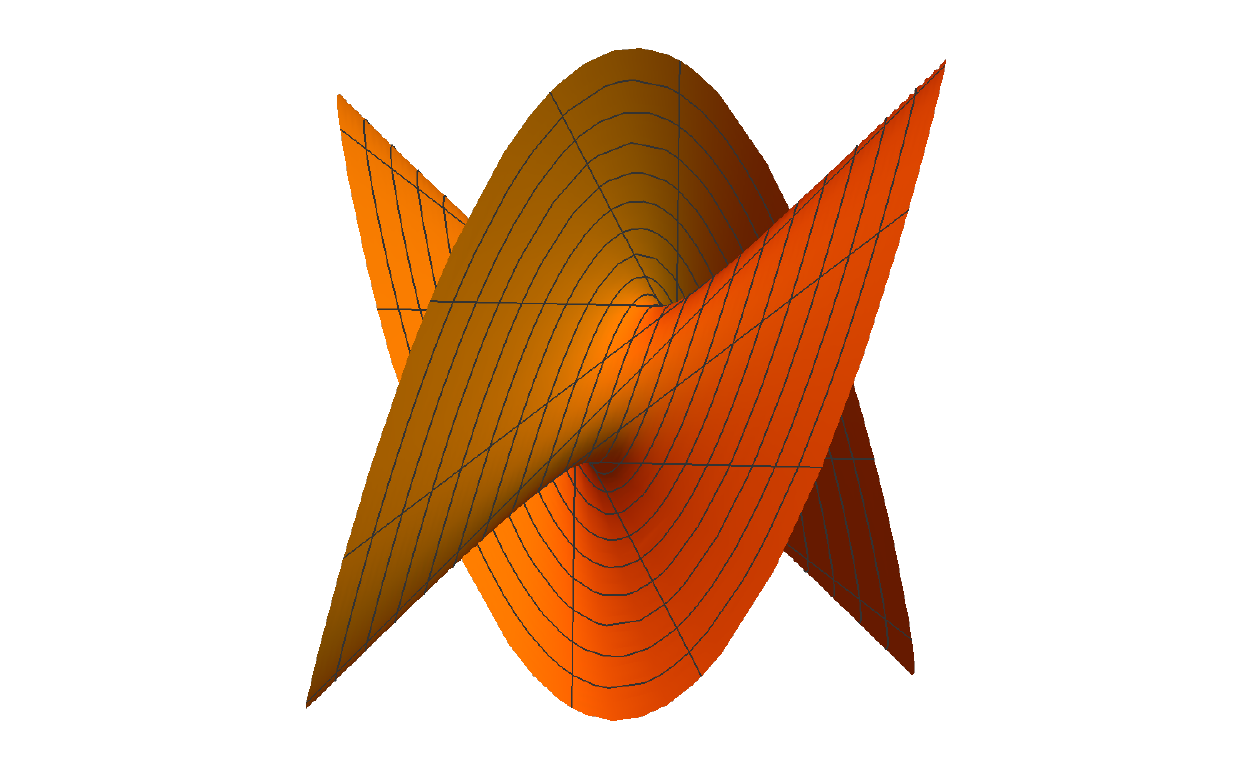}
\includegraphics[scale=0.3]{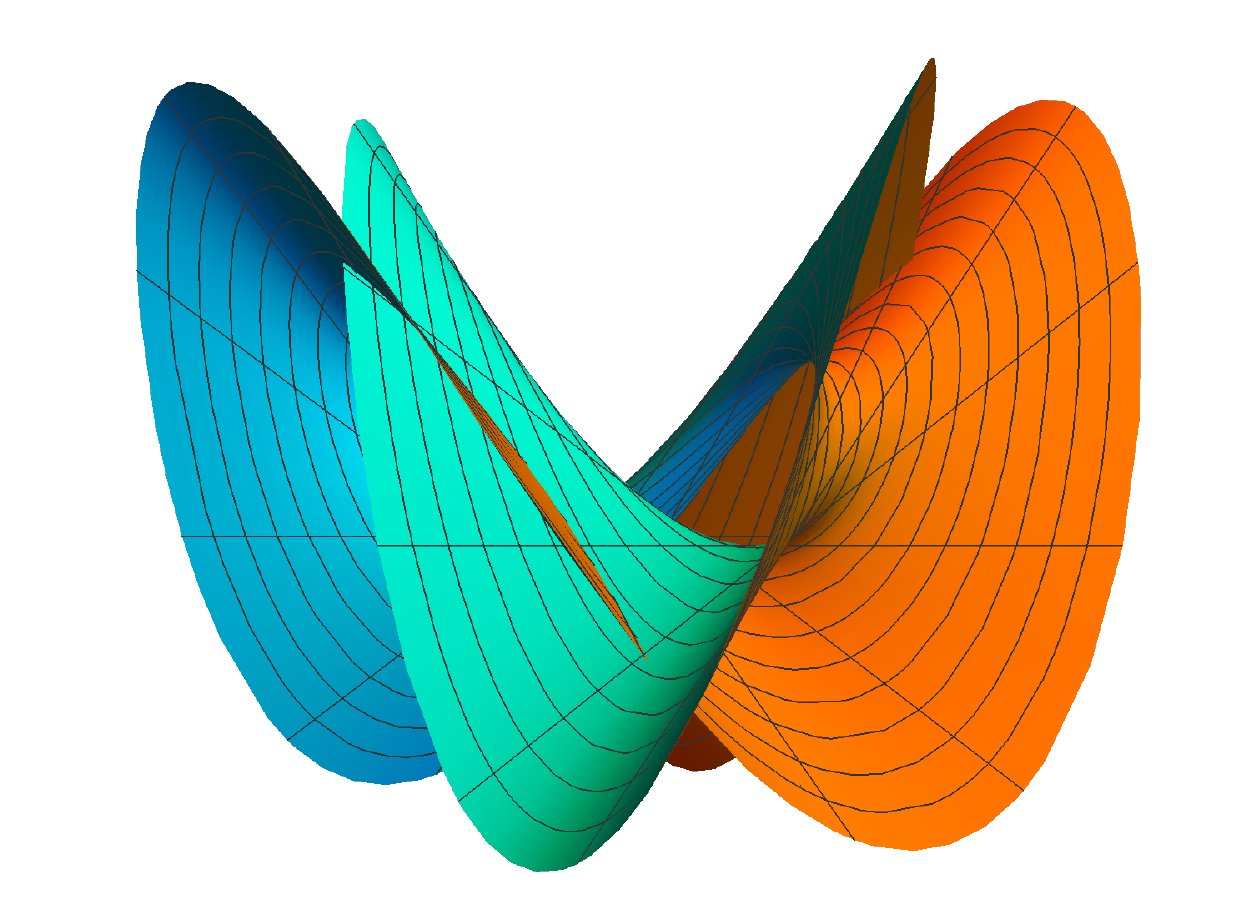}
\caption{The $n=2$ surface (\ref{twormn}) at  $t=-3/2$ (top) and $+3/2$ (bottom) viewed from the same angles as in Fig. \ref{hh}.}
\label{ww}
\end{figure}

The $n=2$ surface for $|t|>1$ resembles more or less a wormhole and has simple self-intersections, see Fig. \ref{ww} and Appendix D.
For sufficiently large $|t|$ also higher $n$-surfaces have simple self-intersections because the domain splits at $h=(n-1)\cos n\phi$ into two
embedded surfaces.  We only need to finish the eversion  recovering a sphere.

\section{Unfolding the wormhole}

\begin{figure}
\includegraphics[scale=0.3]{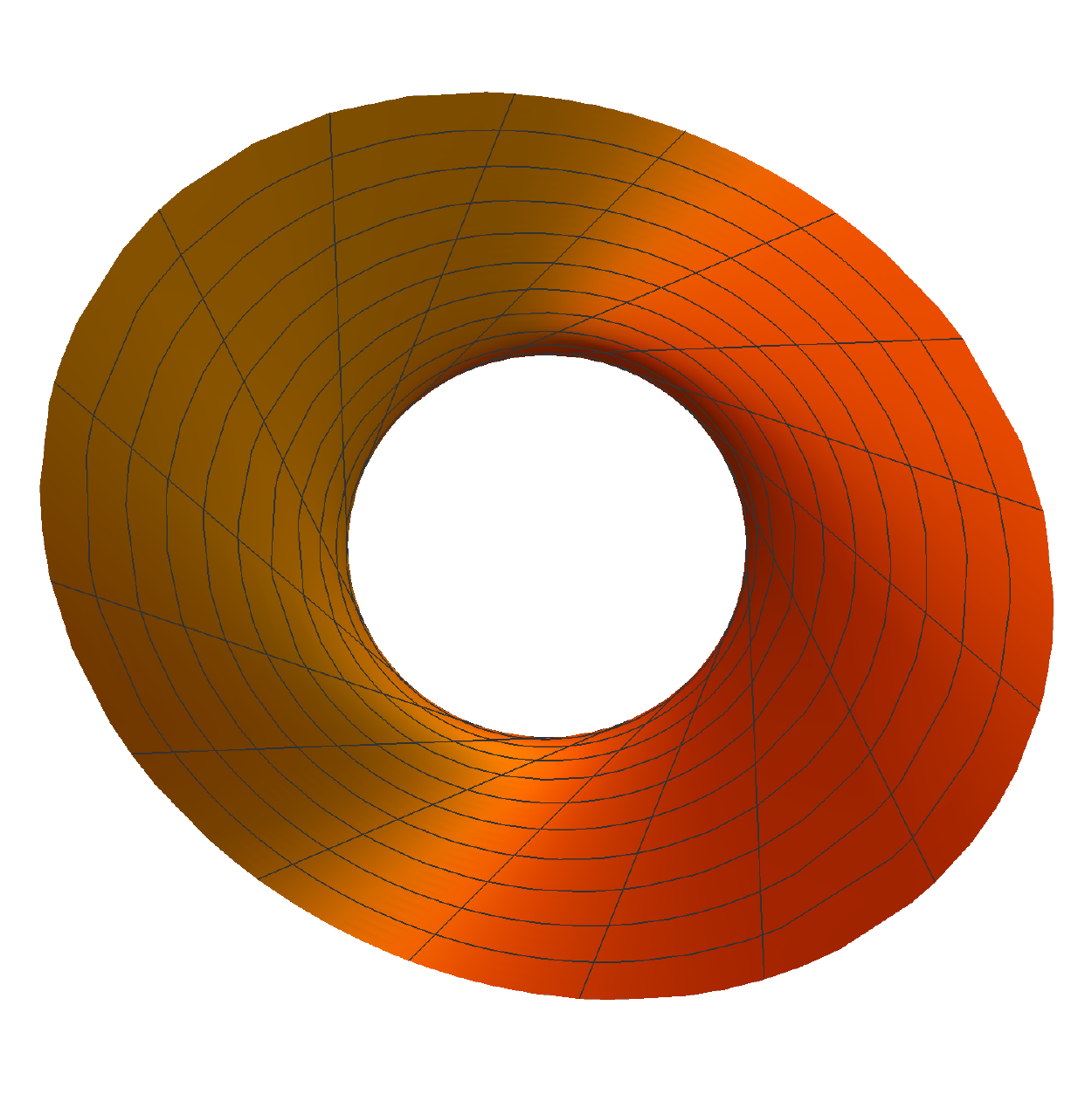}
\includegraphics[scale=0.3]{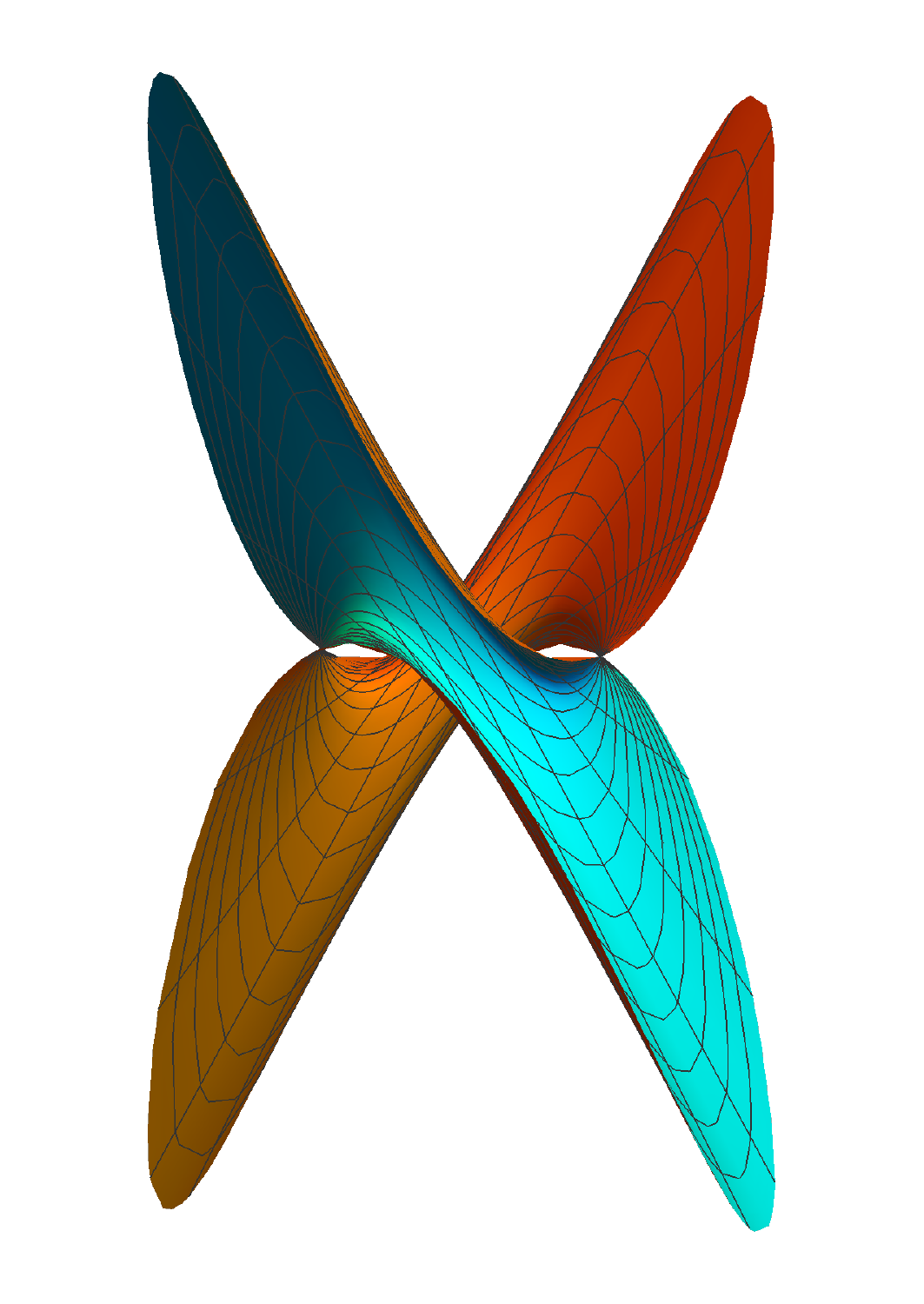}
\includegraphics[scale=0.3]{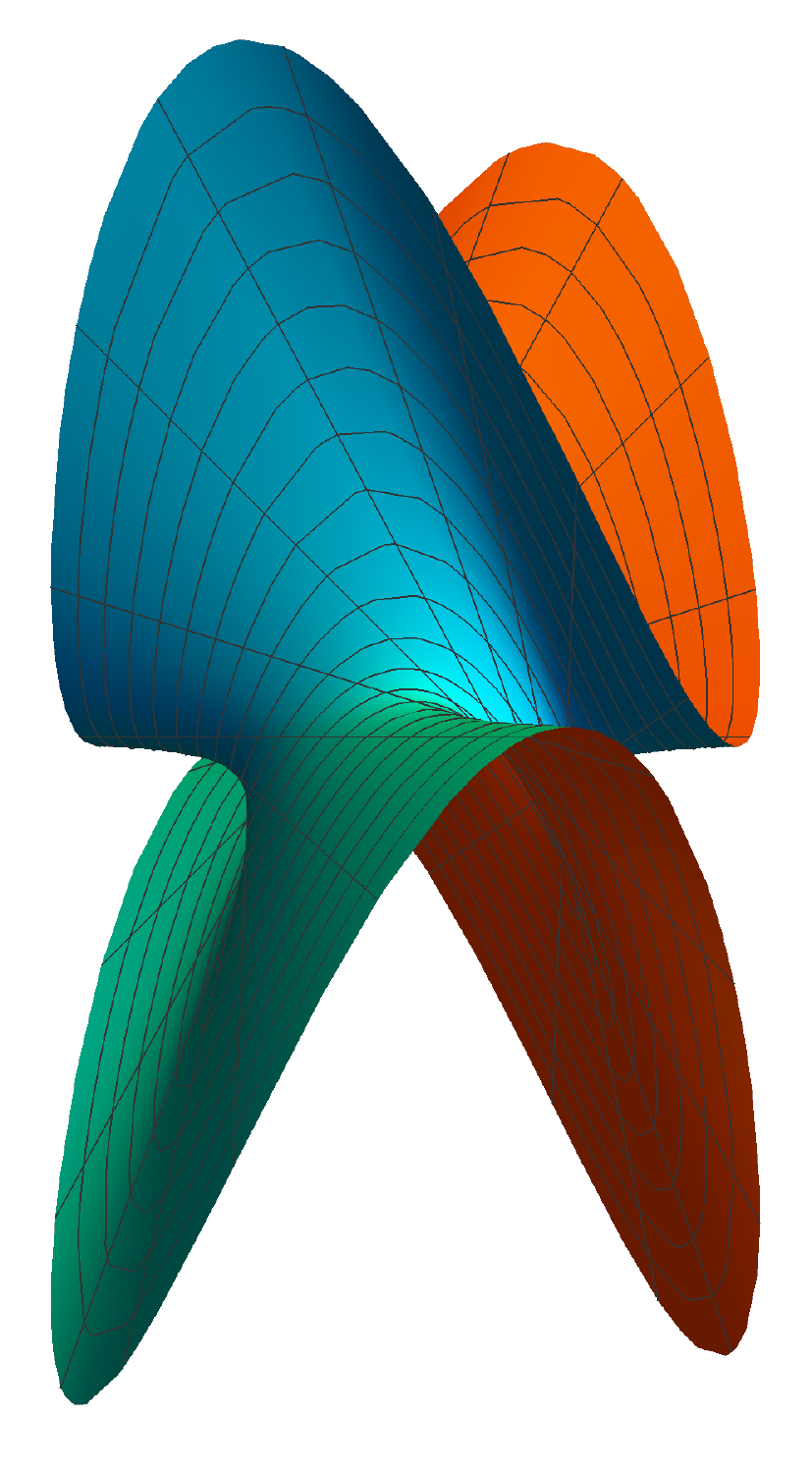}
\includegraphics[scale=0.3]{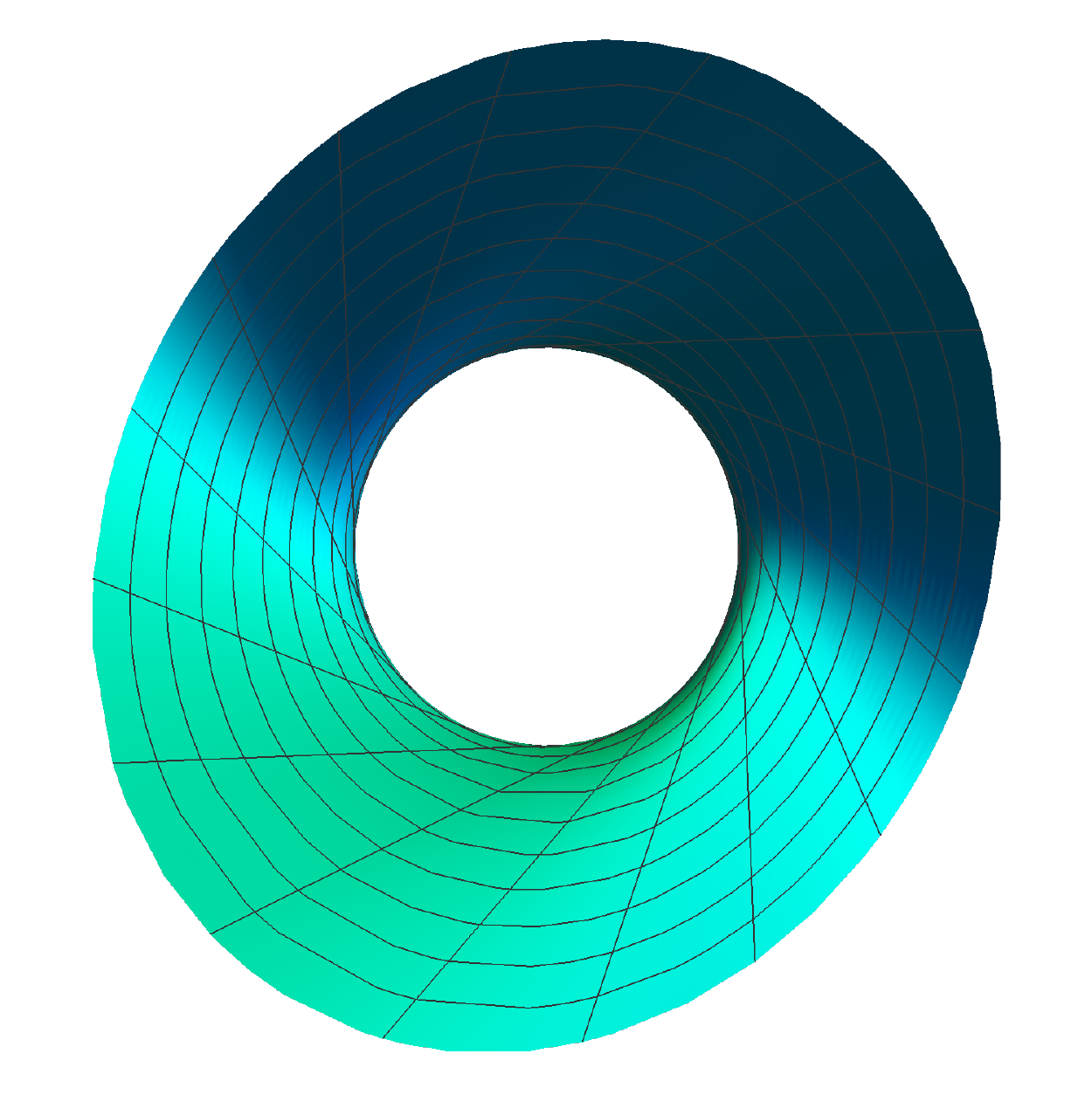}
\includegraphics[scale=0.3]{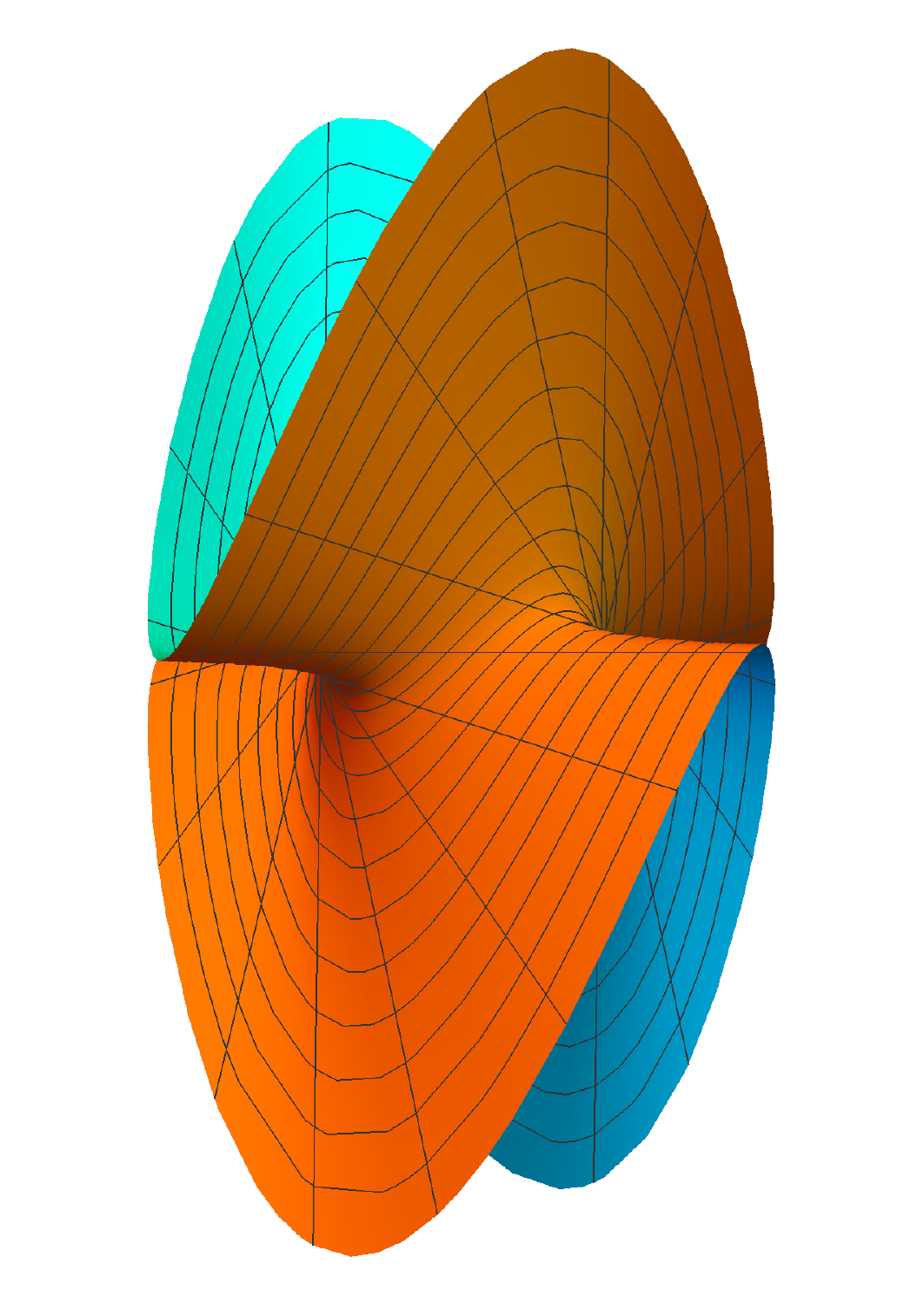}
\includegraphics[scale=0.3]{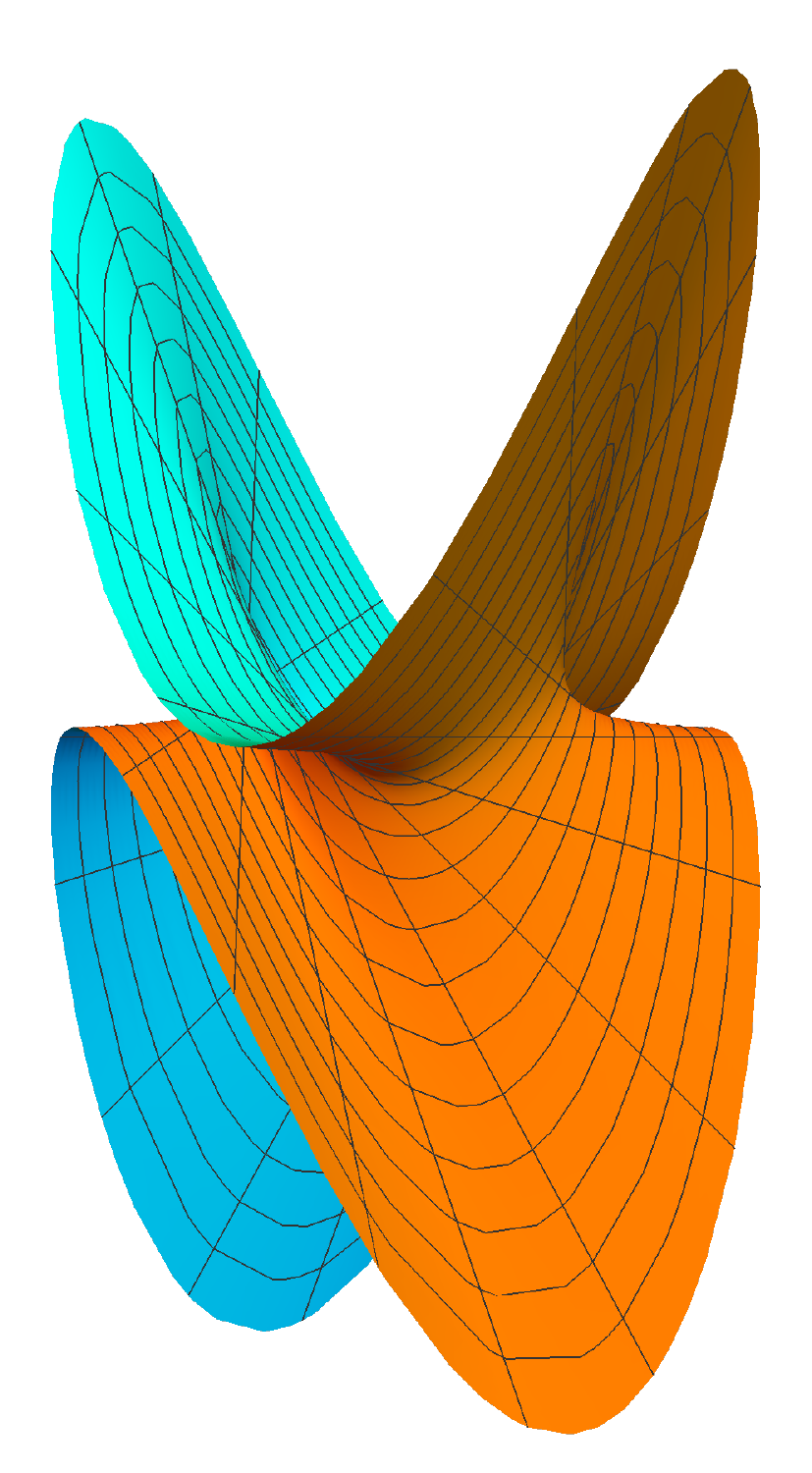}
\caption{The surface (\ref{reg}) at  $t=-3/2$ (top) and $+3/2$ (bottom) ad $q=2/3$ viewed from the same angles as in Fig. \ref{hh}.}
\label{uw}
\end{figure}

Towards full eversion we need to free the surface from self-intersections. To this end we further generalize (\ref{twormn})
into
\begin{equation}
\begin{matrix}
x=t\cos \phi+p\sin(n-1)\phi-h\sin\phi\\
y=t\sin\phi+p\cos(n-1)\phi+h\cos\phi\\
z=h\sin n\phi-(t/n)\cos n\phi-qth
\end{matrix}\label{tworm2}
\end{equation}
with $q\geq 0$.
It is smooth on condition (Appendix A)
\begin{equation}
(n-1)p(1-q|t|)+qt^2>0.\label{conpqt}
\end{equation}
For $n=2$ it is sextic (Appendix B) and has regular self-intersections (Appendix D).
We start from $|t|>1$ (fixed), $p=1$ and $q=0$ to end at $p=0$, $qt=\pm 1$, namely
\begin{equation}
\begin{matrix}
x=t\cos\phi-h\sin\phi\\
y=t\sin\phi+h\cos\phi\\
z=h\sin n\phi-(t/n)\cos n\phi\mp h
\end{matrix}\label{reg}
\end{equation}
depicted in Fig. \ref{uw} for $n=2$.
At $qt=\pm 1$ the intersection disappears at infinity ($D_{01}$ or $D_{21}$ point, Appendix C).
It is convenient to keep $p=1-|qt|\geq 0$ when (\ref{conpqt}) is always satisfied.

\section{Inversion}

As our goal, the traditional $S^2$ sphere is given by the equation $|\vec{R}|^2=X^2+Y^2+Z^2=R^2$ (with a constant radius $R>0$),
we will parametrize it by $\phi\in [-\pi,\pi]$ (periodic) and $\theta\in [-\pi/2,\pi/2]$ using $X=R\cos\theta\cos\phi$, $Y=R\cos\theta\sin\phi$, $Z=R\sin\theta$, and
 map onto the cylinder, $\phi=\phi$, $h=\omega\sin\theta/\cos^n\theta$, with some $\omega>0$.

\begin{figure}
\includegraphics[scale=0.3]{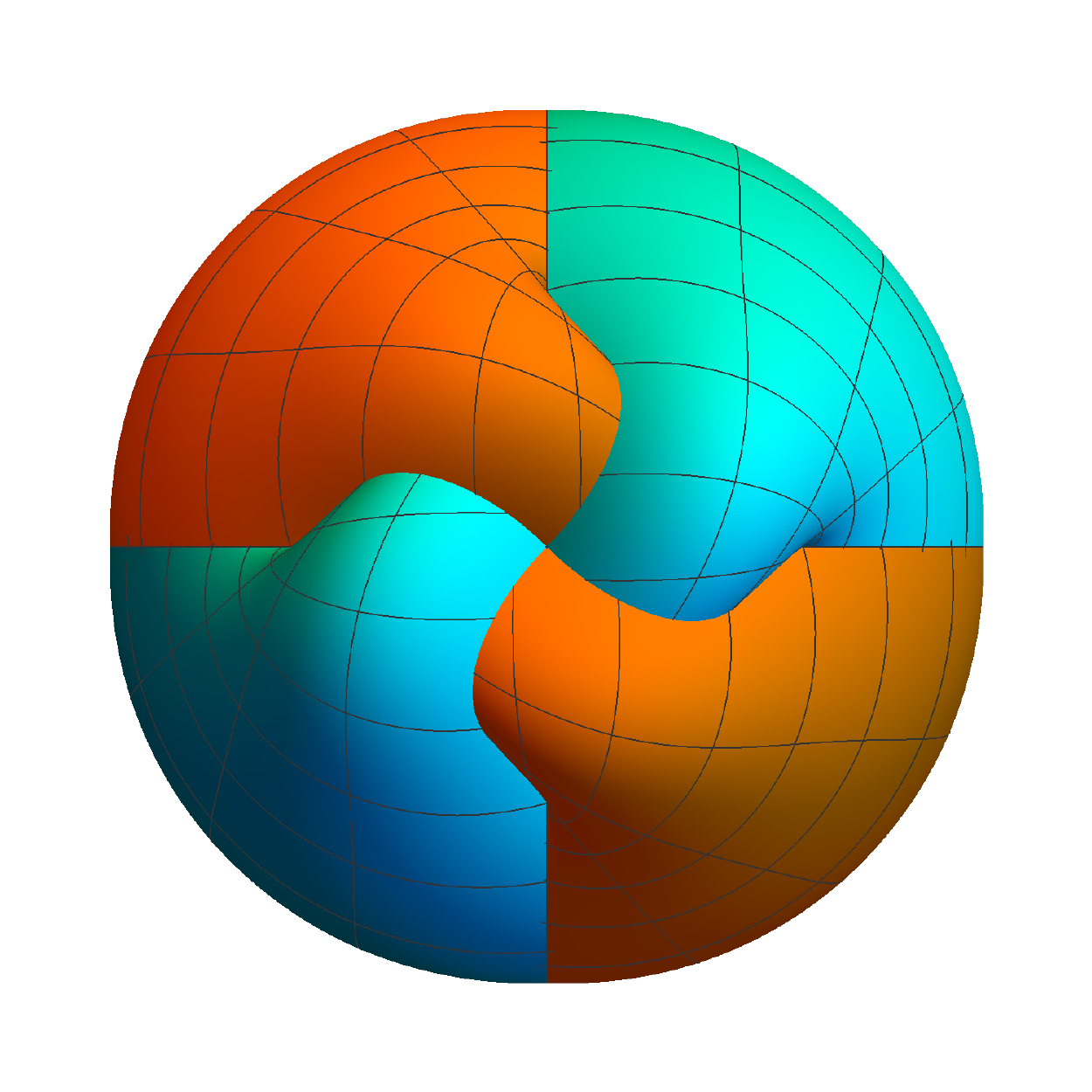}
\includegraphics[scale=0.3]{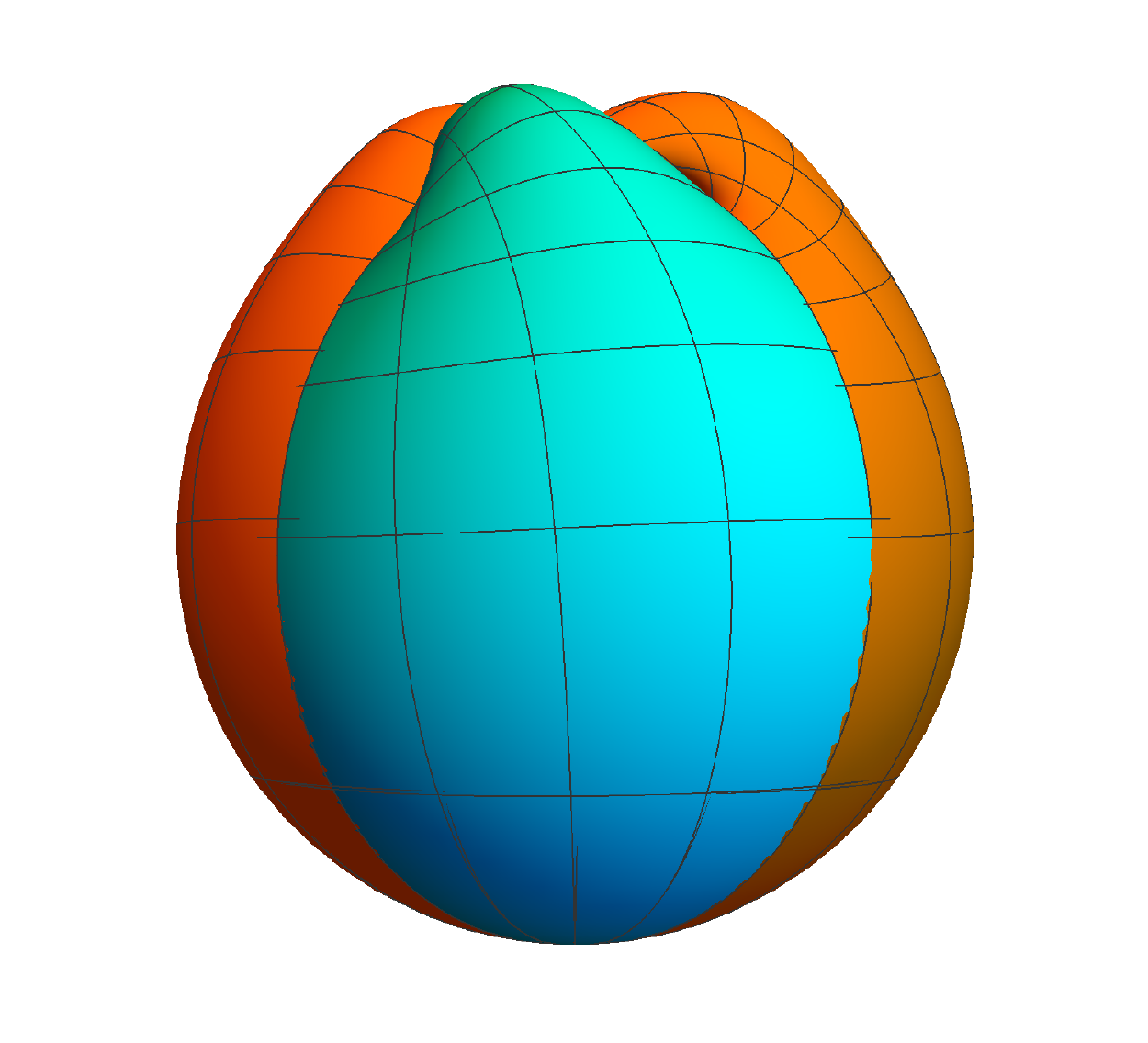}
\includegraphics[scale=0.3]{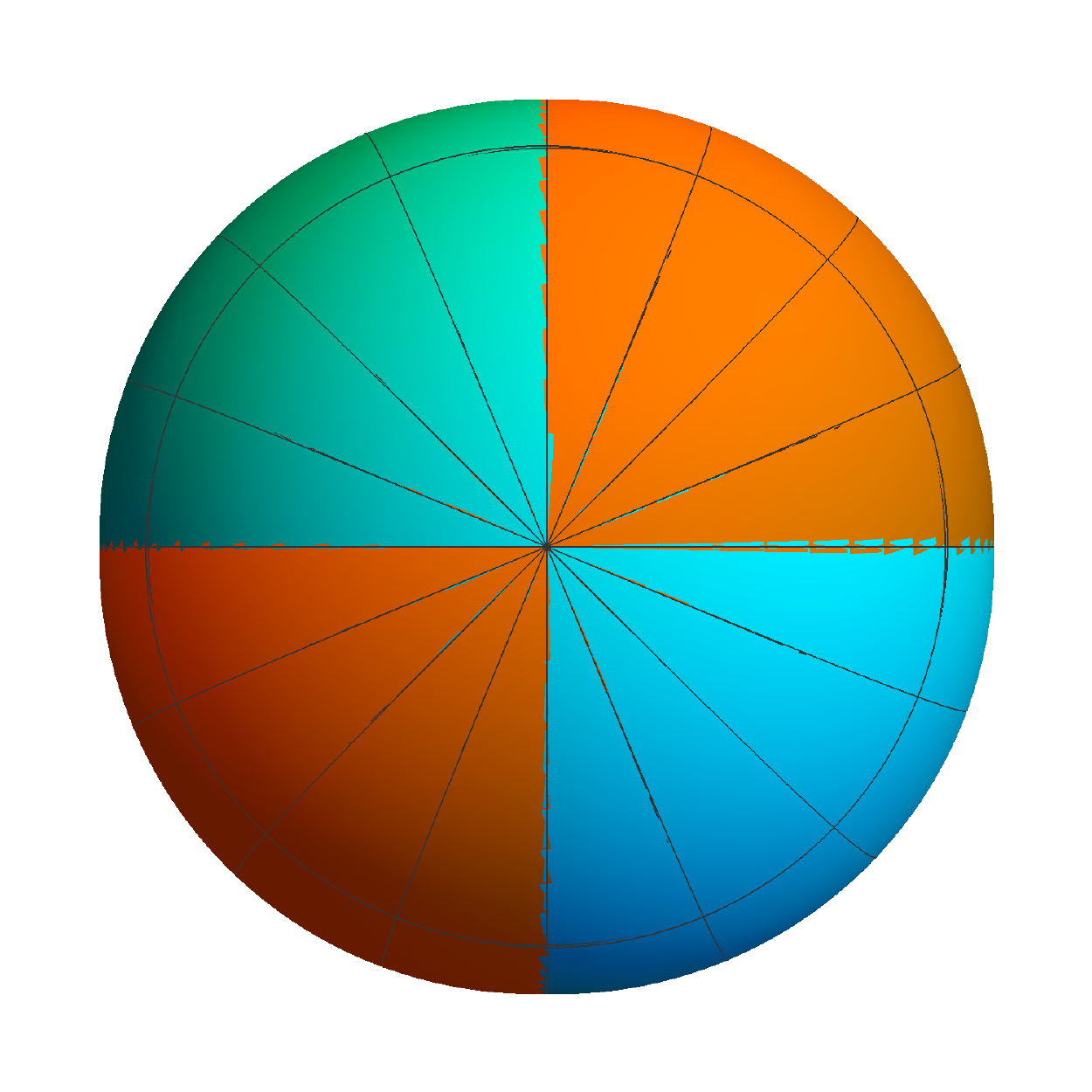}
\caption{The inversion of  $n=2$ halfway surface (\ref{halfn}) using (\ref{wsm1}) and (\ref{wsm2}) with, $\alpha=1$, $\beta=1/25$, $\xi=\eta=1$, $\lambda=1$, $\omega=2$, viewed from the top, side and bottom}
\label{s-t0}
\end{figure}

\begin{figure}
\includegraphics[scale=0.3]{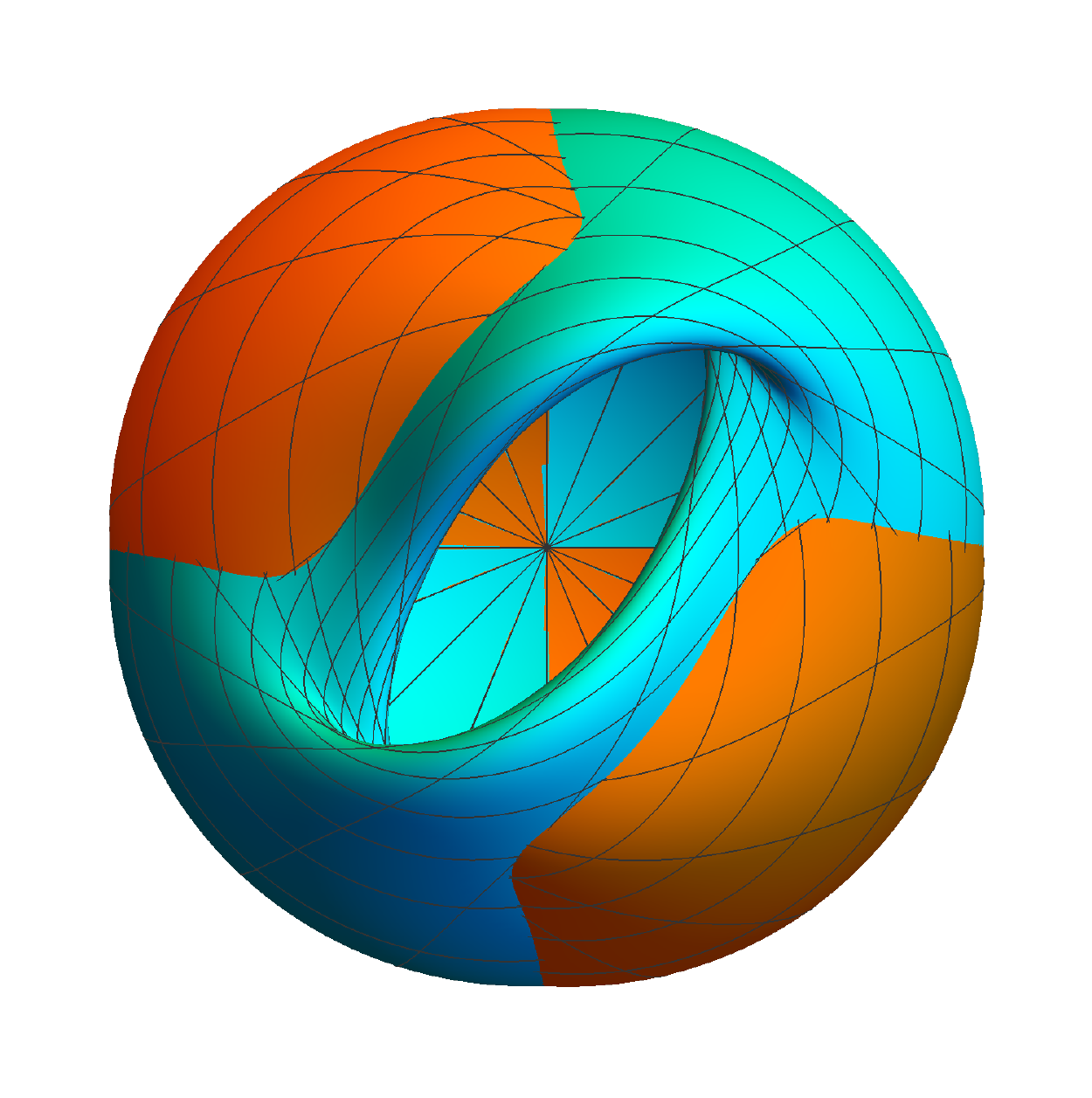}
\includegraphics[scale=0.3]{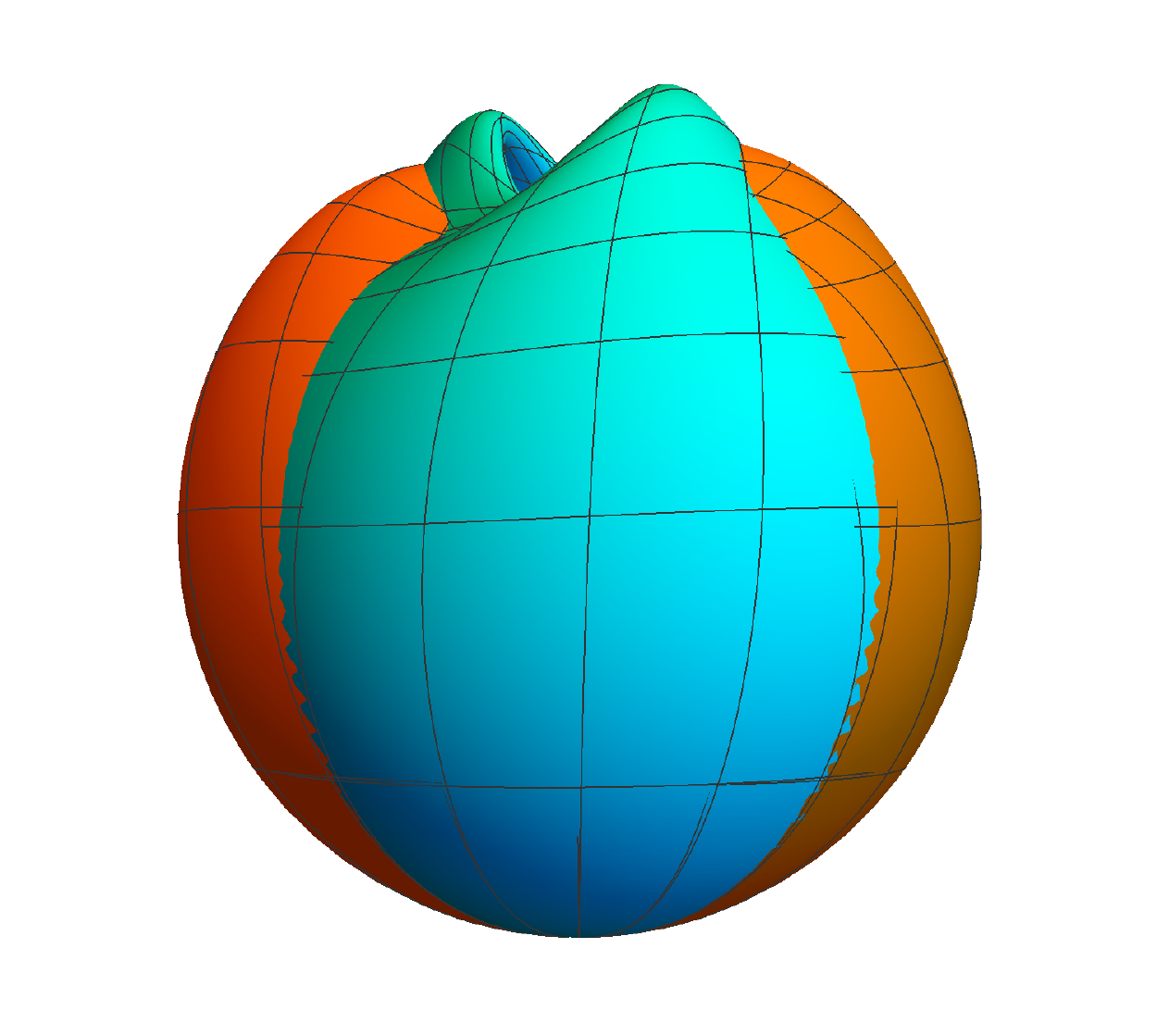}
\includegraphics[scale=0.3]{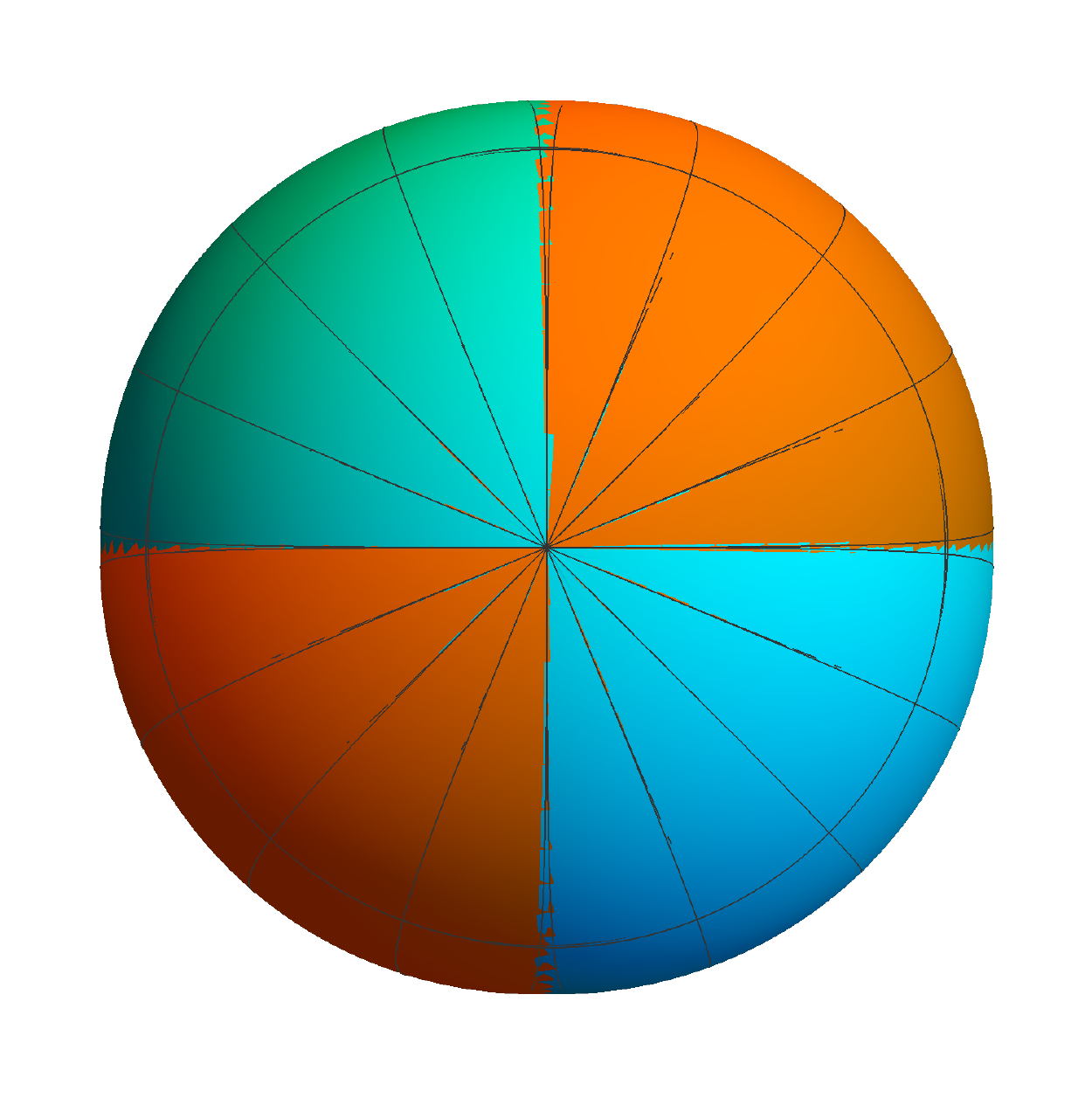}
\caption{The inversion of $n=2$ surface (\ref{tworm2}) using (\ref{wsm1}) and (\ref{wsm2}) with $t=3/2$, $\alpha=1$, $\beta=1/25$, $q=0$, $\xi=\eta=1$, $\lambda=1$, $\omega=2$, Viewpoints as in Fig. \ref{s-t0}}
\label{s-t32}
\end{figure}

\begin{figure}
\includegraphics[scale=0.3]{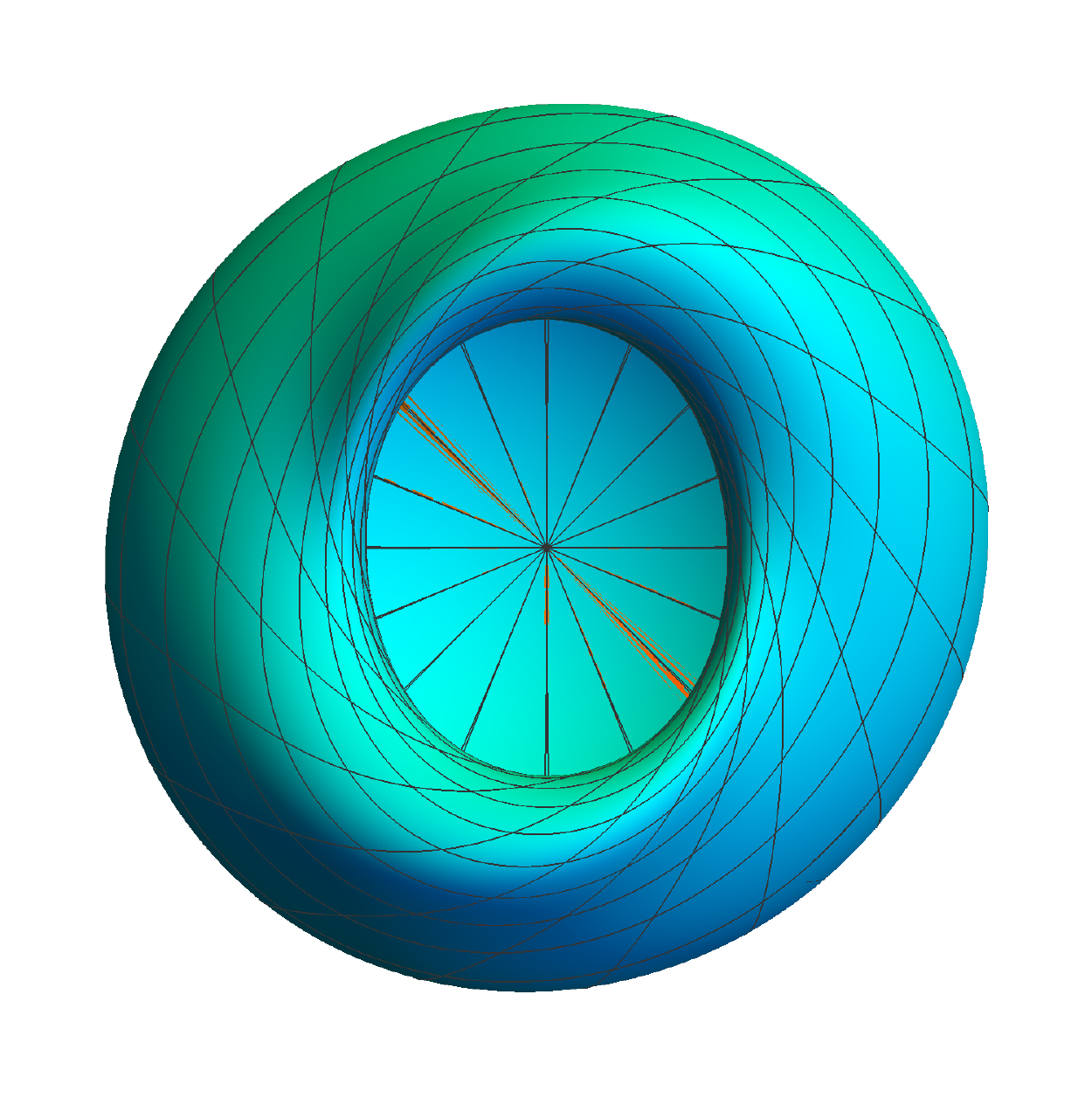}
\includegraphics[scale=0.3]{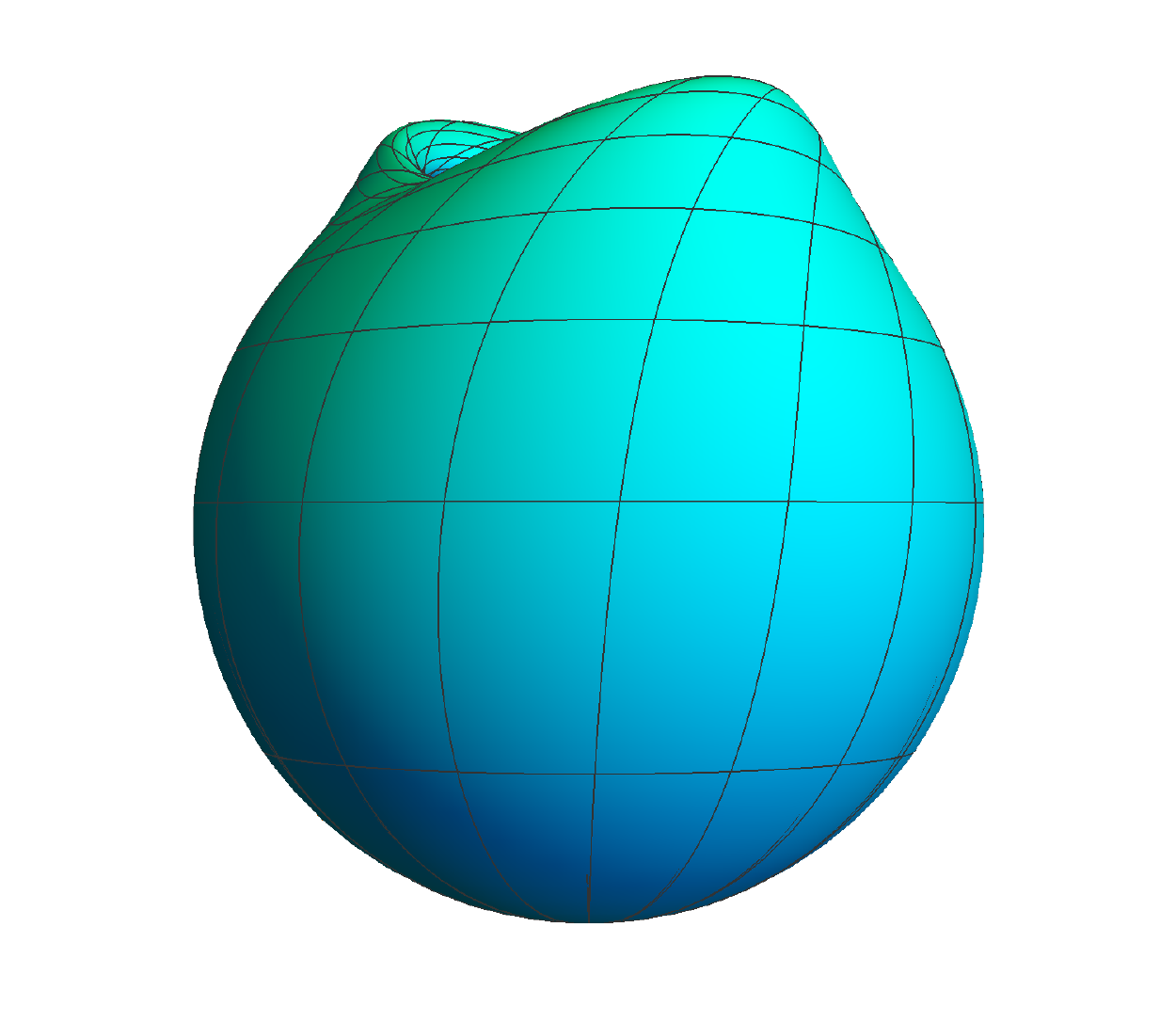}
\includegraphics[scale=0.3]{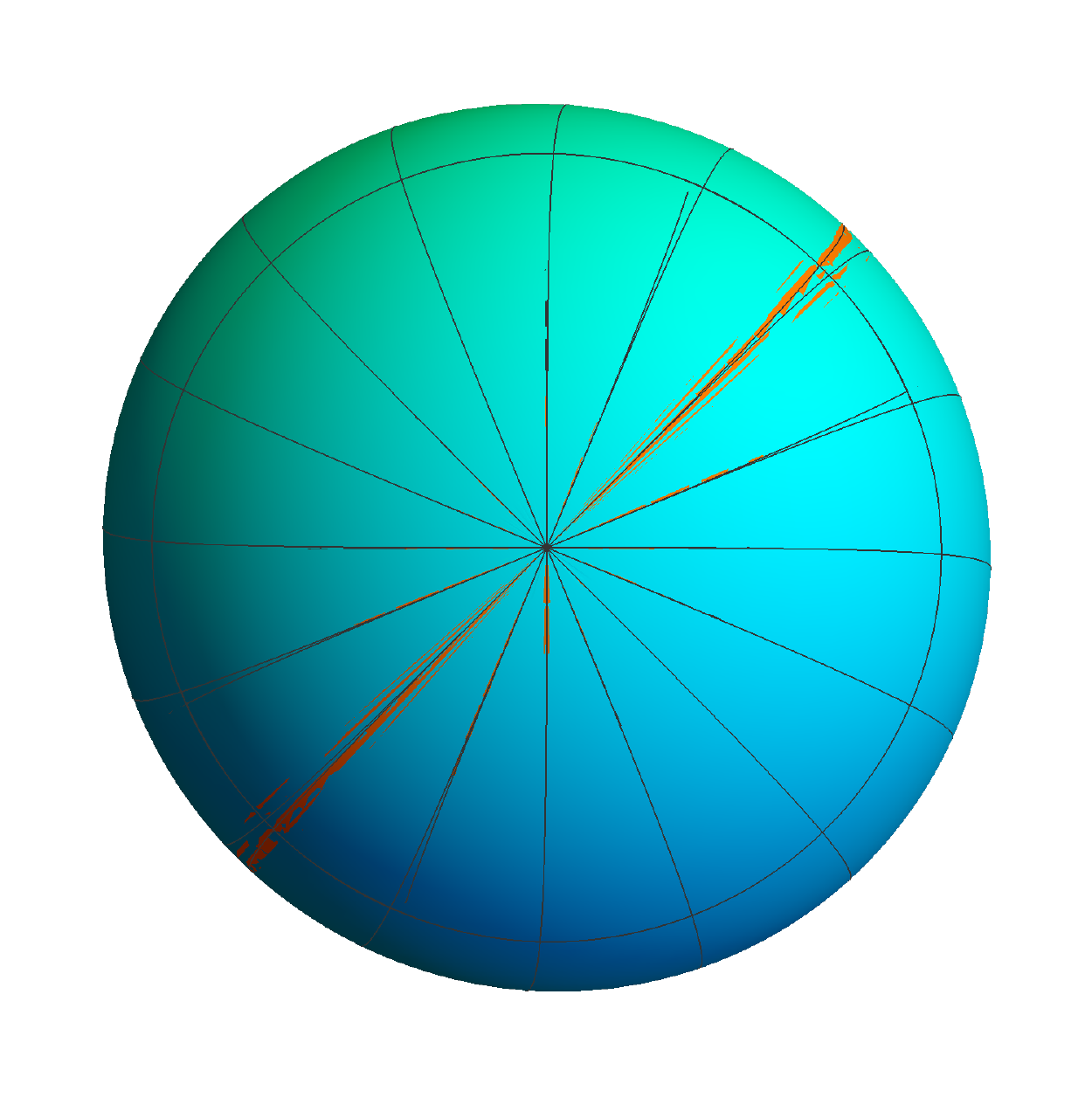}
\caption{The inversion of $n=2$ surface (\ref{reg}) using (\ref{wsm1}) and (\ref{wsm2}) with $t=3/2$, $\alpha=1$, $\beta=1/25$, $q=2/3$, $\xi=\eta=1$, $\lambda=1$, $\omega=2$, Viewpoints as in Fig. \ref{s-t0}}
\label{s-t32q}
\end{figure}

\begin{figure}
\includegraphics[scale=0.3]{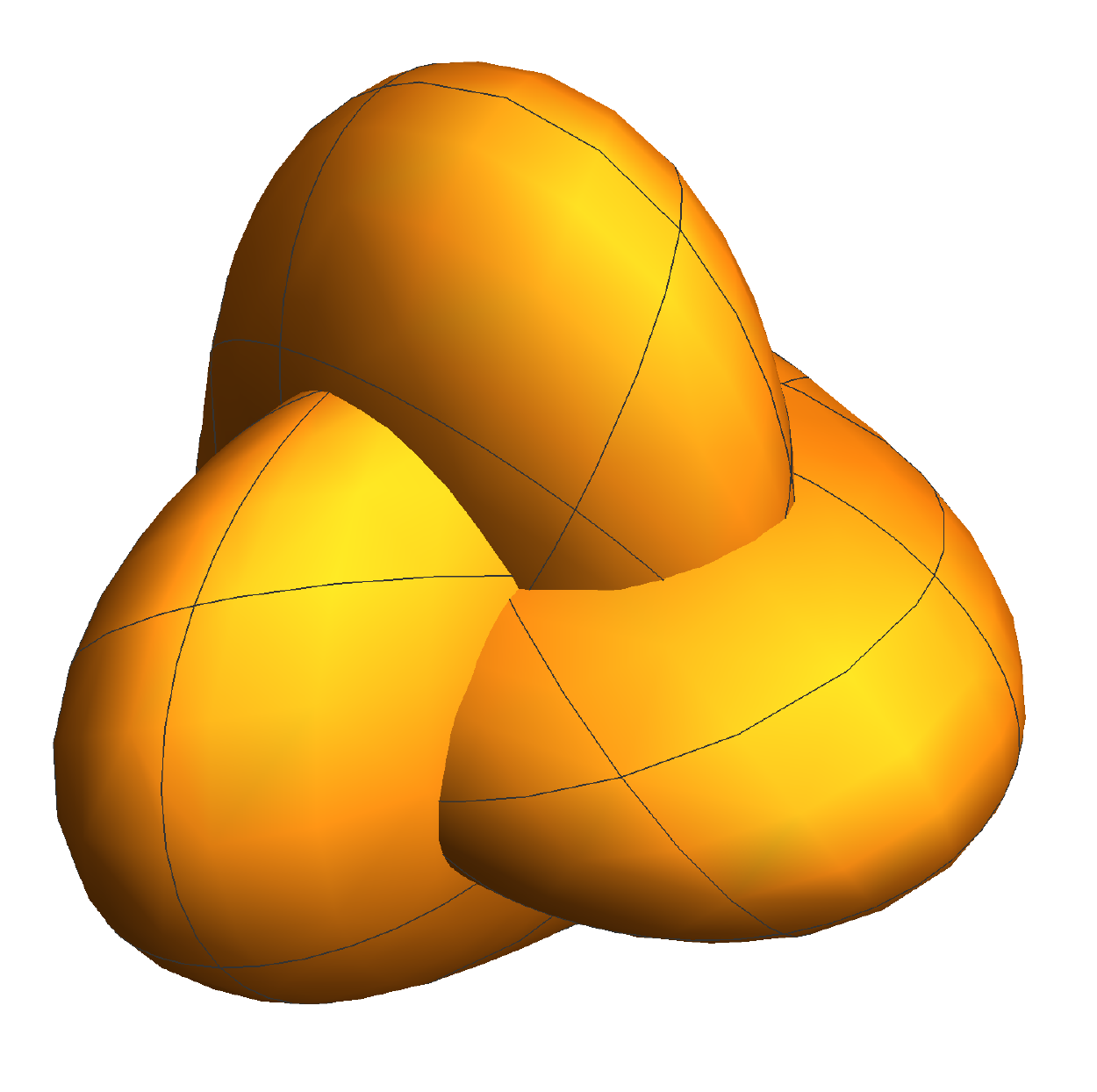}
\includegraphics[scale=0.3]{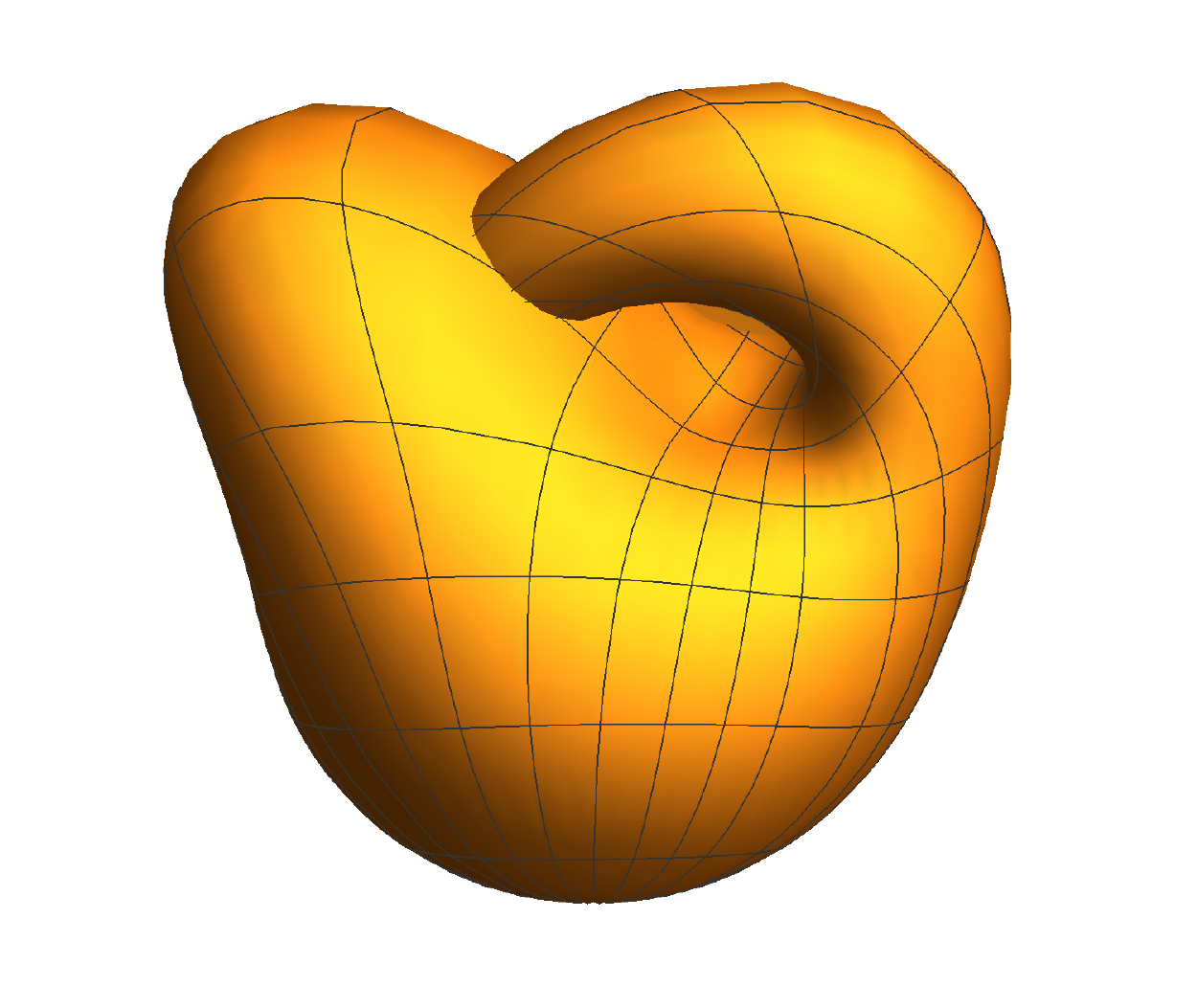}
\includegraphics[scale=0.3]{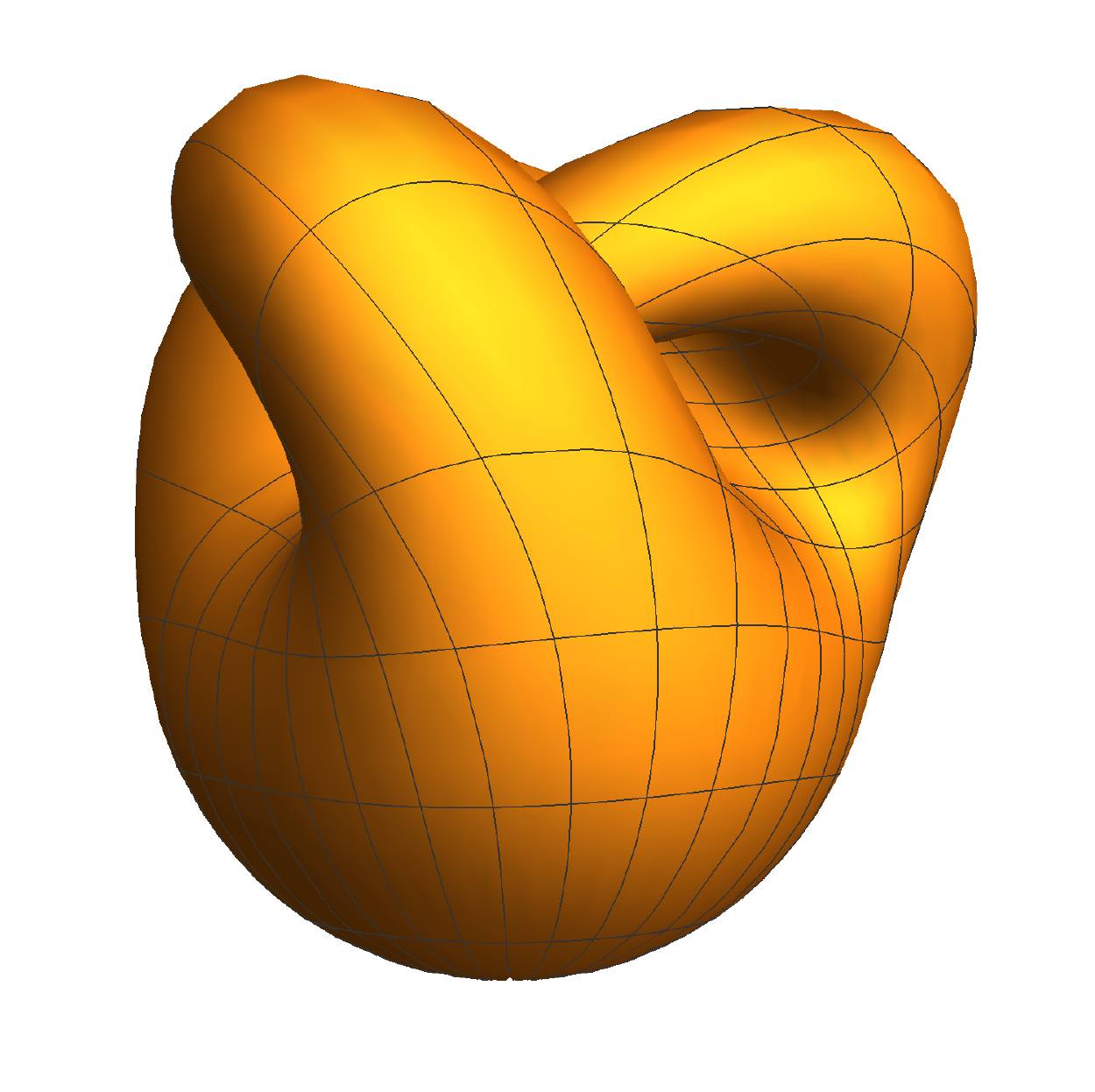}
\caption{The  closed Boy surface (\ref{halfn}) for $n=3$ using (\ref{wsm1}) and (\ref{wsm2}) with $\xi=\eta=1$, $\omega=2$, $\alpha=1$, $\beta =1/4$ viewed as in Fig. \ref{hh}.}
\label{boycl}
\end{figure}

We will close the wormhole at infinity with help of stereographic projection \cite{morin1}. 
The cylinder is mapped
as follows. We add damping at large distances for smoothness
\begin{eqnarray}
&&x'=x(\xi+\eta(x^2+y^2))^{-\kappa},\nonumber\\
&&y'=y(\xi+\eta(x^2+y^2))^{-\kappa},\nonumber\\
&&z'=z/(\xi+\eta(x^2+y^2))\label{wsm1}
\end{eqnarray}
with $x,y,z$ defined by (\ref{tworm2}), some $\xi,\eta\geq 0$ (for $|t|\leq 1$ we keep $\xi>0$)
and $\kappa=(n-1)/2n$ and then
\begin{eqnarray}
&&x''=x' e^{\gamma z'}/(\alpha+\beta(x^{\prime 2}+y^{\prime 2})),\nonumber\\
&&y''=y'e^{\gamma z'}/(\alpha+\beta(x^{\prime 2}+y^{\prime 2})),\label{wsm2}\\
&&z''=\frac{\alpha-\beta(x^{\prime 2}+y^{\prime 2})}{\alpha+\beta(x^{\prime 2}+y^{\prime 2})}\frac{e^{\gamma z'}}{\gamma}-\gamma^{-1}
\frac{\alpha-\beta}{\alpha+\beta}\nonumber
\end{eqnarray}
for $\alpha,\beta\geq 0$ and $\gamma=2\sqrt{\alpha\beta}$. Both mappings are smooth i.e.  $C^\infty$ class for even $n$ or $t=0$ and odd $n$ while $C^1$  in other cases (see Appendix G also how to extend $C^1$ to $C^\infty$), except the case $\alpha=0$ and $\xi=0$ for $|t|\leq 1$. 
The (geometric) mean radius of the inversion sphere is $\gamma^{-1}$. The case $\xi=1$, $\eta=0$, $\alpha=1$, $\beta\to 0$ corresponds to original open wormhole. Although one could replace (\ref{wsm2}) by e.g. a standard inversion $\vec{r}''=(\vec{r}'-\vec{r}_0)/|\vec{r}'-\vec{r}_0|^2$ for some $\vec{r}_0$ away from the wormhole (preferably on $z$ axis) we stick to (\ref{wsm2}) which preserves inversion symmetry. 
One can see that the points $h\to\pm\infty$ are smoothly (see Appendix G) mapped onto $(0,0,-\sqrt{\alpha/\beta}(\alpha+\beta))$ event at $x,y\to\infty$, $z\to 0$, which is a $D_1$ point at $n=2$.
To close the wormhole we need $\eta,\beta>0$,
$\alpha,\beta,\xi,\eta$ can depend on $t$ and we want $\beta=1$, $\alpha=0$, $\xi=0$, $|t|>1$, in the final stage, meaning inversion of $xy$ plane,
\begin{eqnarray}
&&x''=x'/(x^{\prime 2}+y^{\prime 2}),\nonumber\\
&&y''=y'/(x^{\prime 2}+y^{\prime 2}),\label{wsm3}\\
&&z''=-z'\nonumber
\end{eqnarray}
or
\begin{eqnarray}
&&x''=\frac{\eta^\kappa x}{(x^2+y^2)^{1-\kappa}},\nonumber\\
&&y''=\frac{\eta^\kappa y}{(x^2+y^2)^{1-\kappa}},\label{wsm4}\\
&&z''=-\frac{z/\eta}{x^2+y^2},\nonumber
\end{eqnarray}
 which completes the inversion process.
 The full eversion map reads 
\begin{equation}
(\theta,\phi)\to(h,\phi)\to\vec{r}=(x,y,z)\to\vec{r}'=(x',y',z')\to\vec{r}''=(x'',y'',z'').
\end{equation}
 The inverted stages $t=0$, $t=3/2$ and $q=2/3$ of the $n=2$ surface are depicted in Figs. \ref{s-t0}, \ref{s-t32}, \ref{s-t32q}.
 The inverted, closed Boy surface  at $t=0$ and $n=3$ is depicted in Fig. \ref{boycl}.
 
\section{Unfolding the sphere}

From the end of inversion (\ref{wsm4}), we only need to change $z$ or $z''$ while keeping constant $x,y$ or $x'',y''$.
 We obtained almost a sphere, except that it is squeezed and twisted.
In the last step, after reaching $\alpha=0$, $\beta=1$, $\xi=0$ at $|t|>1$ in eqs. (\ref{wsm1}) and (\ref{wsm2}),
let us consider a new parameter $\lambda\in[0,1]$
\begin{eqnarray}
&&x=(t(1-\lambda+\lambda\cos^n\theta)\cos\phi-\lambda\omega\sin\theta\sin\phi)/\cos^n\theta,\nonumber\\
&&y=(t(1-\lambda+\lambda\cos^n\theta)\sin\phi+\lambda\omega\sin\theta\cos\phi)/\cos^n\theta,\label{xyl}
\end{eqnarray}
corresponding by (\ref{wsm4}) to
\begin{eqnarray}
&&x''=\eta^\kappa\cos\theta\frac{t(1-\lambda+\lambda\cos^n\theta)\cos\phi-\lambda\omega\sin\theta\sin\phi}{(t^2(\lambda\cos^n\theta+(1-\lambda))^2+\lambda^2\omega^2\sin^2\theta)^{1-\kappa}},\nonumber\\
&&y''=\eta^\kappa\cos\theta\frac{t(1-\lambda+\lambda\cos^n\theta)\sin\phi+\lambda\omega\sin\theta\cos\phi}{(t^2(\lambda\cos^n\theta+(1-\lambda))^2+\lambda^2\omega^2\sin^2\theta)^{1-\kappa}}
\end{eqnarray}
for $\lambda\in[0,1]$
so that
\begin{equation}
x^{\prime\prime 2}+y^{\prime\prime 2}=\eta^{2\kappa}\cos^2\theta(t^2(\lambda\cos^n\theta+(1-\lambda))^2+\lambda^2\omega^2\sin^2\theta)^{-1/n}\label{xxyy}
\end{equation}
is a growing function of $\cos^2\theta$ (see Appendix G). Finally,
\begin{eqnarray}
&&z=\lambda(\omega\sin\theta(\sin n\phi-qt)/\cos^n \theta-(t/n)\cos n\phi)\nonumber\\
&&-(1-\lambda)\eta^{1+\kappa}t|t|^{2\kappa}\sin\theta/\cos^{2n}\theta\label{zlc}
\end{eqnarray}
and plug it into (\ref{wsm4}).
For $\lambda=1$ we recover previous stage (\ref{reg}) while $\lambda=0$ is the final sphere of radius $R=\eta^{\kappa}|t|^{-1/n}$.
As shown in Appendix G, the mapping is smooth.


\section{Summary}

Using ruled surfaces, we have formulated the complete sphere eversion process in terms of direct analytic parameterization, given in (\ref{twormn}, (\ref{tworm2}), (\ref{reg}), (\ref{xyl}) and (\ref{zlc}) with mappings (\ref{wsm1}) and (\ref{wsm2}), keeping minimum of topological events, generalized also to Boy surface.
The process can be elegantly visualized in any computer modeling software with variable parameters $t$, $q$, $\xi$, $\eta$, $\alpha$, $\beta$, $\lambda$, $\omega$. We presented the suggested variation of these parameters in computer visualization in Table \ref{tta}. The process could be read from the bottom to the top for $t<0$ and back from the top to the bottom for $t>0$. We change linearly the only parameters differing in subsequent rows. Parameters $\beta$ and $\omega$ are arbitrary positive and $p=1-|qt|$. In the top stage we change $t$ linearly from $-1/Q$ to $+1/Q$ with $Q<1$, e.g. $Q=2/3$. Some stages can be combined if we keep the parameters within  smoothness ranges (we gave earlier them explicitly). We also believe the presented equations will help in 3D modeling of sphere eversion
in a simple, controlled and clear way.

\section*{Acknowledgement}
Research partially supported by NCN Grant UMO-2016/21/B/ST1/01489.

\begin{table}
\begin{tabular}{|c|c|c|c|c|c|c|}
\hline
&$|t|$&$q$&$\xi$&$\eta$&$\alpha$&$\lambda$\\
\hline
closed wormhole&$<1/Q$&$0$&$1$&$>0$&$>0$&$1$\\
\hline
unfolded wormhole&$1/Q$&$Q$&$1$&$>0$&$>0$&$1$\\
\hline
inverted wormhole&$1/Q$&$Q$&$0$&$1$&$0$&$1$\\
\hline
sphere&$1/Q$&$Q$&$0$&$1$&$0$&$0$\\
\hline
\end{tabular}
\caption{Suggested values of parameters used in visualization of complete sphere eversion, with $Q<1$ and $p=1-|qt|$, using eqs. (\ref{tworm2}),
(\ref{wsm1}), (\ref{wsm2}), (\ref{xyl}) and (\ref{zlc}).
}\label{tta}
\end{table}

\appendix
\renewcommand{\thesection}{\Alph{section}}

\section{Smoothness of the cylinder map}

We will find the condition when the parametric surfaces (\ref{tworm2}) and its special cases (\ref{twormn}) and (\ref{halfn}) are smooth. For convenience we will use complex representation $w=x+iy$ and $u=e^{i\phi}=\cos\phi+i\sin\phi$
with conjugation $\bar{w}=x-iy$.
Then (\ref{tworm2}) reads
\begin{equation}
\vec{r}=(w,z)=(tu+i(hu+p\bar{u}^{n-1}), h\sin n\phi-(t/n)\cos n\phi-qth).\label{twormnc}
\end{equation}
 We have tangent vectors
\begin{eqnarray}
&&\vec{r}_h=(iu,\sin n\phi-qt)\\
&&\vec{r}_\phi=(itu-hu+p(n-1)\bar{u}^{n-1},nh\cos n\phi+t\sin n\phi)\nonumber
\end{eqnarray}
and $\vec{n}=(n_x+in_y,n_z)=\vec{r}_h\times\vec{r}_\phi=(iz_h w_\phi-iz_\phi w_h,\mathrm{Im}(\bar{w}_h w_\phi))$ giving
\begin{eqnarray}
&&(n_x+in_y)\bar{u}=nh\cos n\phi+qt^2+i(p(n-1)\bar{u}^n-h)(\sin n\phi-qt),\nonumber\\
&&n_z=h-p(n-1)\cos n\phi.\label{nor}
\end{eqnarray}
Now $n_z=0$ when $h=p(n-1)\cos n\phi$, which plugged into $n_{x,y}$ gives
\begin{equation}
(n_x+in_y)\bar{u}=(n-1)p(n-(n-1)\sin^2n\phi-qt\sin n\phi)+qt^2
\end{equation}
The smallest value is at $\sin n\phi=\mathrm{sgn} t$ giving the condition (\ref{conpqt}).

\section{Surface $n=2$ equation}

We will derive surface and self-intersection equations by transforming the parametric form (\ref{halfn}) at $n=2$.
Comparing $xy$ with $z$ we get
\begin{equation}
2xyh=(1-h^2)z\label{zzz}
\end{equation}
From trigonometric unity $\cos^2\phi+\sin^2\phi=1$ we have
\begin{equation}
\frac{x^2}{(1-h)^2}+\frac{y^2}{(1+h)^2}=1
\end{equation}
which we transform into
\begin{equation}
(x^2+y^2)(1+h^2)+2h(x^2-y^2)=(1-h^2)^2.
\end{equation}
Replacing the right hand side with(\ref{zzz}) we get
\begin{equation}
(x^2+y^2)(1+h^2)+2h(x^2-y^2)=4h^2x^2y^2/z^2.\label{xyz}
\end{equation}
Both (\ref{zzz}) and (\ref{xyz}) are quadratic equations with respect to $h$. 

Eliminating $h$ leads to an equation for the surface. If $h$ is a common root of two quadratic equations $a_ih^2+b_ih+c_i=0$, $i=1,2$ then
$(c_1a_2-c_2a_1)^2=(a_1b_2-a_2b_1)(c_2b_1-c_1b_2)$, leading to the equation
\begin{eqnarray}
&&(z^3(x^2-y^2)-z^2xy(x^2+y^2)+4x^3y^3)(xy(x^2+y^2)+z(x^2-y^2))\nonumber\\
&&=(z^2(x^2+y^2)-2x^2y^2)^2
\end{eqnarray}
which reduces to the surface of degree 6
\begin{equation}
4xyz(x^2-y^2)+4x^2y^2(x^2+y^2-1)=4z^4+z^2(x^2+y^2)(x^2+y^2-4)
\end{equation}

For an arbitrary $t$ the surface (\ref{twormn}) at $n=2$ is still sextic (degree 6). To show it we combine (\ref{twormn})
to
\begin{eqnarray}
&&A=2xy=2th\cos 2\phi+(1+t^2-h^2)\sin 2\phi+2t,\nonumber\\
&&B=x^2-y^2=(t^2-1-h^2)\cos 2\phi-2th\sin 2\phi-2h,\nonumber\\
&&C=x^2+y^2=2h\cos 2\phi+2t\sin 2\phi+1+t^2+h^2.\label{xye}
\end{eqnarray}
Note that $C=\sqrt{A^2+B^2}$.
Eliminating linearly $\sin$ and $\cos$ we get
\begin{equation}
(t^2+h^2-1)^2=(h^2+t^2+1)C+2hB-2tA.\label{polA}
\end{equation}
On the other hand finding $\sin 2\phi$ and $\cos 2\phi$ from (\ref{xye}) and plugging them into $z=h\sin 2\phi-(t/2)\cos 2\phi$
we get
\begin{equation}
(2z+th)(h^2+t^2-1)=htC-tB-2hA.\label{polB}
\end{equation}
Now $h$ must be a common root of (\ref{polA}) and (\ref{polB}) so the greatest common divisor of polynomials in $h$ must be of degree at least $1$.
Applying Euclidean algorithm, dividing polynomials (\ref{polA}) (degree 4) by (\ref{polB}) (degree 3), then (\ref{polB}) by the remainder (degree 2), and finally both remainders (one of degree 1), the result is of degree $0$ so it must vanish.
The algorithm can be easily performed by Mathematica, and the result (without common factors) is the final surface equation

\begin{eqnarray}
&&4((x^2+y^2)(t^2+4)-8txy)(z^2(x^2+y^2)-4x^2y^2)\nonumber\\
&&+4z(x^2-y^2)(8tz^2-4xy(7t^2+4))\nonumber\\
&&+3t(t^2+4)(4z(x^4-y^4)+3t(x^2+y^2)^2)\nonumber\\
&&+16z^2(4z^2+((x^2+y^2)(7t^2-4)-4txy(t^2+3)))\nonumber\\
&&+4x^2y^2(3t^4-32t^2+16)+16xyt(2(x^2+y^2)-t^2)(1-2t^2)\nonumber\\
&&=8zt(x^2-y^2)(t^4-15t^2+8)\\
&&+4t^2((3t^4-5t^2+4)(x^2+y^2)-(t^2-1)(8z^2+t^4-t^2)).\nonumber
\end{eqnarray}

The surface remains sextic for (\ref{tworm2}) when we subtract $qht$ from $z$ ($p=1$ by scaling for simplicity). Then
we only need to replace $z$ by $z+qht$ in (\ref{polB}) to get a rather long result

\begin{eqnarray}
&&4 (2 q+1)^4 t^8-16 q (2 q+1)^2 x y t^7-4 (16
   (x^2+y^2+2) q^4+\nonumber\\
&&8(3 x^2+3
   y^2+8) q^3-(x^4+2(y^2-12)
   x^2+y^4-24 y^2-48) q^2+\nonumber\\
&&2(7 x^2+7
   y^2+8) q+3 x^2+3 y^2+2) t^6+8(2
   q(4 q^2+3 q+2) y x^3+\nonumber\\
&&(4
   q^2-1) z x^2+2 y (16 q^4+4
   (y^2+11) q^3+(3 y^2+40)
   q^2+\nonumber\\
&&(2 y^2+15) q+2) x+(1-4
   q^2) y^2 z) t^5+(-64
   (x^2+y^2-1) q^4\nonumber\\
&&-32(x^4+5(2
   y^2+1) x^2+y^4+5 y^2-4) q^3-4
   (x^6+(3 y^2-1) x^4\nonumber\\
&&+(3 y^4+78
   y^2+8) x^2+y^6-y^4+8 y^2-32 z^2-24)
   q^2+\nonumber\\
&&4(7 x^4+(14-22 y^2) x^2+7
   y^4+14 y^2+32 z^2+8) q+9 x^4+9 y^4+
\nonumber\\
&&20
   y^2+32 z^2+10 x^2(3 y^2+2)+4)
   t^4+4(4 q (2 q-1) y x^5+
\nonumber\\
&&(3-8 q^2+2
   q) z x^4+4 y(8 q^3+(4
   y^2-2) q^2
-(2 y^2+13) q-4)
   x^3
\nonumber\\
&&+2(16 q^3+36 q^2+44 q+15) z x^2+4
   y(8(y^2-1) q^3+2
   (y^4-y^2-6) q^2
\nonumber\\
&&-(y^4+13 y^2+4
   z^2+6) q-4 y^2-4 z^2-1) x
\nonumber\\
&&-(2 q+1) zy^2
   (16 q^2-4(y^2-7) q+3
   (y^2+10))) t^3-4(4
   q^2 x^6
\nonumber\\
&&+(4(3 y^2-5) q^2-4
   (8 y^2+7) q+4 y^2-z^2-9) x^4+4
   (12 q+7) y z x^3
\nonumber\\
&&+2(2 y^4-(z^2-7)
   y^2-14 z^2-4 q(4 y^4+5 y^2+3
   z^2-2)+\nonumber\\
&&q^2(6 y^4-20 y^2+8
   z^2+8)+2) x^2
-4 (12 q+7) y^3 z x+4
   q^2 y^6
\nonumber\\
&&+8 (2 q z+z)^2-y^4(20 q^2+28
   q+z^2+9)+4 y^2(4(z^2+1)
   q^2+
\nonumber\\
&&(4-6 z^2) q-7 z^2+1))
   t^2-16(4 q y x^5-(2 q+3) z x^4+
\nonumber\\
&&2 y (z^2-4
   y^2+q(4 y^2-2)-1) x^3-2 z
   (z^2-4 q-2) x^2+
\nonumber\\
&&2 y(2 q
   y^4+(z^2-2 q-1) y^2+6 (2 q+1)
   z^2) x+
\nonumber\\
&&y^2 z(3 y^2+2 z^2+2 q
   (y^2-4)-4)) t-16
   ((4 y^2-z^2) x^4+4 y z
   x^3\nonumber\\
&&+(4 y^4-2(z^2+2) y^2
+4
   z^2) x^2
-4 y^3 z x-z^2 (y^4-4 y^2+4
   z^2))=0.
\end{eqnarray}

\section{Topological events}

It is known that the sphere eversion must go through certain special points related to self-intersections, although just specifying these points is insufficient.
Nevertheless, to understand the complexity of eversion, one should be able to capture these points. The points represent equivalence classes for surfaces
modified smoothly in the arbitrarily small neighborhood of the point. Therefore the precise shape of the surface is irrelevant.

\begin{figure}
\includegraphics[scale=0.3]{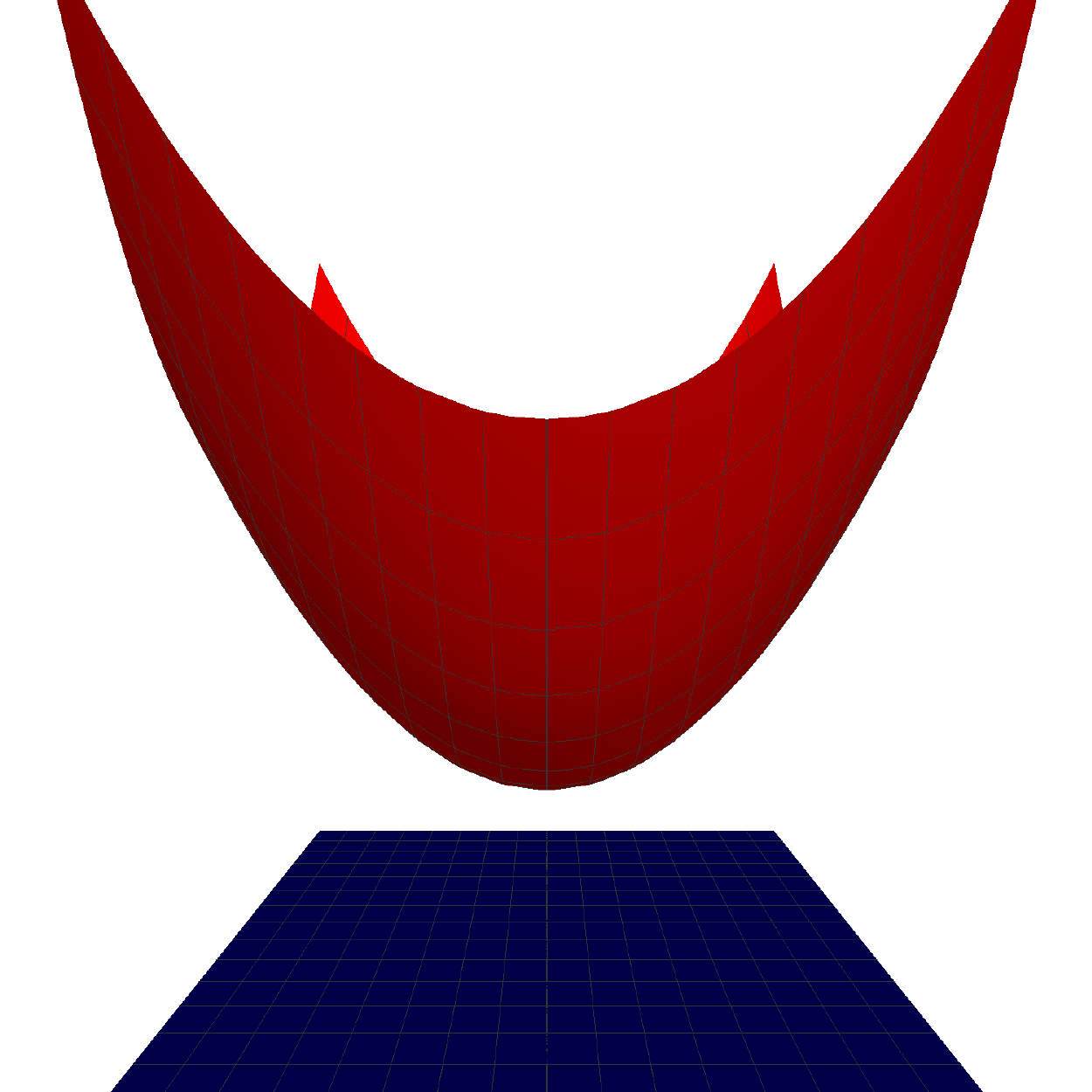}
\includegraphics[scale=0.3]{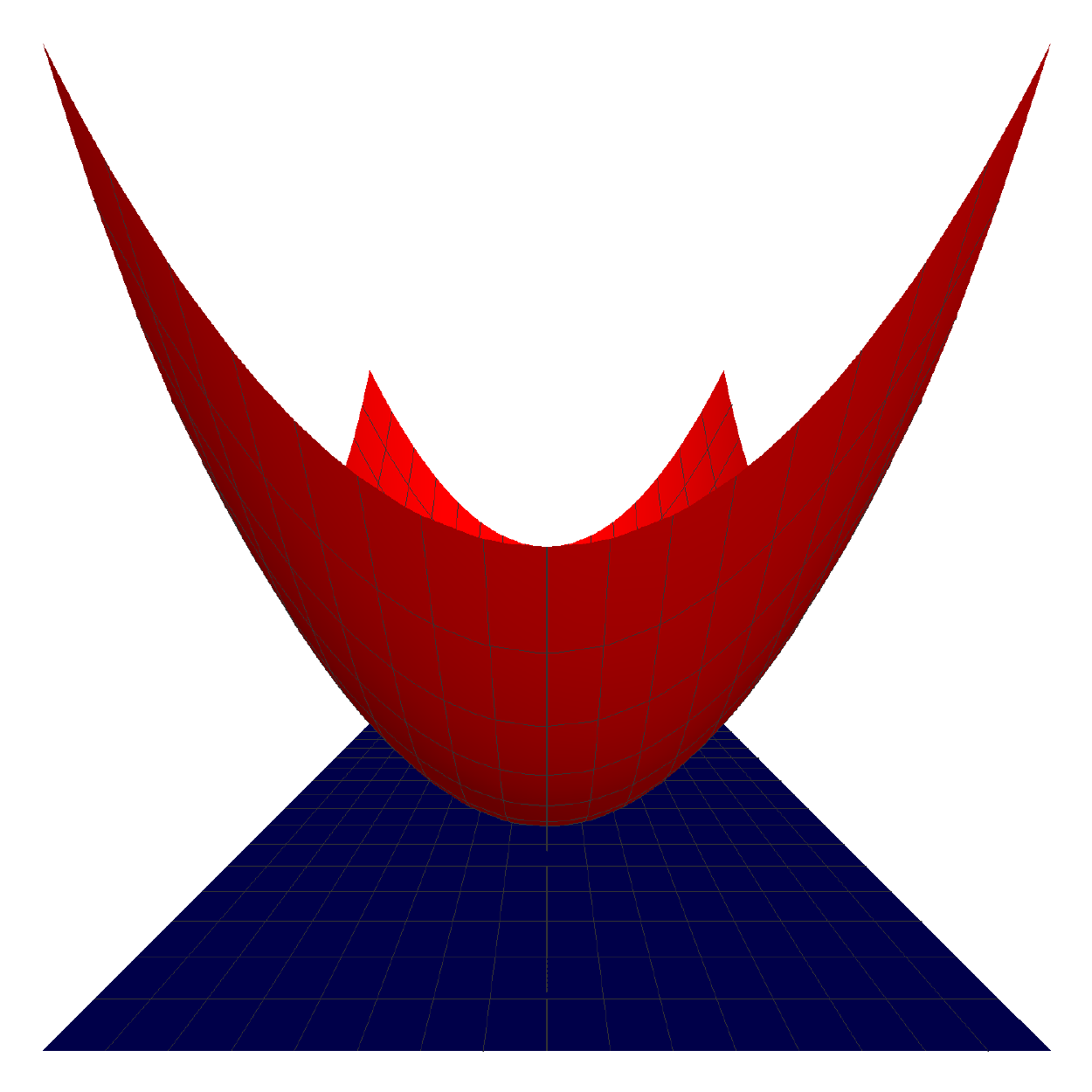}
\includegraphics[scale=0.3]{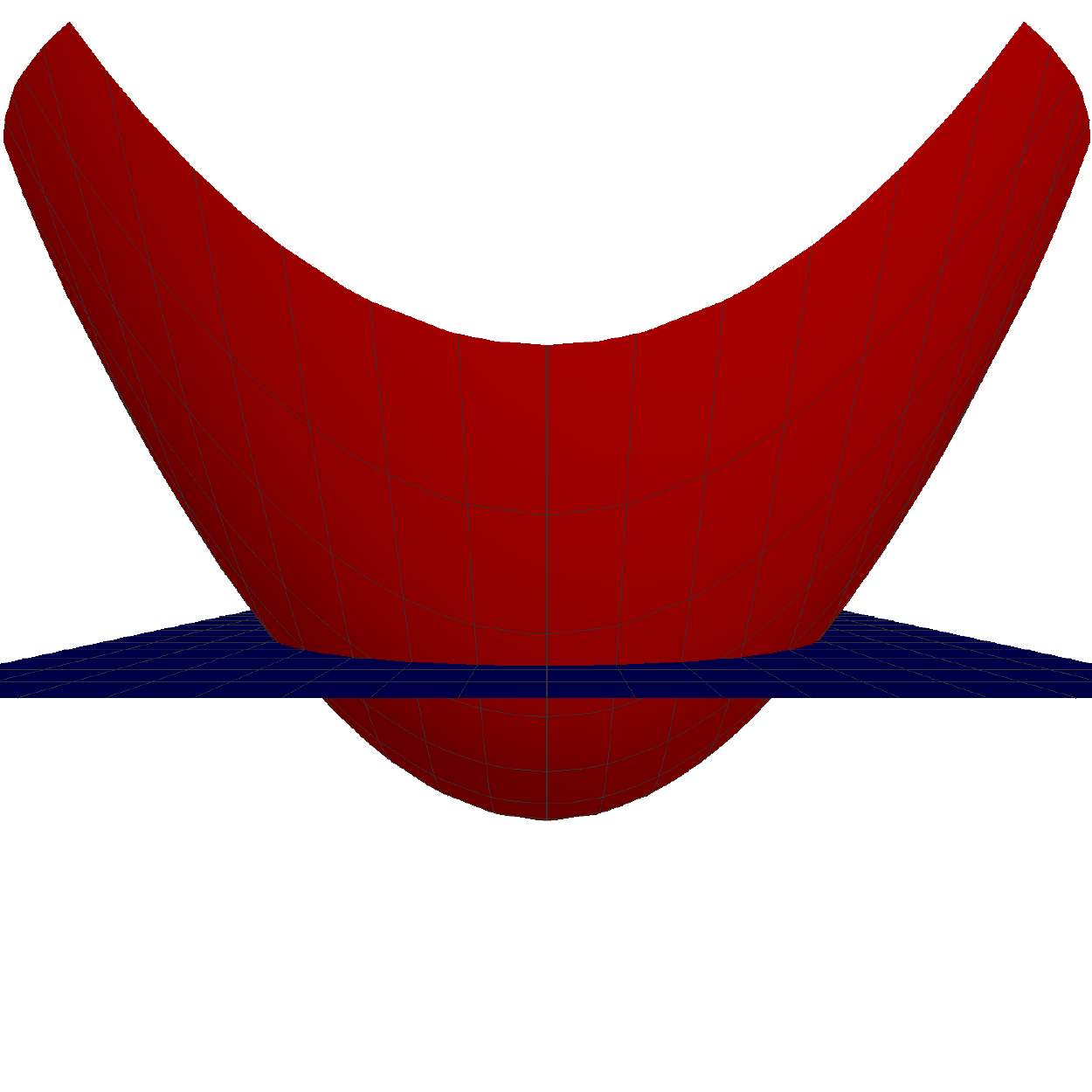}
\caption{Going through a $D_0$ ($D_2$) point in the middle, from the left to the right (viceversa)}
\label{d02}
\end{figure}

\begin{figure}
\includegraphics[scale=0.3]{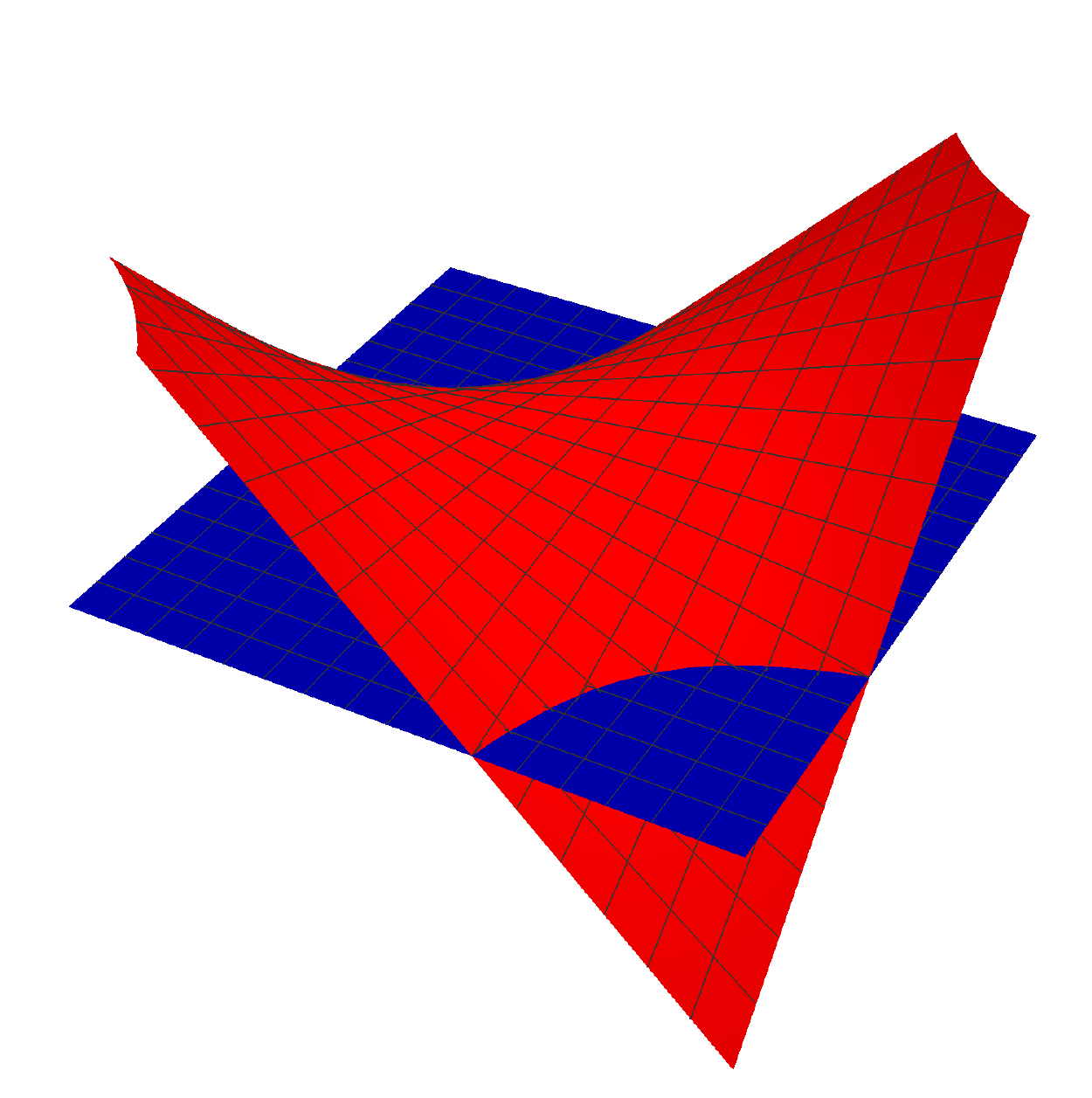}
\includegraphics[scale=0.3]{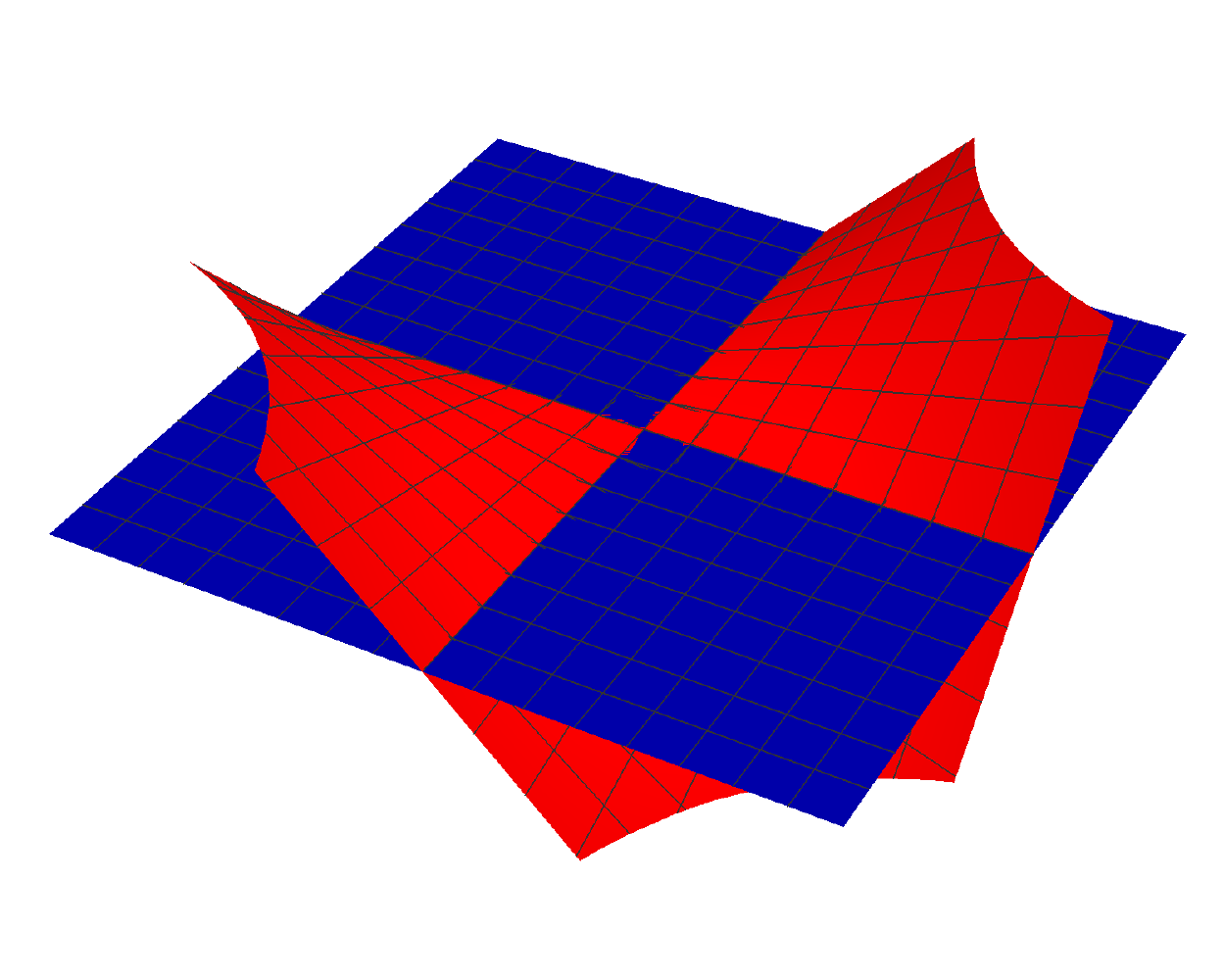}
\includegraphics[scale=0.3]{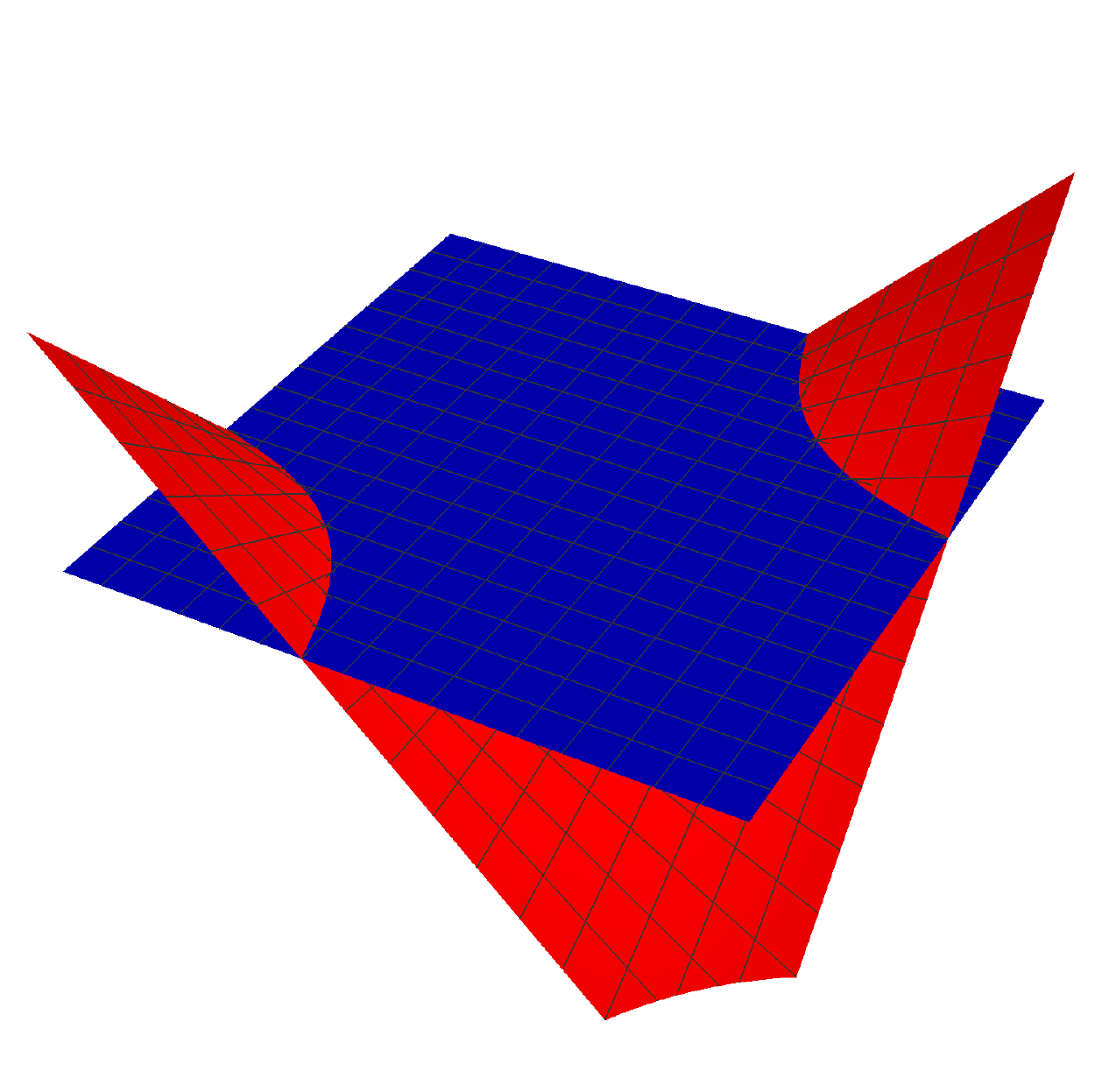}
\caption{Going through a $D_1$ point in the middle}
\label{d1}
\end{figure}

Following previous studies \cite{morin2,morin3,apery}, we start with description of the point $D_0$ and $D_2$ ($D$ stands for double), depicted in Fig. \ref{d02}. This is essentially the same point but traversed forwards and backwards, respectively, when changing $t$. Precisely, we have one surface $z=0$ and the second one moving as
$z=x^2+y^2-t$. For $t<0$ the surfaces do not intersect, for $t=0$ the touch in one point $(0,0,0)$ and for $t>0$ the intersect at the loop $z=0$, $x^2+y^2=t$. The point $D_0$ denotes moving from $t<0$ to $t>0$ and $D_2$ viceversa. Therefore $D_0$ ($D_2$) is a birth (death) of self-intersection loop. Another important point is $D_1$, the saddle, Fig. \ref{d1}. Here again one surface is fixed, $z=0$ while the second is moving, $z=xy-t$. The intersection equation $xy=t$ and $z=0$ gives two disjoint hyperbolas for $t\neq 0$ and two straight lines $x=0$ and $y=0$ at $t=0$. The point $D_1$ allows to switch between two intersection lines.

Other $D$-type events can occur in combination of $D_{0/2}$ and $D_1$ events, e.g. $D_{01}$ or $D_{21}$, between no intersections and two disconnected intersections. It is realized by a fixed surface $z=0$ and the moving one $z=x^2-y^2+t(1+x^2+y^2)(x^2+y^2)$ at $t=\pm 1$.

\begin{figure}
\includegraphics[scale=0.3]{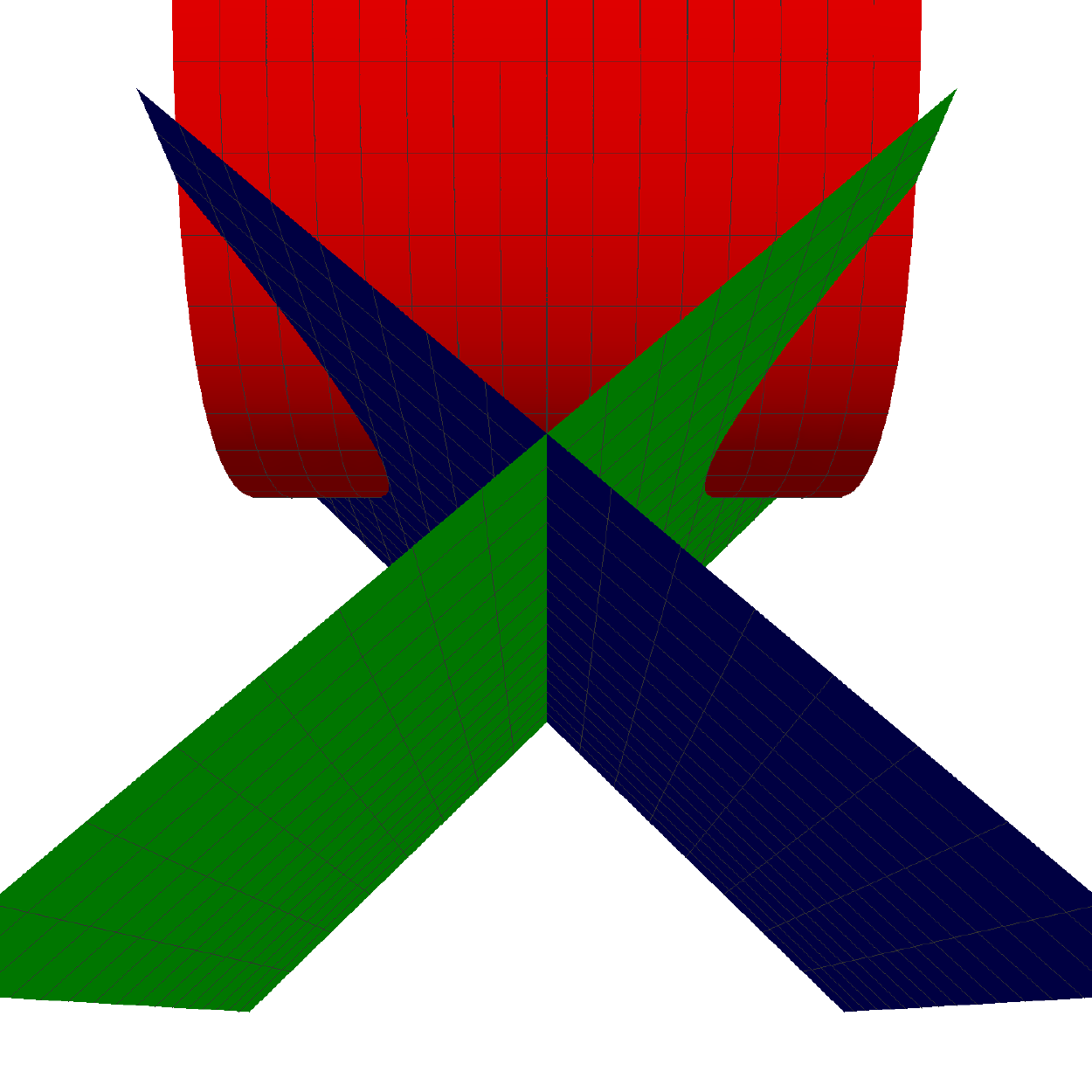}
\includegraphics[scale=0.3]{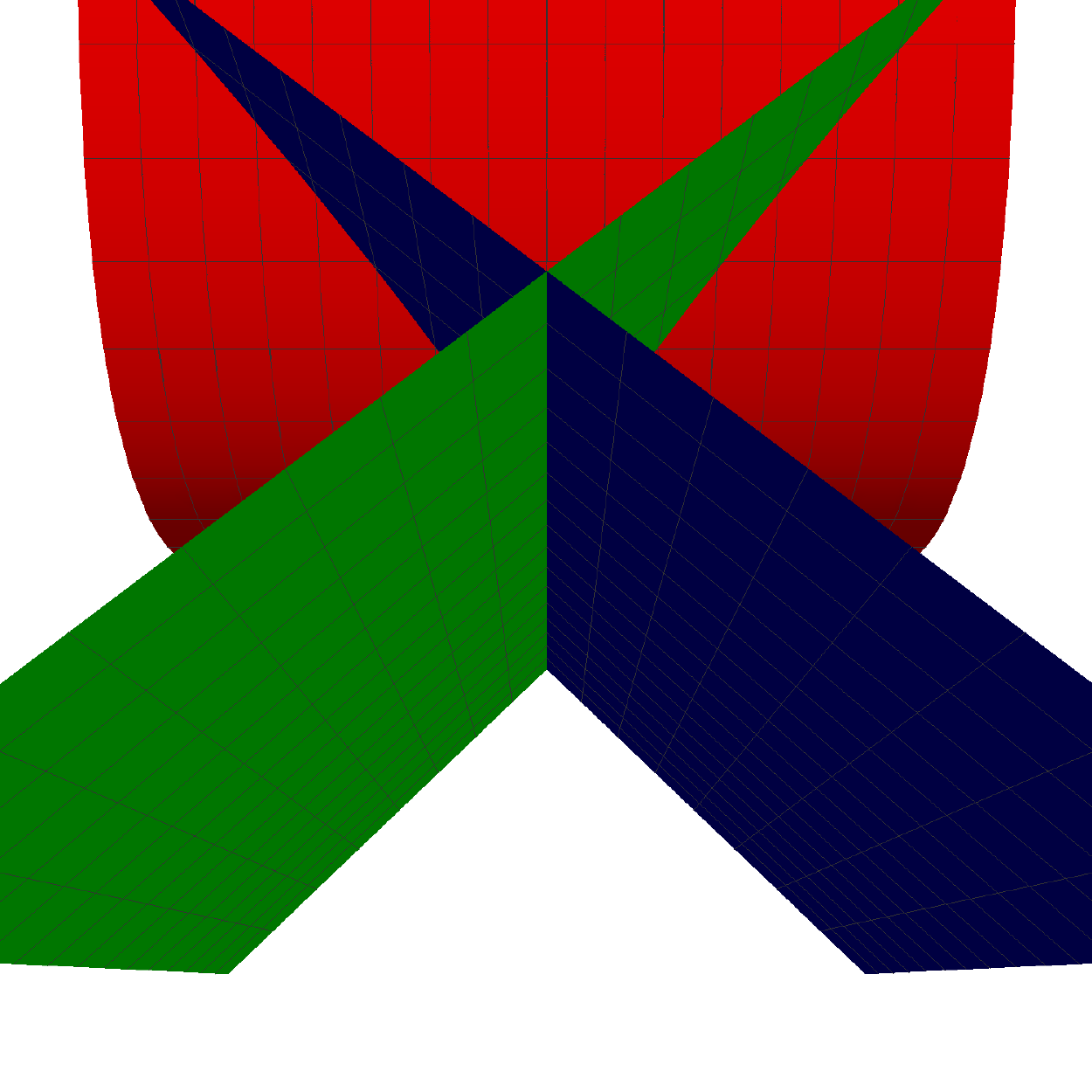}
\includegraphics[scale=0.3]{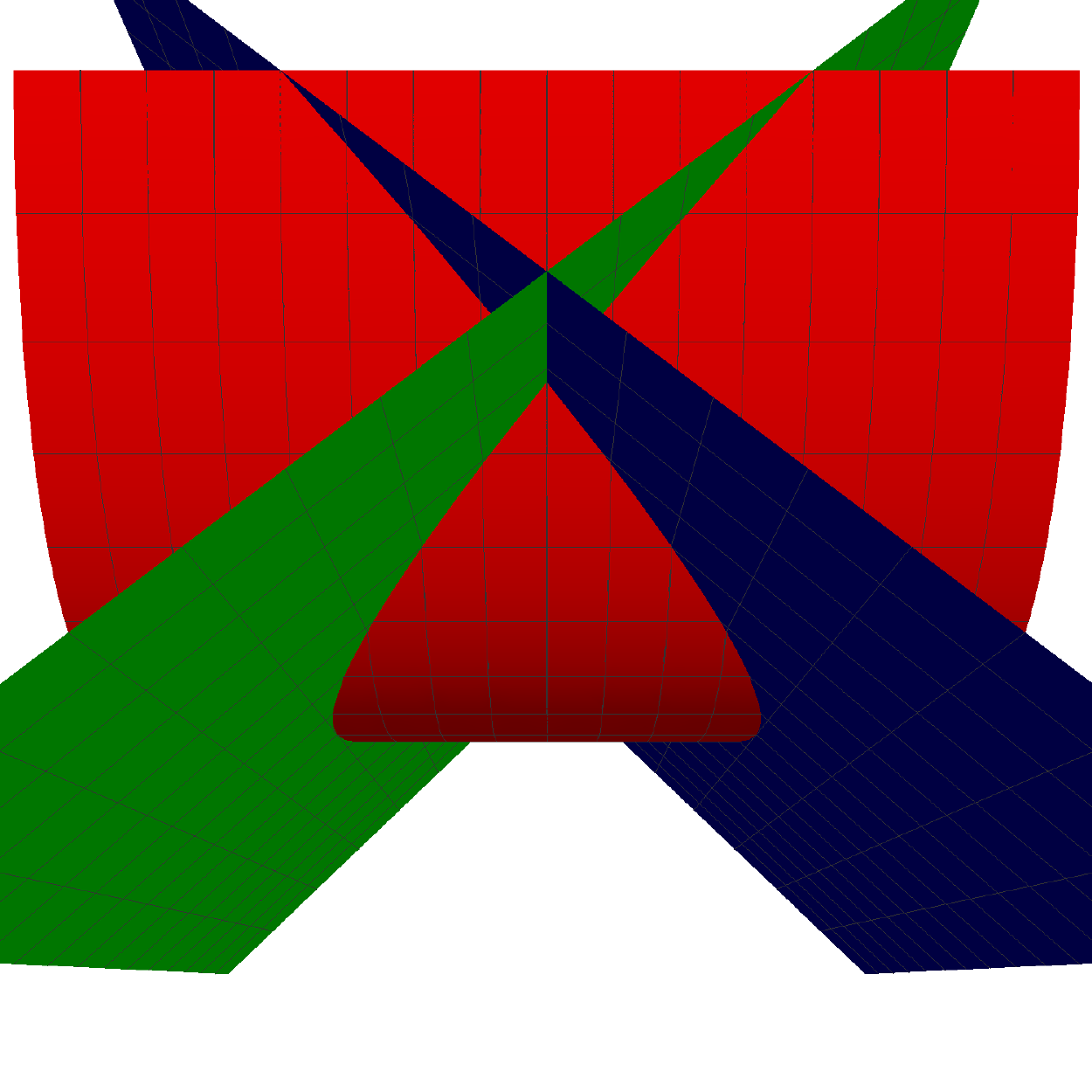}
\caption{Going through a $T_+$ ($T_-$) point in the middle, from the left to the right (viceversa)}
\label{tpm}
\end{figure}

The next important point is birth (death) of triple points $T_\pm$ (three smooth surfaces usually can intersect in one point). We can have two fixed surfaces $|x|=|y|$ and one moving $x=z^2-t$. For $t<0$ we have two disjoint moving intersection parabolas $\pm y=x=z^2-t$ (and fixed intersection line at $x=y=0$). At $t=0$ the parabolas touch each other and the line at $(0,0,0)$. At $t>0$ the parabolas and the line intersect at two triple points $x=0=0$, $z=\pm\sqrt{t}$. Therefore $T_+$ ($T_-$) describes movement from $t<0$ to $t>0$ (viceversa)
and essentially means birth (death) of a pair of triple points (they must always come in pairs), Fig. \ref{tpm}.

\begin{figure}
\includegraphics[scale=0.3]{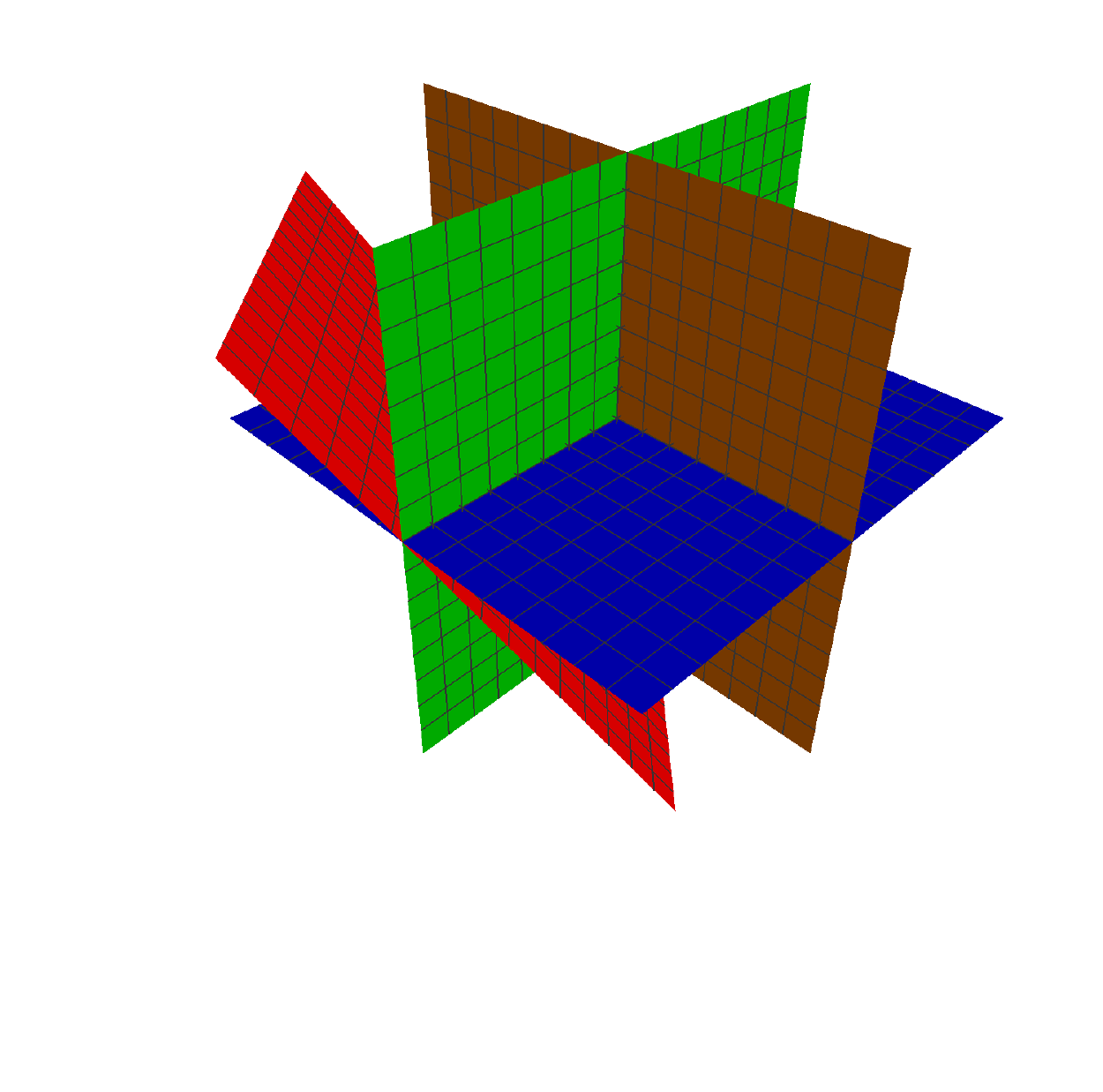}
\includegraphics[scale=0.3]{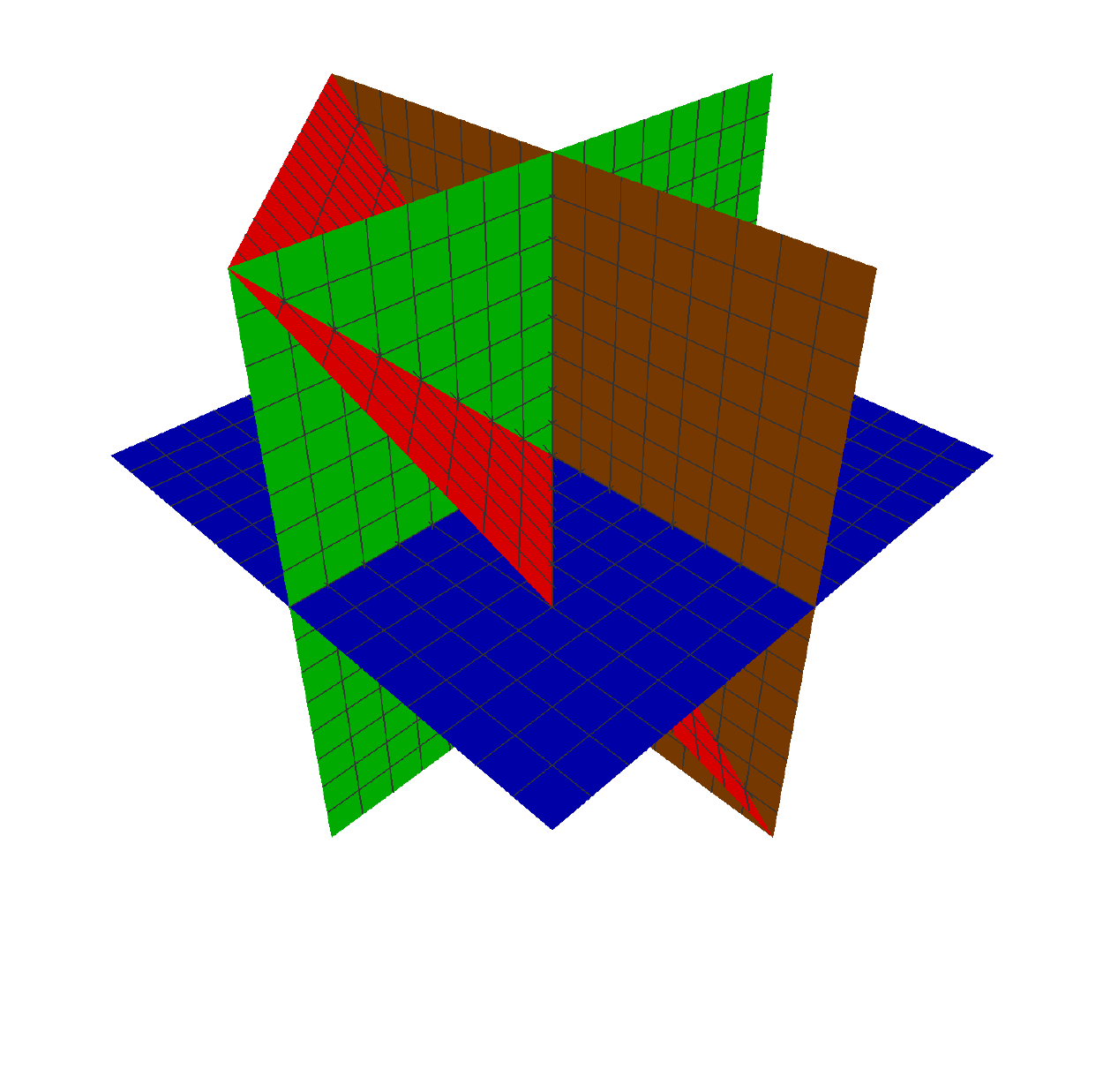}
\includegraphics[scale=0.3]{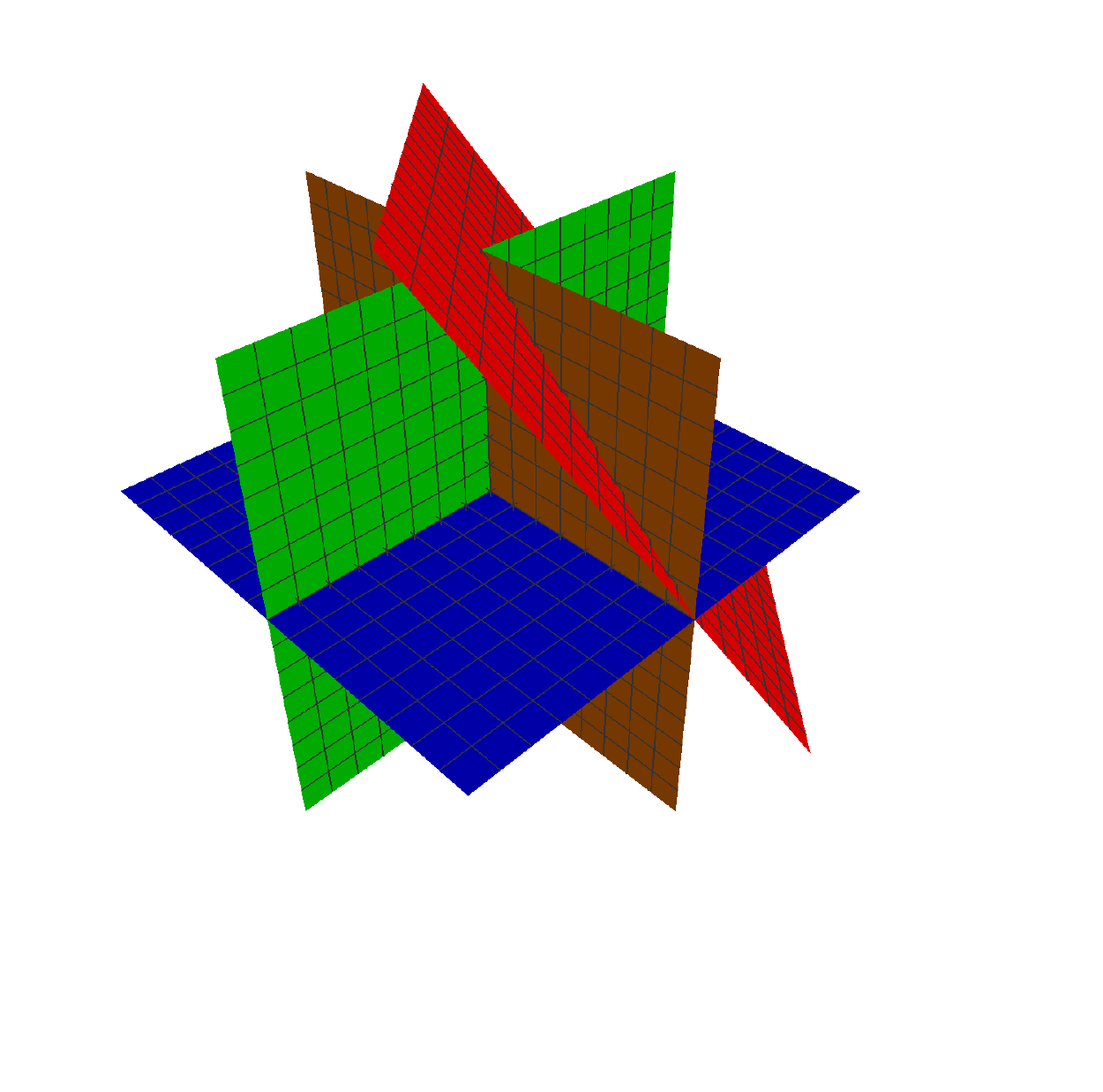}
\caption{Going through a $Q$ point in the middle}
\label{qq}
\end{figure}

The last critical point is the situation when \emph{four} surfaces meet at a single point. Although it is uncommon in a stationary
immersion, we need this point in the dynamics ($t$-dependence) of sphere eversion at the halfway moment. The point can be described e.g. taking three fixed surfaces, say $x=0$, $y=0$ and $z=0$ and one moving, e.g. $x+y+z=t$. Then only at $t=0$ all the four meet at $(0,0,0)$.
This point is called $Q$ (quadruple), Fig. \ref{qq}.

One can encode a sphere eversion by specifying subsequent topological events, occurring during the eversion process, although the encoding is not sufficient to describe the complete eversion. The most known encoding is as follows \cite{morin3, apery}. One starts with two $D_0$ (order irrelevant), then two $T_+$ (order irrelevant). The halfway moment corresponds to $Q$ and four $D_1$ (order irrelevant) and one additional $D_1$ which is rather independent and can be put anywhere between the first $D_0$ and the last $D_2$. The second half of the process is just reversing the former events. Namely, two $T_-$ and then two $D_2$. 

\section{Self-intersections in general}

General intersection lines can be derived as follows. Let two different pairs $(h_1,\phi_1)$ and $(h_2,\phi_2)$ map 
to the same point ($x,y,z$)
Then (\ref{tworm2}) using notation from Appendix A gives
\begin{eqnarray}
&&w=tu_1+ip\bar{u}^{n-1}_1+ih_1u_1=tu_2+ip\bar{u}^{n-1}_2+ih_2u_2\\
&&z=h_1\sin n\phi_1-(t/n)\cos n\phi_1-qth_1=\nonumber\\
&&h_2\sin n\phi_2-(t/n)\cos 2\phi_2-qth_2.\nonumber
\end{eqnarray}
Subtracting middle and right hand sides we get
\begin{eqnarray}
&&t(u_1-u_2)+ip(\bar{u}^{n-1}_1-\bar{u}^{n-1}_2)=\nonumber\\
&&i(h_1+h_2)(u_2-u_1)/2+i(h_2-h_1)(u_1+u_2)/2,\\
&&(h_1+h_2)(\sin n\phi_1-\sin n\phi_2)/2+(h_1-h_2)(\sin n\phi_1+\sin n\phi_2)/2\nonumber\\
&&=(t/n)(\cos n\phi_1-\cos n\phi_2)+qt(h_1-h_2).\nonumber
\end{eqnarray}
Adding middle and right hand sides we get
\begin{eqnarray}
&&2w=t(u_1+u_2)+ip(\bar{u}^{n-1}_1+\bar{u}^{n-1}_2)\nonumber\\
&&+i(h_1+h_2)(u_1+u_2)/2+i(h_1-h_2)(u_1-u_2)/2,\\
&&2z=(h_1+h_2)(\sin n\phi_1+\sin n\phi_2)/2+(h_1-h_2)(\sin n\phi_1-\sin n\phi_2)/2\nonumber\\
&&-(t/n)(\cos n\phi_1+\cos n\phi_2)-qt(h_1+h_2).\nonumber
\end{eqnarray}
Let us introduce $2h_\pm=h_1\pm h_2$ and $2\phi_\pm=\phi_1\pm \phi_2$
so that $h_1=h_++h_-$, $h_2=h_+-h_-$, $\phi_1=\phi_++\phi_-$, $\phi_2=\phi_+-\phi_-$.
Due to $2\pi$ periodicity of $\phi_{1,2}$ we can restrict $(\phi_+,\phi_-)\in[0,2\pi]\times[0,\pi[$.
With
\begin{eqnarray}
&&\sin m\phi_1\pm\sin m\phi_2=2\sin m\phi_\pm\cos m\phi_\mp,\\
&&\cos m\phi_1+\cos m\phi_2=2\cos m\phi_+\cos m\phi_-,\nonumber\\
&&\cos m\phi_1-\cos m\phi_2=-2\sin m\phi_+\sin m\phi_-,\nonumber
\end{eqnarray}
we get
\begin{eqnarray}
&&itu_+\sin\phi_-+p\bar{u}^{n-1}_+\sin(n-1)\phi_-=\nonumber\\
&&h_+u_+\sin\phi_--ih_-u_+\cos\phi_-,\label{hph}\\
&&h_+\sin n\phi_-\cos n\phi_++h_-\sin n\phi_+\cos n\phi_-=\nonumber\\
&&(-t/n)\sin n\phi_+\sin n\phi_-+qth_-\nonumber
\end{eqnarray}
and
\begin{eqnarray}
&&w=tu_+\cos\phi_-+ip\bar{u}_+^{n-1}\cos(n-1)\phi_-\nonumber\\
&&+ih_+u_+\cos\phi_--h_-u_+\sin\phi_-,\label{avg}\\
&&z=h_+\sin n\phi_+\cos n\phi_-+h_-\sin n\phi_-\cos n\phi_+\nonumber\\
&&-(t/n)\cos n\phi_+\cos n\phi_--qth_+.\nonumber
\end{eqnarray}

From the first equation of (\ref{hph}) we get
\begin{eqnarray}
&&h_+\sin\phi_-=p\cos n\phi_+\sin(n-1)\phi_-,\nonumber\\
&&h_-\cos\phi_-=p\sin n\phi_+\sin(n-1)\phi_--t\sin\phi_-.\label{hcs}
\end{eqnarray}
If $\phi_-=0$ then $h_-=0$ meaning that the points are the same (modulo $2\pi$ for $\phi$) so we will require that $\phi_-\neq 0$ so that $h_+=p\cos n\phi_+\sin(n-1)\phi_-/\sin\phi_-$. If $\phi_-=\pi/2$ then we consider two cases, even or odd $n$.
For odd $n$, $t=h_+=0$, obtaining the doubly covered Boy surface (with periodicity $\phi\to\phi+\pi$, $h\to -h$ and opposite orientation).
For even $n=2k$ we get $t=(-1)^{k-1}p\sin n\phi_+$, $h_+=(-1)^{k-1}p\cos n\phi_+$. If $p=0$ then $t=0$ which violates smoothness condition so $p\neq 0$. Then the third line of (\ref{hph}) reduces to $th_-(q+1/p)=0$. Since $q\geq 0$ and $p>0$ it follows $h_-=0$ or $t=0$. If $h_-=0$ then
$x=y=0$, $h=h_+$ and $z=h((1-1/n)\sin n\phi_+-qt)$ but only if $|t|\leq p$. 
 Plugging into (\ref{nor}) one can check directly that the normals are not parallel, except the case $|t|=p$ ($D_{0,2}$ point).
For $t=0$ we get arbitrary $h_-$ and $\sin 2k\phi_+=0$, $z=0$, $w=-h_-u_+$, meaning straight lines in $xy$ plane distributed at
equal angles (e.g. $x$ and $y$ axes for $n=2$) $\phi_1=\pi j/n+\pi/2$, 
$\phi_2=\pi j/n-\pi/2$, $h_+=(-1)^{j+k-1}p$. From now on we assume $\phi_-\neq 0,\pi/2$ so that (\ref{hcs}) plugged into the third line of (\ref{hph}) results in
\begin{eqnarray}
&&p\cos^2n\phi_+\sin n\phi_-\sin(n-1)\phi_-\cos\phi_-/\sin^2\phi_-\nonumber\\
&&+p\sin^2n\phi_+\cos n\phi_-\sin(n-1)\phi_-/\sin\phi_-=\nonumber\\
&&t\sin n\phi_+(\cos n\phi_--(1/n)\sin n\phi_-\cos\phi_-/\sin\phi_-)\nonumber\\
&&+qt(p\sin n\phi_+\sin(n-1)\phi_-/\sin\phi_--t)\label{eqqn}
\end{eqnarray}
At $t=0$ (halfway) we get 
\begin{equation}
\cos 2n\phi_+\sin(n-1)\phi_-=-\sin(n+1)\phi_-
\end{equation}
or $\sin(n-1)\phi_-=0$. In particular for $n=2$ we have
$2\cos^2\phi_-=\sin^22\phi_+$ so then, unless $\sin 2\phi_+=0$,
\begin{eqnarray}
&&x=\mp \sqrt{2}p\cos\phi_+\cos 2\phi_+,\nonumber\\
&&y=\pm \sqrt{2}p\sin\phi_+\cos 2\phi_+,\nonumber\\
&&z=p\sin 4\phi_+/2,
\end{eqnarray}
whose projection in $xy$ plane is a quadrifolium. The intersection set includes also $x$ and $y$ axes in the case $\sin 2\phi_+=0$.
The points of intersection of the quadrifolium with $x$ and $y$ axes  ($\varphi=0,\pm \pi/2,\pi$ and $r=\sqrt{2}$, $z=0$) define $D_1$ points.

For $n=3$ we get
\begin{equation}
2\cos 2\phi_-=-\cos 6\phi_+\label{eqpm}
\end{equation}
and
\begin{equation}
w=-pe^{-2i\phi_+}\sin 6\phi_+, z=-(p/4)\sin 12\phi_+.
\end{equation}
The projection onto $xy$ plane is trifolium $(-\cos 2\phi_+\sin 6\phi_+,\sin 2\phi_+\sin 6\phi_+)$.

For $n=2$ (\ref{eqqn}) is equivalent to
\begin{equation}
p(\sin^22\phi_+-2\cos^2\phi_-)=t\sin 2\phi_+\sin^2\phi_--qt(p\sin 2\phi_+-t)\label{eqq}
\end{equation}
and finally
\begin{eqnarray}
&&x\cos\phi_-=2p\sin\phi_+\cos^2\phi_-+\cos\phi_+(t-p\sin 2\phi_+),\nonumber\\
&&y\cos\phi_-=2p\cos\phi_+\cos^2\phi_-+\sin\phi_+(t-p\sin 2\phi_+),\nonumber\\
&&2z=p\sin 4\phi_+-t\cos 2\phi_+(1+2\sin^2\phi_-)-2qtp\cos 2\phi_+,\label{sss}
\end{eqnarray}

Now we can keep $\phi_+$ as an independent variable
\begin{equation}
\cos^2\phi_-=\frac{(\sin 2\phi_++qt)(p\sin 2\phi_+-t)}{2p-t\sin 2\phi_+}.
\end{equation}
Since $\cos^2$ lies in $[0,1]$ we have valid disjoint intervals $\sin 2\phi_+\in[-1,-qt]$ and $[t,p]$ for $t>0$
while $[-qt,1]$ and $[-p,t]$ for $t<0$.
Plugging them in to
(\ref{sss}) we get
\begin{eqnarray}
&&x=
\cos\phi_-\frac{\cos\phi_+(t\sin 2\phi_+-2p\cos 2\phi_+)+2pqt\sin\phi_+}{\sin 2\phi_++qt},\nonumber\\
&&y=\cos\phi_-\frac{\sin\phi_+(t\sin 2\phi_++2p\cos 2\phi_+)+2pqt\cos\phi_+}{\sin 2\phi_++qt}.
\label{ss}
\end{eqnarray}
 For $|t|>p$  one of the intervals disappears because of $D_{0,2}$ events and the second one disappears if $|qt|\geq 1$. 
There also no other $D$ events.
Let us show this first for $\sin 2\phi_+\in[-1,-qt]$ and $t>0$ (without loss of generality).
The function $w(S=-\sin 2\phi_+)=x^2+y^2$ is equal
\begin{equation}
\frac{(pS+t)((t^2+4pqt^2)S^2+4p(p+tS)(1-S^2)+4p^2q^2t^2)}{(S-qt)(2p+St)}
\end{equation}
while $-(S-qt)^2(2p+St)^2w'(S)$ is equal
\begin{eqnarray}
&&\tilde{q}S^2t^4+8p^4(2S^3+\tilde{q}(1+\tilde{q}^2)-3\tilde{q}S^2)+\nonumber\\
&&pSt^3(3S^3+4\tilde{q}+S(2+4\tilde{q}^2)-6\tilde{q}S^2)+\nonumber\\
&&2p^2t^2(4S^5+4S(1+3\tilde{	q}^2)+S^3(6+4\tilde{q}^2)+2\tilde{q}(1-\tilde{q}^2)-13\tilde{q}S^2-8\tilde{q}S^4)\nonumber\\
&&+4p^3t(2+2\tilde{q}^2+7S^4+S^2(1+7\tilde{q}^2)-14\tilde{q}S^3)\label{der}
\end{eqnarray}
with $\tilde q=qt$. We will show that, under certain condition, the above expression is positive for $t>0$ and $S\in[0,1]$ (without loss of generality). 
Then $w=x^2+y^2$ is monotonic in $S$ and there is no loop causing $D_{0/2}$.

In each term of (\ref{der}) the negative parts are located at the end.
We will show that the other parts overrule the negativity, using Cauchy inequality $a+b\geq 2\sqrt{ab}$ and our assumptions, including $|qt|<1$. 
In the first line $\tilde{q}(1+\tilde{q}^2)\geq 2\tilde q^2$ and
$2\tilde{q}^2+2\tilde{S}^3\geq 4\tilde{q}S^{3/2}\geq 4\tilde{q}S^2$.
In the second line $3S(2+4\tilde{q}^2)\geq 4\sqrt{6} \tilde{q}S$.
In the third line $4S+4S^5\geq 8S^3$ and $8S^3+12S\tilde{q}^2\geq 8\sqrt{6}\tilde{q}S^2$
which is $\geq 13\tilde{q}S^2$ and $6+4\tilde{q}^2\geq 4\sqrt{6}\tilde{q}\geq 8\tilde{q}S$.
In the last line $7S^4+S^2\geq 2\sqrt{7}S^3$ and $2\sqrt{7}S^3+7\tilde{q}^2S^2\geq
2^{3/2}7^{3/4}\tilde{q}S^{5/2}\geq 10\tilde{q}S^3$ and $2+2\tilde{q}^2\geq 4\tilde{q}S^3$.

On the other hand for the case $p=1$, $q=0$, $t\in]0,1]$ and $\sin 2\phi_+\in[t,1]$ we can introduce $g(k=\tan\phi_+)=y/x$, with $k>0$, given by
\begin{equation}
2k\frac{tk+1-k^2}{tk-1+k^2}=-2k+4t-4t\frac{tk-1}{tk-1+k^2}.
\end{equation}
Now $g'(k)<0$ because
$
(tk-1+k^2)^2>2t((tk-1)(t+2k)-t(tk-1+k^2))
$ equivalent to $(tk-1+k^2)^2>2tk(tk-2)$ or $k^4+1+2tk^3+2tk>t^2k^2+2k^2$.
It follows from $k^4+1\geq 2k^2$ and $2tk^3+2tk\geq 2tk^2\geq t^2k^2$. So now $y/x$ is monotonic in $\tan\phi_+$. The loops can be only connected in 3 cases: $\sin 2\phi_+=t$ for $p=1$ and $q=0$,
$\sin 2\phi_+=\pm 1$, and $\sin 2\phi_+=-qt$. In the first case from (\ref{eqq}) we get $(t^2-2)\cos^2\phi_-=0$ so $\cos\phi_-=x=y=0$, the case already considered (this is the main central loop with $D_{0/2}$ at $(0,0,0)$). In the last case $x$ or $y$ diverges. In the second case $|x|=|y|$. We will find all solutions of $x=y$ (without loss of generality).
From (\ref{ss}) we get
\begin{equation}
(\cos\phi_+-\sin\phi_+)(t\sin 2\phi_+-2pqt)=2p(\cos\phi_++\sin\phi_+)\cos 2\phi_+
\end{equation}
Defining $\tilde{\phi}=\phi_++\pi/4$ we have $\cos\phi_+-\sin\phi_+=\sqrt{2}\cos\tilde\phi$,
$\cos\phi_++\sin\phi_+=\sqrt{2}\sin\tilde{\phi}$, $\sin 2\phi_+=-\cos 2\tilde{\phi}$, $\cos 2\phi_+=\sin 2\tilde{\phi}$
so that
\begin{equation}
\cos\tilde{\phi}(-t\cos 2\tilde\phi-2pqt)=2p\sin\tilde{\phi}\sin 2\tilde{\phi}.
\end{equation}
One solution is $\cos\tilde{\phi}=0$ corresponding to the second case. The other solution would require
\begin{equation}
-2pqt=2p(1-\cos 2\tilde{\phi})+t\cos 2\tilde{\phi}.
\end{equation}
For $t>0$ we have $\cos 2\tilde\phi>0$ so equality is impossible.  No other $D$ point occurs.

\section{Properties of the open Boy sufrace}

We will find the equation of the (Boy) surface (\ref{halfn}) with $n=3$ corresponding to (\ref{twormnc}) with $t=0$.
We have $wu^2=ihu^3+i$ so $\mathrm{Re}wu^2=-z$ and $\mathrm{Re}(w\bar{u}-i\bar{u}^3)=0$. This gives
\begin{equation}
wu^4+2zu^2+\bar{w}=0,\;iu^6+\bar{w}u^4+wu^2-i=0
\end{equation}
Eliminating $u^6$ and $u^4$ in the second equation by means of the first one we get
\begin{equation}
u^2(w^3-i|w|^2+4iz^2-2z|w|^2)=iw^2+\bar{w}(|w|^2-2iz)
\end{equation}
Since $|u^2|=1$, taking square of moduli we can write
\begin{eqnarray}
&&|w|^6+|w|^4+16z^4+4z^2|w|^4+\nonumber\\
&&(8z^2-2|w|^2)\mathrm{Im}w^3-4z|w|^2\mathrm{Re}w^3-8|w|^2z^2=\nonumber\\
&&|w|^4+4z^2|w|^2+|w|^6-4z\mathrm{Re}w^3-2|w|^2\mathrm{Im}w^3
\end{eqnarray}
reducing by division by $4z$ to
\begin{equation}
4z^3+z|w|^4-3z|w|^2+2z\mathrm{Im}w^3+(1-|w|^2)\mathrm{Re}w^3=0
\end{equation}
or explicitly
\begin{equation}
4z^3+z(x^2+y^2)(x^2+y^2-3)+2z(3x^2y-y^3)+(1-x^2-y^2)(x^3-3xy^2)=0.
\end{equation}

\section{Determining $T_\pm$ events and triple points}

The value of $t$ corresponding  to events $T_\pm$ for (\ref{twormn}) and $n=2$ can be found with help of symmetry of the surface.
Namely, the $T$ points must be located at $x=y$ or $x=-y$ and $z=0$ and be common for 3 different points $(h,\phi)$.
One family of such points is given by $h=0$ and $\phi=\pm\pi/4,\pm 3\pi/4$. We now only need to find a \emph{different} family mapping to the same $\vec{r}$. The condition $z=0$ means $2h=t\cot 2\phi$. Let us first consider $x=y$ which means
\begin{equation}
(t-1)(\cos\phi-\sin\phi)=h(\cos\phi+\sin\phi)
\end{equation}
equivalent to
\begin{equation}
(t-1)\cos\psi=h\sin\psi=-t\cos\psi\frac{\sin^2\psi}{\cos 2\psi}
\end{equation}
with $\psi=\phi+\pi/4$.
Discarding the case $\cos\psi=0$ we obtain
\begin{equation}
t-1=-t\frac{\sin^2\psi}{\cos 2\psi}
\end{equation}
and 
\begin{equation}
t^{-1}=\frac{\cos^2\psi}{\cos 2\psi}.\label{ttt}
\end{equation}
On the other hand, the point must match the other family with $x=y=\pm (t+1)/\sqrt{2}$ so
\begin{equation}
x+y=\pm(t+1)\sqrt{2}=(t+1)(\sin\phi+\cos\phi)+h(\cos\phi-\sin\phi)
\end{equation}
equivalent to
\begin{equation}
\pm (t+1)=(t+1)\sin\psi+h\cos\psi.
\end{equation}
Substituting $h$, we get
\begin{equation}
\pm (t+1)=\sin\psi(t+1-t\cos^2\psi/\cos 2\psi)=t\sin\psi.
\end{equation}
Squaring yields
\begin{equation}
(t+1)^2=t^2\sin^2\psi.
\end{equation}
We combine it with (\ref{ttt}) to get
\begin{equation}
(t+1)^2=t^2(1-1/(2-t))=t\frac{1-t}{2-t}
\end{equation}
and finally $(2-t)(t+1)^2=t^2(1-t)$ equivalent to $t^2-3t-2=0$. From the two solutions $2t=3\pm\sqrt{17}$ we have to exclude
$t>0$ because then $\sin^2\psi=(1+t^{-1})^2>1$ and we are left with $t=(3-\sqrt{17})/2$.
For $x=-y$ the analysis is analogous and yields opposite $t=(\sqrt{17}-3)/2$.

We shall also show how to find triple points (intersections of three parts of the surface) between $T_\pm$ events.
From (\ref{twormn}) and $n=2$ we infer
\begin{equation}
\cos\phi=\frac{xt-y(1-h)}{h^2+t^2-1},\sin\phi=\frac{yt-x(1+h)}{h^2+t^2-1}
\end{equation}
and
\begin{eqnarray}
&&(z+qth)(h^2+t^2-1)^2=2h(xt-y(1-h))(yt-x(1+h))\\
&&-(t/2)(xt-y(1-h))^2-(yt-x(1+h))^2.\nonumber
\end{eqnarray}
It is clear that three different values of $h$ must correspond to a triple point
From trigonometric unity $\sin^2\phi+\cos^2\phi=1$ and opening brackets we get equations
\begin{equation}
(h^2+t^2-1)^2=(x^2+y^2)(t^2+1+h^2)-2(x^2-y^2)h-4txy
\end{equation}
and
\begin{eqnarray}
&&z(h^2+t^2-1)^2=\nonumber\\
&&2h(xy(t^2+1-h^2)-x^2t(1+h)-y^2t(1-h))\nonumber\\
&&-(t/2)((x^2-y^2)(t^2-1-h^2)-2h(x^2+y^2)+4xyht).
\end{eqnarray}
Let us parametrize $x=\sqrt{w}\cos(\psi/2)$, $y=\sqrt{w}\sin(\psi/2)$ with $w\geq 0$.
Then, the above equations reduce to
\begin{equation}
(h^2+t^2-1)^2=w(t^2+1+h^2-2h\cos\psi-2t\sin\psi)\label{www}
\end{equation}
and
\begin{eqnarray}
&&z(h^2+t^2-1)^2=w[h(t^2+1-h^2)\sin\psi-2ht(1+h\cos\psi)\nonumber\\
&&-(t/2)((t^2-1-h^2)\cos\psi-2h+2ht\sin\psi)]=\nonumber\\
&&w[h(1-h^2)\sin\psi-(3h^2+t^2-1)(t/2)\cos\psi-ht].
\end{eqnarray}
By transformation $\psi\to \pi-\psi$, $h,z\to -h,-z$ and $\psi\to \psi+2\pi$ we get 4 triple points from a single one.
Comparing left hand sides we get
\begin{eqnarray}
&&z(t^2+1+h^2-2h\cos\psi-2t\sin\psi)=\nonumber\\
&&h(1-h^2)\sin\psi-(3h^2+t^2-1)(t/2)\cos\psi-ht.\label{zz1}
\end{eqnarray}
This cubic equation in $h$ must have all 3 roots at triple point (no nontrivial quadruple point occurs).
Since quartic equation (\ref{www}) must have the same roots (\ref{zz1}) is its divisor,
there must exist numbers $A,B$ such that
\begin{eqnarray}
&&(h^2+t^2-1)^2-w(t^2+1+h^2-2h\cos\psi-2t\sin\psi)\equiv\\
&&(Ah+B)[z(t^2+1+h^2-2h\cos\psi-2t\sin\psi)\nonumber\\
&&-h(1-h^2)\sin\psi+(3h^2+t^2-1)(t/2)\cos\psi+ht]
\end{eqnarray}
for all $h$. Therefore all coefficients of both sides must be equal.
It leads to 5 equations
\begin{equation}
A\sin\psi=1,
\end{equation}
\begin{equation}
z'+(3t/2)A\cos\psi+B\sin\psi=0,\label{zpp}
\end{equation}
\begin{equation}
2(t^2-1)-w'=-2z'\cos\psi+(3t/2) B\cos\psi+A(t-\sin\psi),\label{wpp}
\end{equation}
\begin{equation}
2w'\cos\psi=z'(t^2+1-2t\sin\psi)+B(t-\sin\psi)+At\cos\psi(t^2-1)/2,\label{h11}
\end{equation}
\begin{equation}
(t^2-1)^2-w'(t^2+1-2t\sin\psi)=Bt\cos\psi (t^2-1)/2,\label{h00}
\end{equation}
where we denoted $w'=w+Bz$ and $z'=Az$.
Replacing $2(t^2-1)$ by $2(t^2-1)A\sin\psi$ and $(t^2-1)^2$ by $(t^2-1)^2A\sin\psi$, equations (\ref{zpp}), (\ref{wpp}),(\ref{h11})
and (\ref{h00}) become a set of 4 homogeneous linear equations of 4 variables $A,B,w',z'$ with nonzero solution ($A\neq 0$)
so the determinant of the following matrix
\begin{equation}
\begin{pmatrix}
(3t/2)\cos\psi&\sin\psi&0&1\\
(2t^2-1)\sin\psi-t&-(3t/2)\cos\psi&-1&2\cos\psi\\
\cos\psi(t-t^3)/2&\sin\psi-t&2\cos\psi&2t\sin\psi-1-t^2\\
(t^2-1)^2\sin\psi&\cos\psi(t-t^3)/2&2t\sin\psi-1-t^2&0
\end{pmatrix}
\end{equation}
equal
\begin{equation}
t^6+12t^4+4t^2-8st^3(4+t^2)+2s^2t^2(2+7t^2)+8s^3t(1+t^2)-s^4(4+3t^2+4t^4)\label{dets}
\end{equation}
vanishes (here $s=\sin\psi$).
The determinant is of degree $4$ in $s$ and $6$ in $t$. For $t=0$ it reduces to $-s^4$. For $t\neq 0$ it has two roots, $s_+>0$
and $s_-<0$. To show this, note that it goes to $-\infty$ at $s\to\pm\infty$ and it is positive at $s=0$. They are the only roots, which is due to the fact that there is only one local maximum of (\ref{dets}) as a function of $s$ because the first derivative
\begin{equation}
-8t^3(4+t^2)+4st^2(2+7t^2)+24s^2t(1+t^2)-4s^3(4+3t^2+4t^4)\label{dfir}
\end{equation}
has a single root. From Cardano method the cubic equation $as^3+bs^2+cs+d=0$
it is true if its discriminant $\Delta=18abcd-4b^3d+b^2c^2-4ac^3-27a^2d^2$ is negative.
Here $\Delta$ is equal
\begin{equation}
-2^{10}t^6(t^2-1)^2(3388+4796 t^2+4735 t^4+2084t^6+432t^8).
\end{equation}
For $s=\pm 1$ (\ref{dets}) reduces to $(t\mp 2)(t\mp 1)^3(t^2\mp 3t -2)$ allowing to recover once again $T_\pm$ events.
From continuity, the triple point between $Q$ and $T$ events will correspond to only one of the roots, with $ts<0$. We discard the other root, which is the easiest to prove by continuity. That other root will remain between $0$ and $1$ for $t\in ]0,1[$
(the case of negative $t$ is analogous). For $t=1$ we have $s=1$ but then $A=1$, $w'=0$ and $z'=-B$ giving $w=-B^2$ so $B=0=w$.
In addition, triple points can disappear only in pairs. We have already found $T_\pm$ events. The only left  possibility of pair disappearance is at $z$ axis, namely $(0,0,z)$ corresponding to $w=0$. For $w=0$ and $t\neq 0$ our equations reduce to $h=\pm \sqrt{1-t^2}$, $\cot\phi=(h-1)/t$ and 
\begin{equation}
z=\frac{4h(h-1)t-t((h-1)^2-t^2)}{2((h-1)^2+t^2)}=t+2\frac{h^2-1}{h^2-2h+1+t^2}.
\end{equation}
The last fraction is different for different signs of $h$ except $h=\pm 1$ but then $t=0$.
For the same $h$ there are only two inequivalent $\phi$ and $\pi+\phi$ giving the same $z$, so there is no new $T_\pm$ event nor other triple points by continuity.

\section{Smoothness of closed wormhole}

To show that the mappings (\ref{wsm1}) and (\ref{wsm2}) are smooth we will prove that their Jacobi matrices are nondegenerate,
the Jacobian (determinant of derivatives) is nonzero.
The Jacobi matrix of (\ref{wsm1}) in $xyz$ basis is 
\begin{equation}
\begin{pmatrix}
\displaystyle\frac{\xi+\eta(x^2(1-2\kappa)+y^2)}{(\xi+\eta|w|^2)^{1+\kappa}}&
\displaystyle-\frac{2\kappa\eta xy}{(\xi+\eta|w|^2)^{1+\kappa}}&0\\
\displaystyle-\frac{2\kappa\eta xy}{(\xi+\eta|w|^2)^{1+\kappa}}&\displaystyle\frac{\xi+\eta(y^2(1-2\kappa)+x^2)}{(\xi+\eta|w|^2)^{1+\kappa}}&0\\
\ast&\ast&(\xi+\eta|w|^2)^{-1}
\end{pmatrix}
\end{equation}
whose determinant is equal
\begin{equation}
\frac{\xi^2+\xi\eta|w|^2(2-2\kappa)+\eta^2|w|^4(1-2\kappa)}{(\xi+\eta|w|^2)^{3+2\kappa}}
\end{equation}
 while the Jacobi matrix  of (\ref{wsm2}) in $x'y'z'$ basis is
\begin{equation}
e^{\gamma z'}
\begin{pmatrix}\displaystyle
\frac{\alpha+\beta(y^{\prime 2}-x^{\prime 2})}{(\alpha+\beta(x^{\prime 2}+y^{\prime 2}))^2}&\displaystyle-\frac{2\beta x'y'}{(\alpha+\beta(x^{\prime 2}+y^{\prime 2}))^2}&\displaystyle\frac{\gamma x'}{\alpha+\beta(x^{\prime 2}+y^{\prime 2})}\\
\displaystyle-\frac{2\beta x'y'}{(\alpha+\beta(x^{\prime 2}+y^{\prime 2}))^2}&\displaystyle\frac{\alpha+\beta(x^{\prime 2}-y^{\prime 2})}{(\alpha+\beta(x^{\prime 2}+y^{\prime 2}))^2}&\displaystyle\frac{\gamma y'}{\alpha+\beta(x^{\prime 2}+y^{\prime 2})}\\
\displaystyle-\frac{ 4\alpha\beta x'/\gamma}{(\alpha+\beta(x^{\prime 2}+y^{\prime 2}))^2}&\displaystyle-\frac{ 4\alpha\beta y'/\gamma}{(\alpha+\beta(x^{\prime 2}+y^{\prime 2}))^2}&
\displaystyle\frac{\alpha-\beta(x^{\prime 2}+y^{\prime 2})}{\alpha+\beta(x^{\prime 2}+y^{\prime 2})}
\end{pmatrix}
\end{equation}
with determinant $e^{\gamma z'}(\alpha+\beta(x^{\prime 2}+y^{\prime 2}))^{-2}$.

We will show that the surface  is smooth at $h\to\infty$  parameterizing
 $h=\omega\sin\theta/\cos^n\theta$ for $\theta\in[-\pi/2,\pi/2]$.
Similarly to the previous Appendices, we consider a general combination
 (\ref{tworm2}).
Then, with  $C=\cos^2\theta$, $W=uC^{1/2}$, $Z=\sin\theta$ ($|W|^2+Z^2=1$),
\begin{eqnarray}
&&R=C^n|w|^2=\\
&&Z^2\omega^2+C^n(t^2+p^2)+2tp C^{n/2}\mathrm{Im}W^n+ 2p\omega Z\mathrm{Re}W^n\nonumber
\end{eqnarray}
or
\begin{equation}
R=C^n|w|^2=\lambda^2\omega^2Z^2+t^2(\lambda^2 C^n+2(\lambda-\lambda^2)C^{n/2}+(1-\lambda)^2)
\end{equation}
in the case of (\ref{xyl}) and (\ref{xxyy})
and $R'=C^n(\xi+\eta|w|^2)=C^n\xi+\eta R$.
Now
\begin{eqnarray}
&&w''=wC^{(n+1)/2}\frac{R^{\prime \kappa}}{C\alpha R^{\prime 2\kappa}+\beta R}e^{\gamma z'},\\
&&z''=\frac{C\alpha R^{\prime 2\kappa}-\beta R}{C\alpha R^{\prime 2\kappa}+\beta R}\frac{e^{\gamma z'}}{\gamma}-\gamma^{-1}
\frac{\alpha-\beta}{\alpha+\beta}.\nonumber
\end{eqnarray}
The front factor is
\begin{equation}
wC^{(n+1)/2}=tWC^{n/2}+ip\bar{W}^{n-1}C+i\omega ZW
\end{equation}
and
\begin{equation}
z'R'=\omega Z\mathrm{Im}W^n-(C^{n/2}t/n)\mathrm{Re}W^n-qt\omega ZC^{n/2}.
\end{equation}
or
\begin{equation}
wC^{(n+1)/2}=tW(\lambda C^{n/2}+(1-\lambda))+i\lambda\omega ZW
\end{equation}
and
\begin{equation}
z'R'=\lambda(\omega Z\mathrm{Im}W^n-tC^{n/2}(\mathrm{Re}W^n/n+q\omega Z))-(1-\lambda)\eta^{1+\kappa}t|t|^{2\kappa}Z.
\end{equation}
in the case of (\ref{xxyy}) and (\ref{zlc}).
Note also the saddle-like shape of the $\theta=\pm \pi/2$ point (in opposite $z$ directions) disappearing at $|qt|=1$.
The surface is there $C^\infty$ at even $n$ but only $C^1$ at odd $n$ because of $C^{n/2}$ factors which
disappear only at $t=0$. Nevertheless, we can enforce $C^\infty$ smoothness replacing every $C^{n/2}$ by $(C^n+\epsilon)^{1/2}$ for sufficiently small positive $\epsilon$.

We will show that (\ref{xxyy}) is a growing function of $C=\cos^2\theta\in[0,1]$. It will suffice to show that the inverse $n$th power of (\ref{xxyy}) is decreasing, namely
\begin{equation}
(t^2(\lambda C^{n/2}+(1-\lambda))^2+\lambda^2\omega^2(1-C))C^{-n}.
\end{equation}
Its derivative with respect to $C$ is
\begin{equation}
t^2(-n\lambda(1-\lambda)C^{-n/2-1}-n(1-\lambda)^2C^{-n-1}+\lambda^2\omega^2((n-1)C^{-n}-nC^{-n-1}).
\end{equation}
Multiplied by $C^{n+1}$, it transforms into
\begin{equation}
-nt^2(1-\lambda)^2-n\lambda^2\omega^2+(n-1)\lambda^2\omega^2)C-nt^2\lambda(1-\lambda)C^{n/2}.
\end{equation}
Without the last negative term, the remaining linear function is negative at the endpoints $0,1$  and so in all $[0,1]$.

Now, if 
$
w=(a(h)+ib(h))u
$
and $|w|^2=a^2+b^2$ is a growing/decreasing function of $h$  at $a^2+b^2>0$ and $z(h,\phi)$ is smooth then the surface as a function of $(h,\phi)$ is smooth except possibly $h=0$. Tangent vectors are
\begin{equation}
w_h=(a'+ib')u,
\end{equation}
(here $a'=da/dh$, $b'=db/dh$) while
\begin{equation}
w_\phi=(ia-b)u.
\end{equation}
Then $n_z=\mathrm{Im}w^\ast_h w_\phi=(a^2+b^2)'/2$. Moreover, at $h=0$ smoothness is guaranteed by nonvanishing
$\lambda qt\omega+(1-\lambda)\eta^{1+\kappa}|t|^{2\kappa}$.

\end{document}